\documentclass[a4paper,reqno, oneside]{amsart}
\usepackage[english]{babel}
\usepackage{amsmath, amssymb, amsthm, amscd}
\usepackage{enumerate}
\usepackage{palatino}
\usepackage{mathpazo}
\usepackage{paralist}
\usepackage[a4paper]{geometry}

\usepackage[all]{xy}

\usepackage{subcaption} 
\usepackage{caption} 
\usepackage{floatrow}

\usepackage{appendix}

\usepackage{tikz} 
\usepackage{tikz-cd}

\usepackage{dashbox}

\makeatletter
\newcommand{\ostar}{\mathbin{\mathpalette\make@circled\ast}}
\newcommand{\make@circled}[2]{%
  \ooalign{$\m@th#1\smallbigcirc{#1}$\cr\hidewidth$\m@th#1#2$\hidewidth\cr}%
}
\newcommand{\smallbigcirc}[1]{%
  \vcenter{\hbox{\scalebox{0.77778}{$\m@th#1\bigcirc$}}}%
}
\makeatother

\renewcommand{\setminus}{\smallsetminus}

\let\oldtocsection=\tocsection

\let\oldtocsubsection=\tocsubsection

\let\oldtocsubsubsection=\tocsubsubsection

\renewcommand{\tocsection}[2]{\hspace{0em}\oldtocsection{#1}{#2}}
\renewcommand{\tocsubsection}[2]{\hspace{1em}\oldtocsubsection{#1}{#2}}
\renewcommand{\tocsubsubsection}[2]{\hspace{2em}\oldtocsubsubsection{#1}{#2}}

\theoremstyle{plain}
\newtheorem{theorem}{Theorem}[section]
\newtheorem{prop}[theorem]{Proposition}
\newtheorem{lemma}[theorem]{Lemma}
\newtheorem{cor}[theorem]{Corollary}

\newtheorem*{TheoremA}{Theorem A}
\newtheorem*{TheoremB}{Theorem B}
\newtheorem*{TheoremC}{Theorem C}
\newtheorem*{TheoremD}{Theorem D}
\newtheorem*{TheoremE}{Theorem E}

\theoremstyle{definition}
\newtheorem{definition}[theorem]{Definition}

\theoremstyle{remark}
\newtheorem{remark}[theorem]{Remark}
\newtheorem{example}[theorem]{Example}

\numberwithin{equation}{section}

\begin{document}
\setlength{\parindent}{0.cm}

\title[Extensions of the loop product and coproduct]{Extensions of the loop product and coproduct, the space of antipodal paths and resonances of closed geodesics}

\author{Maximilian Stegemeyer}
\address{Mathematisches Institut, Universit\"at Freiburg, Ernst-Zermelo-Straße 1, 79104 Freiburg, Germany}
\email{maximilian.stegemeyer@math.uni-freiburg.de}
\date{\today}


\begin{abstract}
    We study the space of paths in a closed manifold $M$ with endpoints determined by an involution $f\colon M\to M$.
    If the involution is fixed point free and if $M$ is $2$-connected then this path space is the universal covering space of the component of non-contractible loops of the free loop space of $M/\mathbb{Z}_2$.
    On the homology of said path space we study string topology operations which extend the Chas-Sullivan loop product and the Goresky-Hingston loop coproduct, respectively.
    We study the case of antipodal involution on the sphere in detail and use Morse-Bott theoretic methods to give a complete computation of the extended loop product and the extended coproduct on even-dimensional spheres.
    These results are then applied to prove a resonance theorem for closed geodesics on real projective space.
\end{abstract}
\maketitle

\tableofcontents

\section{Introduction}

Studying the topology of the free loop space of a closed manifold is of great importance for many problems in geometry and topology.
For the classical question of understanding the closed geodesics in a closed Riemannian manifold $M$ one finds that the closed geodesics are precisely the critical points of the energy functional on the free loop space $\Lambda M$.
Morse-theoretic methods have been developed to study how many closed geodesics there are and what properties they have.

Another area where the free loop space plays a central role is string topology. String topology studies algebraic structures on the homology of the free loop space of a closed manifold.
The most prominent operation is the \textit{Chas-Sullivan product}, or \textit{loop product} which is of the form
$$   \wedge_{\mathrm{CS}} \colon \mathrm{H}_{i}(\Lambda M)\otimes \mathrm{H}_j(\Lambda M)\to \mathrm{H}_{i+j-n}(\Lambda M)    $$
where $n  = \mathrm{dim}(M)$ .
This was first studied in \cite{chas:1999}.
An operation which is in some sense dual to the loop product is the \textit{Goresky-Hingston product}. 
This is a product of the form
$$   \ostar_{\mathrm{GH}} \colon \mathrm{H}^i(\Lambda M,M)\otimes \mathrm{H}^j(\Lambda M,M)\to \mathrm{H}^{i+j+n-1}(\Lambda M,M) .  $$
Here, $M$ is embedded into $\Lambda M$ as the constant loops.
This product was introduced in \cite{goresky:2009} based on ideas by Sullivan in \cite{sullivan:2004}.
In recent years the dual coproduct to the Goresky-Hingston product has received a lot of attention.
This \textit{string topology coproduct}, or \textit{loop coproduct} is a map 
$$  \vee\colon \mathrm{H}_i(\Lambda M,M) \to \mathrm{H}_{i+1-n}(\Lambda M\times \Lambda M,\Lambda M\times M\cup M\times \Lambda M) .   $$
It turns out that the coproduct is not homotopy invariant, see \cite{naef2024string}.

A connection between string topology and the theory of the closed geodesics was first made by Goresky and Hingston in \cite{goresky:2009}. One way of relating the loop products to closed geodesics is to study \textit{critical values} of homology and cohomology classes.
For a homology class $X\in\mathrm{H}_{\bullet}(\Lambda M)$ or a cohomology class $\varphi\in\mathrm{H}^{\bullet}(\Lambda M,M)$ one can define numbers $\mathrm{cr}(X)$ and $\mathrm{cr}(\varphi)$ with respect to the energy functional on $\Lambda M$ and these numbers are critical values of the energy functional.
The critical values  behave nicely with respect to the products, we have
$$    \mathrm{cr}(X\wedge_{\mathrm{CS}} Y) \leq \mathrm{cr}(X) + \mathrm{cr}(Y) \quad \text{and}\quad    \mathrm{cr}(\varphi\ostar_{\mathrm{GH}} \psi)\geq \mathrm{cr}(\varphi)  + \mathrm{cr}(\psi)$$
for $X,Y\in\mathrm{H}_{\bullet}(\Lambda M)$ and $\varphi,\psi\in\mathrm{H}^{\bullet}(\Lambda M,M)$.
Hingston and Rademacher use these inequalities as well as the explicit computation of the loop product and the Goresky-Hingston product for spheres to show a \textit{resonance theorem} for closed geodesics on spheres.
This theorem, \cite[Theorem 1.1]{hingston:2013}, states that for a given Riemannian metric on $\mathbb{S}^n$ there are constants $\alpha,\beta > 0$ such that for all $X\in\mathrm{H}_{\bullet}(\Lambda \mathbb{S}^n)$ we have
$$  | \mathrm{deg}(X) - \alpha \, \mathrm{cr}(X) | \leq \beta .   $$
This result can be related to the index growth of closed geodesics, see e.g. \cite[Theorem 1.3]{hingston:2013} and \cite[Section 2]{hingston:2013} for a discussion of the resonance theorem.
It is a natural question whether a similar resonance theorem also hold for other spaces.

Not only can string topology give insights into closed geodesics, but one can also use knowledge about closed geodesics to infer properties of string topology operations as is done in the seminal paper by Goresky and Hingston \cite{goresky:2009}.
In \cite{kupper2022string} and \cite{stegemeyer2024string} the Chas-Sullivan product of symmetric spaces and the Goresky-Hingston product on complex and quaternionic projective spaces are computed by using Morse-Bott theoretic methods for the energy functional of the symmetric metrics on these spaces.
This is based on results by Ziller in \cite{ziller:1977} who shows that the energy functional of compact symmetric spaces is a \textit{perfect Morse-Bott function}.

In this paper we define and study \textit{extensions} of the loop product and the loop coproduct and we investigate the relation of these operations to closed geodesics in the case of the sphere, respectively of real projective space.
Eventually, we will show a resonance theorem for closed geodesics on real projective space.

Let $M$ be a closed oriented manifold with a smooth involution $f\colon M\to M$.
We consider the path space 
$$  P_f M  = \{\gamma \colon [0,1]\to M\,|\,   \gamma(1) = f(\gamma(0)), \,\,\, \gamma \text{ absolutely continuous} \}  .  $$
If $M$ is $2$-connected and $f$ is fixed point free, then $P_f M$ is the universal covering space of the component of non-contractible loops of $\Lambda (M/\mathbb{Z}_2)$, where the $\mathbb{Z}_2$ action is induced by $f$.
We will define a pairing
\begin{equation}\label{eq_pairing_introduction}
   \wedge_f \colon   \mathrm{H}_i(P_f M)\otimes \mathrm{H}_j(P_f M)\to \mathrm{H}_{i+j-n}(\Lambda M) .  
\end{equation}
The idea behind the pairing $\wedge_f$ in equation \eqref{eq_pairing_introduction} is the following naive observation.
If we have two paths $\gamma,\sigma\in P_f M$ such that $\gamma(1) = \sigma(0)$, then the concatenation $\eta = \gamma \star \sigma$ is a closed loop.
Moreover, we shall define a left- and right-module structure of $\mathrm{H}_{\bullet}(P_f M)$ over the Chas-Sullivan ring $(\mathrm{H}_{\bullet}(\Lambda M),\wedge_{\mathrm{CS}})$.
Together with the loop product these four structures yield the following extension of the loop product.
\begin{TheoremA}[Theorem \ref{theorem_extended_cs_product}]
    Let $M$ be a closed manifold with an involution $f\colon M\to M$.
    Take homology with coefficients in a commutative unital ring $R$.
    The pairing $\wedge_f$ together with the left- and right-module structure over the Chas-Sullivan ring and together with the loop product make the direct sum $\mathrm{H}_{\bullet}(\Lambda M)\oplus \mathrm{H}_{\bullet}(P_f M)$ into a unital associative algebra.
    Moreover, the canonical inclusion $\mathrm{H}_{\bullet}(\Lambda M)\hookrightarrow \mathrm{H}_{\bullet}(\Lambda M)\oplus \mathrm{H}_{\bullet}(P_f M)$ embeds $\mathrm{H}_{\bullet}(\Lambda M)$ as a subalgebra.
\end{TheoremA}

It is convenient to consider the degree shifted direct sum 
$  \widehat{\mathbb{H}}_{\bullet}(M) := \mathrm{H}_{\bullet + n}(\Lambda M)\oplus \mathrm{H}_{\bullet+n}(P_f M)    $
as the extended loop product is then of degree zero.
For homotopic involutions $f,g\colon M\to M$ there are homotopy equivalences $P_f M\simeq P_g M$. 
We show that the induced isomorphisms in homology intertwine the pairings $\wedge_f\colon \mathrm{H}_{\bullet}(P_f M)^{\otimes 2}\to \mathrm{H}_{\bullet}(\Lambda M)$ and $\wedge_g\colon \mathrm{H}_{\bullet}(P_g M)^{\otimes 2}\to \mathrm{H}_{\bullet}(\Lambda M)$ as well as the respective module structures over the Chas-Sullivan ring.
In particular in the case that $f\simeq \mathrm{id}_M \colon M\to M$ the module structures as well as the pairing are determined by the Chas-Sullivan product of $M$.
We use this homotopy invariance to compute the extended loop algebra for odd-dimensional spheres with the antipodal map.
\begin{TheoremB}[Theorem \ref{theorem_extended_product_circle} and Corollary \ref{cor_extendedn_product_sphere}]
    Consider the antipodal map $a\colon \mathbb{S}^n\to \mathbb{S}^n$ and take homology with integer coefficients.
    \begin{enumerate}
        \item For $n = 1$ we have an isomorphism of algebras
        $$   \widehat{\mathbb{H}}_{\bullet}(\mathbb{S}^1) \cong \Lambda_{\mathbb{Z}}(a) \otimes \mathbb{Z}[t,t^{-1}] \quad \text{with} \,\,\, |a|= -1, |t| =0.    $$
        \item For $n = 4k-1$, $k\geq 1$ we have an isomorphism of algebras
        $$   \widehat{\mathbb{H}}_{\bullet}(\mathbb{S}^n) \cong \frac{\mathbb{Z}[A,U,{E},B,V]}{\big(  A^2, B^2, BV - AU, V^2 - U^2, {E}B = A, {E}V = U, {E}^2 -1    \big)}     $$
    where $|A|= |B| = -n$, $|{E}| = 0$ and $|U| = |V| = n-1$. 
    \end{enumerate}
\end{TheoremB}

In order to understand the extended loop product for even-dimensional spheres we then study the homology of the space $P_a \mathbb{S}^n$ using Morse-Bott theory.
Note that the paths in $P_a\mathbb{S}^n$ are paths $\gamma\colon [0,1]\to\mathbb{S}^n$ such that $\gamma(1) = -\gamma(0)$.
We therefore refer to $P_a \mathbb{S}^n$ as the \textit{space of antipodal paths} on the sphere.
For a Riemannian metric $g$ on $\mathbb{S}^n$ we consider the \textit{energy functional} $E\colon P_a \mathbb{S}^n\to [0,\infty)$ given by
$$   E(\gamma)  = \int_0^1 g_{\gamma(t)}(\Dot{\gamma}(t)  ,\Dot{\gamma}(t)) \,\mathrm{d}t \quad \text{for}\,\,\,\gamma\in P_a \mathbb{S}^n .  $$
The critical points of this functional are \textit{antipodal geodesics}, i.e. geodesics $\gamma\colon [0,1]\to\mathbb{S}^n$ such that $\Dot{\gamma}(1) = -\Dot{\gamma}(0)$.
In Section \ref{sec_antipodal_geodesics} we show that the energy functional is a Morse-Bott function if we take $g$ to be the standard metric on $\mathbb{S}^n$.
All critical manifolds of the energy functional for the standard metric are diffeomorphic to the unit tangent bundle $U\mathbb{S}^n$.
In Section \ref{sec_completing_manifolds} we use the notion of \textit{completing manifolds} to show that the energy function for the standard metric is even a \textit{perfect} Morse-Bott function.
In particular this implies that the homology of $P_a\mathbb{S}^n$ can be computed from the homology the critical manifolds.
We show in Theorem \ref{theorem_homology_additively} that there is an isomorphism
$$   \mathrm{H}_{\bullet}(P_a\mathbb{S}^n )  \cong   \bigoplus_{i=0}^{\infty} \mathrm{H}_{\bullet - 2i(n-1)} \mathrm{H}_{\bullet}(U\mathbb{S}^n) .     $$
As the homology of the unit tangent bundle of the sphere is well-known this yields an explicit computation of $\mathrm{H}_{\bullet}(P_a\mathbb{S}^n)$.
We then use this description of the homology of $P_a\mathbb{S}^n$ to determine the extended loop product of even-dimensional spheres with rational coefficients.
It turns out that the product behaves even nicer if we consider the reduced homology of the free loop space.
\begin{TheoremC}[Corollary \ref{cor_extended_product_on_reduced_homology}]
    Consider the antipodal map $a\colon \mathbb{S}^n\to \mathbb{S}^n$ for $n$ even and take homology with rational coefficients.
    The extended loop product induces a product on the direct sum
    $\widetilde{\mathbb{H}}_{\bullet}(\mathbb{S}^n) := \widetilde{\mathrm{H}}_{\bullet + n}(\Lambda \mathbb{S}^n) \oplus \mathrm{H}_{\bullet+ n}(P_a\mathbb{S}^n)$ and there is an isomorphism of algebras
    $$     \widetilde{\mathbb{H}}_{\bullet}(\mathbb{S}^n) \cong \frac{\mathbb{Q}\langle A,U\rangle}{(A^2 , AU + UA)} \quad \text{with}\,\,\,|A| = -n,\,\, |U| = n-1 .       $$
\end{TheoremC}
It is remarkable that - up to signs -  the extended loop product on even-dimensional spheres behaves formally like the Chas-Sullivan product on an odd-dimensional sphere, since for $n$ odd we have by \cite{cohen:2003} an isomorphism of algebras 
$$(\mathrm{H}_{\bullet+n}(\Lambda \mathbb{S}^n),\wedge_{\mathrm{CS}}) \cong  \frac{\mathbb{Q}[A,U]}{(A^2) } \quad \text{with}\,\,\,|A| = -n,\,\, |U| = n-1.  $$

Next, we turn to co-operations on the homology $\mathrm{H}_{\bullet}(\Lambda \mathbb{S}^n)\oplus \mathrm{H}_{\bullet}(P_a\mathbb{S}^n)$.
In the general situation of a closed oriented manifold $M$ with a smooth fixed point free involution $f\colon M\to M$ and for homology with field coefficients we define maps 
$$    \vee_f \colon \mathrm{H}_{\bullet}(\Lambda M)\to \mathrm{H}_{\bullet}(P_f M)\otimes \mathrm{H}_{\bullet}(P_f M)  \quad \text{and}\quad   \vee_{f,l}\colon \mathrm{H}_{\bullet}(P_f M)\to \mathrm{H}_{\bullet}(\Lambda M,M)\otimes \mathrm{H}_{\bullet}(P_f M) $$ as well as $$ \vee_{f,r}\colon \mathrm{H}_{\bullet}(P_f M)\to \mathrm{H}_{\bullet}(P_f M)\otimes \mathrm{H}_{\bullet}(\Lambda M,M)     $$
which we think of as a copairing and left and right comodule structures.
These maps define an extension of the string topology coproduct  $$\vee\colon \mathrm{H}_{\bullet}(\Lambda M,M)\to \mathrm{H}_{\bullet}(\Lambda M,M)\otimes \mathrm{H}_{\bullet}(\Lambda M,M)$$ to a coproduct on the direct sum $\mathrm{H}_{\bullet}(\Lambda M,M)\oplus \mathrm{H}_{\bullet}(P_f M)$.
Dually, these co-operations induce a product structure on the cohomology $\mathrm{H}^{\bullet}(\Lambda M,M) \oplus \mathrm{H}^{\bullet}(P_f M)$ such that the Goresky-Hingston algebra is embedded as a subalgebra.

In the case that the involution $f$ is fixed point free and homotopic to the identity we show that the copairing is a lift of the string topology coproduct, see Theorem \ref{theorem_copairing_as_lift}.
We then use the completing manifolds on even-dimensional spheres to give a full computation of the copairing as well as of the comodule structures, see Theorem \ref{theorem_computation_copairing_even_sphere} and Proposition \ref{prop_computation_comodule_even_sphere}.

In the last section we apply the explicit computation of the extended products to show the following resonance theorem.
\begin{TheoremD}[Theorem \ref{theorem_resonance_even} and Theorem \ref{theorem_resonance_odd}]
    Let $g$ be a Riemannian metric on $\mathbb{S}^n$ and consider the induced energy functional on $\Lambda \mathbb{S}^n\sqcup P_a\mathbb{S}^n$.
    There are constants $\overline{\alpha},\beta > 0$ such that
    $$    | \mathrm{deg}(X) - \overline{\alpha}\,\mathrm{cr}(X)| \leq \beta   $$
     for all $X\in\widehat{\mathrm{H}}_{\bullet}(\mathbb{S}^n) = \mathrm{H}_{\bullet}(\Lambda \mathbb{S}^n) \oplus \mathrm{H}_{\bullet}(P_a\mathbb{S}^n)$.
\end{TheoremD}
Here, the critical values of homology classes are taken with respect to the function $\mathcal{L}\colon \Lambda \mathbb{S}^n\sqcup P_a\mathbb{S}^n\to [0,\infty)$ defined by $\mathcal{L}(\gamma) = \sqrt{E(\gamma)}$.
If $g$ is a $\mathbb{Z}_2$-invariant Riemannian metric on $\mathbb{S}^n$ then 
antipodal geodesics on $\mathbb{S}^n$ are mapped to non-contractible closed geodesic in $\mathbb{R}P^n$.
Similarly, closed geodesic in $\mathbb{S}^n$ are mapped to contractible closed geodesic in $\mathbb{R}P^n$.
In particular, critical levels of the energy functional on $\Lambda \mathbb{S}^n\sqcup P_a\mathbb{S}^n$ correspond bijectively to energies of closed geodesics in $\mathbb{R}P^n$.
The Morse index of a critical point in $\Lambda \mathbb{S}^n\sqcup P_a\mathbb{S}^n$ agrees with the index of the corresponding closed geodesic in $\mathbb{R}P^n$.
Moreover, for a homology class $X\in\mathrm{H}_{\bullet}(\Lambda \mathbb{S}^n)\oplus \mathrm{H}_{\bullet}(P_a\mathbb{S}^n)$ the index of a critical point at level $\mathrm{cr}(X)$ is approximately equal to the degree of $X$.
Hence, the resonance theorem gives us information about the index growth of closed geodesics in $\mathbb{R}P^n$ compared to its length.
See Section \ref{sec_resonance_theorem} for a more detailed discussion of the resonance theorem.

The constant $\overline{\alpha}$ in Theorem D is called the \textit{global mean frequency}.
The global mean frequency depends on the Riemannian metric $g$ on $\mathbb{S}^n$.
If $g$ is a $\mathbb{Z}_2$-invariant metric on $\mathbb{S}^n$ and $\gamma$ is a closed geodesic in $\mathbb{R}P^n$ for the induced metric $g'$ on $\mathbb{R}P^n$ then we let $\alpha_{\gamma}$ be the \textit{average index} which is defined by
$$  \alpha_{\gamma} = \lim_{k\to\infty} \frac{\mathrm{ind}(\gamma^k)}{k} .   $$
Furthermore we define the \textit{mean frequency} of $\gamma$ to be $\overline{\alpha}_{\gamma} = \tfrac{\alpha_{\gamma}}{L(\gamma)}$ where $L(\gamma)$ is the length of $\gamma$.
Theorem D implies the following result which transfers the \textit{density theorem}, \cite[Theorem 1.2]{hingston:2013}, to our situation.
\begin{TheoremE}[Theorem \ref{theorem_density}]
     Let $g$ be a Riemannian metric on $\mathbb{R}P^n$ for $n$ even and let $\overline{\alpha}$ be the global mean frequency.
    Let $\epsilon > 0$ and let $\mathcal{S}$ be the set of closed geodesics $\gamma$ in $\mathbb{R}P^n$ with mean frequency $\overline{\alpha}_{\gamma}\in(\overline{\alpha}-\epsilon,\overline{\alpha}+\epsilon)$.
    Then we have
    $$ \sum_{\gamma\in\mathcal{S}} \frac{1}{\alpha_{\gamma}}  \geq \frac{1}{n-1}  .    $$
\end{TheoremE}

\medskip

This article is organized as follows.
In Section \ref{sec_involution_path_space_string_operations} we study the general situation of involutions on closed manifolds and define the pairing and the module structures on the homology of the path space $P_f M$.
After showing the invariance under homotopies of $f$ we can compute the extended loop product for odd-dimensional spheres.
In Section \ref{sec_antipodal_geodesics} we then turn to the space of antipodal paths on the sphere $P_a\mathbb{S}^n$ and study the energy functional on this space.
Using the energy functional and Morse-Bott theoretic methods we compute the homology of $P_a\mathbb{S}^n$ in Section \ref{sec_completing_manifolds}.
We use the explicit description of the homology to compute the extended loop product of even-dimensional spheres in Section \ref{sec_computation_spheres}.
In Section \ref{sec_extensions_of_coproduct} we turn to dual structures, i.e. to the copairing and the comodule structures.
We show that the copairing can be used to define a lift of the string topology coproduct in case that the involution $f$ is homotopic to the identity.
We then compute the co-operations in the case of even-dimensional spheres.
Finally, in Section \ref{sec_resonance_theorem} we use the computations of the extended loop product and coproduct to show a resonance theorem for closed geodesics on real projective space and a density theorem for closed geodesics on real projective spaces of even dimension.

\medskip

\noindent \textbf{Acknowledgments:} 
The author wants to thank Marius Amann, Philippe Kupper and Hans-Bert Rademacher for helpful discussions and suggestions regarding this manuscript.

Some of the work on this manuscript was carried out during a stay of the author at the Copenhagen Centre for Geometry and Topology supported by the Deutsche Forschungsgemeinschaft (German Research Foundation) – grant agreement number 518920559.
The author is grateful for the support by the Danish National Research Foundation through the Copenhagen Centre for Geometry and Topology (DNRF151).

\section{Involutions of manifolds and extensions of the loop product}\label{sec_involution_path_space_string_operations}

In this section we introduce the space of paths in a manifold with endpoints determined by an involution and study certain product structures on its homology.
In the first two subsections we introduce these path spaces and the string topology operations.
We then show the invariance of these operations under homotopies of the involution and eventually apply this invariance to compute the extended loop product for many odd-dimensional spheres.

\subsection{The space of paths with respect to an involution}\label{subsec_space_of_paths}

Let $M$ be a closed manifold of dimension $n$ and assume that $g$ is a Riemannian metric on $M$.
Denote the unit interval by $I = [0,1]$.
We define the \textit{free path space of $M$} as $$  PM = \big\{ \gamma : I\to M\,|\, \gamma \,\,\text{absolutely continuous}, \,\, \int_0^1 g_{\gamma(t)}(\Dot{\gamma}(t),\Dot{\gamma}(t)) \,\mathrm{d}t < \infty \big\} \,   .$$
See \cite[Definition 2.3.1]{klingenberg:1995} for the notion of absolutely continuous curves in a smooth manifold.
The free path space $PM$ can be endowed with the structure of a Hilbert manifold and it can be shown to be independent of the choice of $g$, see \cite[Section 2.3]{klingenberg:1995}.

We note that the evaluation map $\mathrm{ev}\colon PM\to M\times M, \gamma\mapsto (\gamma(0),\gamma(1))$ is both a submersion as well as a fibration.
We shall consider various submanifolds of $PM$ which are obtained by pulling back the fibration $\mathrm{ev}$ along smooth embeddings $f\colon M\to M\times M$.
In fact, if we take the pullback along the diagonal $\Delta\colon M\to M\times M$ we obtain the \textit{free loop space of} $M$
$$  \Lambda M = \{\gamma \in PM \,|\, \gamma(0 ) = \gamma(1)\}  . $$
The underlying manifold $M$ can be seen as a submanifold of $\Lambda M$ via the identification with the trivial loops, see \cite[Proposition 1.4.6]{klingenberg:78}.
Let $f\colon M\to M$ be a smooth involution.
Note that $f$ is allowed to have fixed points, in particular it can also be the identity map.
Denote by $\Delta_f$ the set
$$   \Delta_f =  \{(x,f(x))\in M\times M\,|\, x\in M\} \subseteq M\times M     $$
which is a codimension $n$ submanifold of $M\times M$.
We pullback the path fibration $\mathrm{ev}\colon PM\to M\times M$ along the inclusion of $\Delta_f$, to obtain the space
$$    P_f M  =  \{ \gamma\in PM\,|\, (\gamma(0),\gamma(1))\in \Delta_f\} .     $$
If we take $f = \mathrm{id}_M$ we recover the free loop space $\Lambda M$.
Note that $P_f M$ comes with evaluation maps $\mathrm{ev}_0\colon P_f M\to M$ and $\mathrm{ev}_1\colon P_f M\to M$ given by
$$   \mathrm{ev}_0(\gamma) = \gamma(0) \quad \text{and}\quad \mathrm{ev}_1(\gamma) = \gamma(1)     .  $$
for $\gamma\in P_f M$
These maps are related by the identity $\mathrm{ev}_1 = f\circ \mathrm{ev}_0$.
Moreover, both maps are submersions as well as fibrations with typical fiber the based loop space $\Omega M$.

As we shall see next the space $P_f M$ can be understood as a universal covering space in certain situations.
Let $M$ be a closed simply connected oriented manifold and let $f\colon M\to M$ be a fixed point free involution.
    Let $N$ be the quotient $N=M/\mathbb{Z}_2$ where the $\mathbb{Z}_2$-action is induced by $f$.
    The free loop space of $N$ has two connected components $\Lambda_0 N$ and $\Lambda_1 N$ where the loops in $\Lambda_0 N$ are contractible and the loops in $\Lambda_1 N$ are not contractible.
\begin{prop}\label{prop_involution_path_space_is_universal_covering}
    Let $M$ be a closed $2$-connected oriented manifold and let $f\colon M\to M$ be a fixed point free involution.
    Let $p\colon M\to N = M/\mathbb{Z}_2$ be the universal $2$-fold covering.
    Then the space $P_f M$ is the universal covering space of $\Lambda_1 N$ with covering map $\pi_1\colon P_f M\to \Lambda_1 N$ given by $\pi(\gamma) = p\circ \gamma$. Similarly, the free loop space $\Lambda M$ is the universal covering space of $\Lambda_0 N$ with covering map $\pi_0\colon \Lambda M\to \Lambda_0 N$ given by $\pi(\gamma) = p\circ \gamma$.
\end{prop}
\begin{proof}
    We first show that $P_f M$ is simply connected.
    Consider the long exact sequence of homotopy groups 
    $$    \ldots \to \pi_1(\Omega M)\to \pi_1(P_f M) \to \pi_1(M) \to  \ldots         $$
    By assumption we have $\pi_1(M) \cong \{0\}$ as well as $\pi_1(\Omega M)\cong \pi_2(M) \cong \{0\}$.
    Consequently, $\pi_1(P_f M)=\{0\}$.
    For the covering property we take an open set $U\subseteq M$ such that $U\cap f(U) = \emptyset$.
    Set $V = p(U)$ and consider the open sets $   P_f M|_U = (\mathrm{ev}_0)^{-1}(U)  $ and $\Lambda N|_V = \mathrm{ev}_{\Lambda_1 N}^{-1}(V)$ where $\mathrm{ev}_{\Lambda_1 N}\colon \Lambda_1 N\to N$ is the evaluation map of the free loop space of $N$.
    One checks that the restriction
    $$   \pi|_{P_fM|_U} \colon P_f M_U \to \Lambda_1 N|_V     $$
    is a homeomorphism and the covering property follows.
    The statement about the component $\Lambda_0 N$ is shown analogously.    
\end{proof}
Throughout this paper the following example is of central importance.
Consider the sphere $\mathbb{S}^n$ and let $f=a\colon \mathbb{S}^n\to \mathbb{S}^n$ be the \textit{antipodal map} $a(p) = -p$.
We note the following special case of the above Proposition.
\begin{cor}\label{cor_antipodal_universal_covering_space}
    Let $n\geq 3$ and let $p\colon \mathbb{S}^n\to\mathbb{R}P^n$ be the projection.
    The space of antipodal paths on the sphere $P_a\mathbb{S}^n$ is the universal covering space of the component of non-contractible loops of the free loop space $\Lambda \mathbb{R}P^n$ with universal covering map $\pi\colon P_a\mathbb{S}^n\to \Lambda \mathbb{S}^n$, $\pi(\gamma) =p\circ \gamma$.
\end{cor}
\begin{remark}
    Note that in the example $M =\mathbb{S}^2$ with $f = a\colon \mathbb{S}^2\to \mathbb{S}^2$ the antipodal map we have $\pi_1(\Lambda \mathbb{S}^2) \neq \{0\}$.
    We shall later see that $\mathrm{H}_1(P_a \mathbb{S}^2) \cong \mathbb{Z}_2$ and hence $P_a \mathbb{S}^2$ is not simply connected either.
    Thus, neither $\Lambda \mathbb{S}^2$ nor $P_a \mathbb{S}^2$ can be universal covering spaces of one of the components of $\Lambda \mathbb{R}P^2$.
\end{remark}

In string topology one important operation is induced by the circle action $\chi\colon \mathbb{S}^1\times \Lambda M\to \Lambda M$.
This action is given by
$$   {\chi}_f(s,\gamma)(t) = \gamma(t+s) \quad \text{for}\,\,\,s\in\mathbb{S}^1,\,\,\gamma\in \Lambda M ,    $$
i.e. the loop gets rotated along itself.
If $f\colon M\to M$ is an involution then there is also a circle action on $P_f M$ which is defined as the map $\chi\colon \mathbb{S}^1\times P_f M\to P_f M$ given by
$$      \chi_f (s,\gamma) = \begin{cases}
    \gamma_a^s , & 0\leq s \leq \tfrac{1}{2} \\
    \gamma_b^s , & \tfrac{1}{2}\leq s \leq 1    
\end{cases}      $$
where
$$  \gamma_a^s ( t )  =  \begin{cases}
    \gamma(t + 2s), & t\in [0,1-2s] \\
    f(\gamma(t+2s-1)), & t\in [1-2s,1]
\end{cases}     $$
and 
$$
\gamma_b^s(t) =  \begin{cases}
    f(\gamma(t+2s-1)), & t\in [0,2-2s]\\
    \gamma(t+2s-2), & t\in [2- 2s,1] .
\end{cases}
$$
Note that for $f = \mathrm{id}_M$ this does not give the usual action $\widetilde{\chi}$ but rather 
$$   \chi_{\mathrm{id}} (s,\gamma)(t) = \gamma(t +2s)   \quad \text{for}\,\,\,s,t\in\mathbb{S}^1,\,\, \gamma\in \Lambda M. $$

\subsection{Loop pairing and module structure}
We shall now construct a pairing 
$$   \wedge_f \colon \mathrm{H}_{\bullet}(P_f M)\otimes \mathrm{H}_{\bullet}(P_f M) \to \mathrm{H}_{\bullet-n} (\Lambda M)      $$
where $n = \mathrm{dim}(M)$.
The idea behind the pairing is the trivial observation, that if we have a pair of paths $\gamma,\sigma\in P_f M$ with $\gamma(1) =\sigma(0)$, then their concatenation $\eta = \mathrm{concat}_2 (\gamma,\sigma)$ satisfies
$$    \eta(1) = \sigma(1) = f(\sigma(0)) = f(\gamma(1)) = f^2(\gamma(0)) = \gamma(0) = \eta(0),    $$
i.e. we end up with a loop in $M$.
Moreover, we shall also define left and a right module structure on $\mathrm{H}_{\bullet}(P_f M)$ over the Chas-Sullivan ring of $M$ following a similar idea.
The definition of the string topology operations in this section stays close to the way that the Chas-Sullivan product is defined in \cite{hingston:2017}.

Let $M$ be a closed oriented manifold and let $f\colon M\to M$ be a smooth involution.
We choose an $f$-invariant Riemannian metric $g$ on $M$ and denote the induced distance function by $\mathrm{d}\colon M\times M\to [0,\infty)$.
Define the space of concatenateable paths as
$$      C_f = (\mathrm{ev}_1\times \mathrm{ev}_0)^{-1} (\Delta M) = \{(\gamma,\sigma)\in P_fM\times P_f M\,|\, \gamma(1) =\sigma(0)\} .     $$
This is a codimension $n$ submanifold of $P_f M\times P_f M$ and thus has a tubular neighborhood.
We shall describe this tubular neighborhood quite explicitly.

Let $U_M$ be the tubular neighborhood of the diagonal in $M$ given by
$$   U_M = \{ (p,q)\in M\times M\,|\, \mathrm{d}(p,q) < \epsilon\}      $$
for some small $\epsilon> 0$.
Denote the Thom class of the normal bundle of the diagonal by 
$   \tau_M\in\mathrm{H}^n(U_M,U_M\setminus M) $. 
As an orientation convention we demand that $\tau_M \cap [M\times M] = \Delta_* [M]$ where $\Delta\colon M\to M\times M$ is the diagonal map, see \cite{hingston:2017}.
We define 
$$    U_C =    (\mathrm{ev}_1\times \mathrm{ev}_0)^{-1}(U_M ) = \{(\gamma,\sigma)\in P_f M\times P_f M\,|\, \mathrm{d}(\gamma(1),\sigma(0))<\epsilon\}.    $$
This is a tubular neighborhood of $C_f$.
We pull back the Thom class $\tau_M$ to get a class
$$   \tau_C =  (\mathrm{ev}_1\times \mathrm{ev}_0)^* \tau_M \in \mathrm{H}^n(U_C,U_C\setminus C_f ) .     $$
We now define a map $R^f\colon U_C\to C_f$ which is homotopic to the retraction $U_C\to C_f$ induced by the tubular neighborhood, respectively the normal bundle.
Let $(\gamma,\sigma)\in U_C$ and define $\sigma' = \mathrm{concat}( \overline{\gamma(1)\sigma(0)}, \sigma, \overline{\sigma(1)\gamma(0)})$.
Here, the path $\overline{pq}$ for $p,q\in M$, $\mathrm{d}(p,q)< \epsilon$ is the unique length-minimizing geodesic connecting $p$ and $q$.
Since the metric is $f$-invariant it follows from $\mathrm{d}(\gamma(1),\sigma(0))< \epsilon$ that also $\mathrm{d}(\sigma(1),\gamma(0))< \epsilon$.
Hence, $\sigma'$ is well-defined and we set $R^f(\gamma,\sigma) = (\gamma,\sigma')$.
This definition is inspired by the retraction maps in \cite{hingston:2017}.
With similar methods as in \cite[Lemma 2.3]{hingston:2017} one shows that the map $R^f$ is indeed homotopic to the retraction given by the data of the tubular neighborhood. 
Finally, there is a map $\mathrm{concat}\colon C_f\to P_f M$ which is defined by
$$   \mathrm{concat}(\gamma,\sigma)(t) = \begin{cases}
    \gamma(2t), & 0\leq t \leq \tfrac{1}{2}
 \\
 \sigma(2t-1), & \tfrac{1}{2}\leq t \leq 1 .
 \end{cases}     $$
 We now define the pairing on the homology of $P_f M$.
\begin{definition}
    Let $M$ be a closed oriented manifold and let $f\colon M\to M$ be an involution.
    Take homology with coefficients in a commutative ring $R$.
    The \textit{loop pairing} on the homology of $P_f M$ is defined as the following composition of maps.
    \begin{eqnarray*}
        \wedge_f \colon  \mathrm{H}_{i}(P_f M)\otimes \mathrm{H}_j(P_f M) &\xrightarrow{  (-1)^{n(n-i)}\times   }& \mathrm{H}_{i+j}(P_f M\times P_f M)
     \\ 
     &\xrightarrow{ \hphantom{coninci}\vphantom{r}\hphantom{conci} }& \mathrm{H}_{i+j}(P_f M\times P_f M,P_f M\times P_f M\setminus C_f) 
     \\
    &\xrightarrow{\hphantom{cii}\text{excision}\hphantom{ici}}&
    \mathrm{H}_{i+j}(U_{C}, U_{C}\setminus C_f) \\
    &\xrightarrow{\hphantom{coci} \tau_{C}\cap \hphantom{coci} }& \mathrm{H}_{i+j-n}(U_{C})
    \\ &\xrightarrow[]{\hphantom{conci}R^f_* \hphantom{conci}} & \mathrm{H}_{i+j-n}(C_f) 
    \\ & \xrightarrow[]{\hphantom{ici}\mathrm{concat}_*\hphantom{ici}} & \mathrm{H}_{i+j-n}(\Lambda M) \, .
    \end{eqnarray*}
\end{definition}

With a similar idea we now define the left and right module pairing.
Define 
$$   D_l^f = (\mathrm{ev}_{\Lambda}\times \mathrm{ev}_0)^{-1}(\Delta M) = \{(\gamma,\sigma)\in \Lambda M\times P_f M\,|\, \gamma(0) =\sigma(0) \}    $$
and 
$$  D_r^f =  (\mathrm{ev}_1\times \mathrm{ev}_{\Lambda})^{-1}(\Delta M) = \{(\gamma,\sigma)\in P_f M\times \Lambda M\,|\, \gamma(1) =\sigma(0) \} .     $$
As before we obtain tubular neighborhoods $U_l$ and $U_r$ of $D_l^f$ and $D_r^f$, respectively by defining
$$  U_l = (\mathrm{ev}_{\Lambda}\times \mathrm{ev}_0)^{-1}(U_M) \quad \text{and}\quad   U_r =  (\mathrm{ev}_1\times \mathrm{ev}_{\Lambda})^{-1}(U_M)  .   $$
We define cohomology classes
$$    \tau_l =  (\mathrm{ev}_{\Lambda}\times \mathrm{ev}_0)^* \tau_M \in \mathrm{H}^n(U_l,U_l\setminus D_l^f) \quad \text{and}\quad  \tau_r = (\mathrm{ev}_1\times \mathrm{ev}_{\Lambda})^* \tau_M \in\mathrm{H}^n(U_r, U_r\setminus D_r^f) .  $$
As before we can define homotopy retraction maps $R_l\colon U_l\to D_l^f$ and $R_r\colon U_r\to D_r^f$.
Moreover, the concatenation map now takes the form
$$   \mathrm{concat} \colon D_l^f \to P_f M \quad \text{as well as}\quad  \mathrm{concat}\colon D_r^f \to P_f M .   $$
In complete analogy to the above we define the following pairings.
\begin{definition}
    Let $M$ be a closed oriented manifold and let $f\colon M\to M$ be an involution.
    Take homology with coefficients in a commutative ring $R$.
    We define the following pairings
    \begin{eqnarray*}
        *_l \colon  \mathrm{H}_{i}(\Lambda M)\otimes \mathrm{H}_j(P_f M) &\xrightarrow{  (-1)^{n(n-i)}\times   }& \mathrm{H}_{i+j}(\Lambda M\times P_f M)
     \\ 
     &\xrightarrow{ \hphantom{coninci}\vphantom{r}\hphantom{conci} }& \mathrm{H}_{i+j}(\Lambda M\times P_f M,\Lambda M\times P_f M\setminus D_l^f) 
     \\
    &\xrightarrow{\hphantom{cii}\text{excision}\hphantom{ic}}&
    \mathrm{H}_{i+j}(U_{l}, U_{l}\setminus D_l^f) 
    \\
    &\xrightarrow{\hphantom{cocii} \tau_{l}\cap \hphantom{coici} }& \mathrm{H}_{i+j-n}(U_{l})
    \\ &\xrightarrow[]{\hphantom{coi}(R_l)_* \hphantom{coci}} & \mathrm{H}_{i+j-n}(D_l^f) 
    \\ & \xrightarrow[]{\hphantom{ici}\mathrm{concat}_*\hphantom{ici}} & \mathrm{H}_{i+j-n}(P_f M) \, 
    \end{eqnarray*}
    and
    \begin{eqnarray*}
        *_r \colon  \mathrm{H}_{i}(P_f M)\otimes \mathrm{H}_j(\Lambda M) &\xrightarrow{  (-1)^{n(n-i)}\times   }& \mathrm{H}_{i+j}(P_f M\times \Lambda M)
     \\ 
     &\xrightarrow{ \hphantom{coninci}\vphantom{r}\hphantom{conci} }& \mathrm{H}_{i+j}(P_f M\times \Lambda M, P_f M\times \Lambda M\setminus D_r^f) 
     \\
    &\xrightarrow{\hphantom{cii}\text{excision}\hphantom{ic}}&
    \mathrm{H}_{i+j}(U_{r}, U_{r}\setminus D_r^f) 
    \\
    &\xrightarrow{\hphantom{cocii} \tau_{r}\cap \hphantom{cocii} }& \mathrm{H}_{i+j-n}(U_{r})
    \\ &\xrightarrow[]{\hphantom{coi}(R_r)_* \hphantom{coci}} & \mathrm{H}_{i+j-n}(D_r^f) 
    \\ & \xrightarrow[]{\hphantom{ici}\mathrm{concat}_*\hphantom{ici}} & \mathrm{H}_{i+j-n}(P_f M) \, .
    \end{eqnarray*}
\end{definition}

Finally, we give the well-known definition of the Chas-Sullivan product.
Again, this is analogous to the above.
Consider the \textit{figure-eight space}
$$  \Lambda M\times_M \Lambda M = (\mathrm{ev}_{\Lambda}\times\mathrm{ev}_{\Lambda})^{-1}(\Delta M) = \{(\gamma,\sigma)\in \Lambda M\times \Lambda M\,|\,  \gamma(0) =\sigma(0)\}    .$$
This has a tubular neighborhood
$   U_{\mathrm{CS}} =  (\mathrm{ev}_{\Lambda}\times\mathrm{ev}_{\Lambda})^{-1}(U_M)   $
and it comes with a Thom class
$$   \tau_{\mathrm{CS}} = (\mathrm{ev}_{\Lambda}\times\mathrm{ev}_{\Lambda})^*\tau_M \in\mathrm{H}^n(U_{\mathrm{CS}},U_{\mathrm{CS}}\setminus \Lambda M\times_M\Lambda M) .   $$
Moreover, we again have a retraction $R_{\mathrm{CS}}\colon U_{\mathrm{CS}}\to \Lambda M\times_M \Lambda M$ and the concatenation takes $\Lambda M\times_M\Lambda M$ to $\Lambda M$.
\begin{definition}
    Let $M$ be a closed oriented $n$-dimensional manifold.
The \textit{Chas-Sullivan product} is defined as the composition
\begin{eqnarray*}
     \wedge\colon  \mathrm{H}_{i}(\Lambda M)\otimes \mathrm{H}_j(\Lambda M) &\xrightarrow{  (-1)^{n(n-i)}\times   }& \mathrm{H}_{i+j}(\Lambda M \times \Lambda M)
     \\
     &\xrightarrow{ \hphantom{coninci}\vphantom{r}\hphantom{conci} }& \mathrm{H}_{i+j}(\Lambda M\times \Lambda M, \Lambda M\times \Lambda M\setminus \Lambda M\times_M \Lambda M)
     \\
    &\xrightarrow{\hphantom{cii}\text{excision}\hphantom{ic}}&
    \mathrm{H}_{i+j}(U_{\mathrm{CS}}, U_{\mathrm{CS}}\setminus \Lambda M\times_M\Lambda M) \\
    &\xrightarrow{\hphantom{coci} \tau_{\mathrm{CS}}\cap \hphantom{coci} }& \mathrm{H}_{i+j-n}(U_{\mathrm{CS}})
    \\ &\xrightarrow[]{\hphantom{coi}(R_{\mathrm{CS}})* \hphantom{coi}} & \mathrm{H}_{i+j-n}(\Lambda M\times_M \Lambda M) 
    \\ & \xrightarrow[]{\hphantom{ici}\mathrm{concat}_*\hphantom{ici}} & \mathrm{H}_{i+j-n}(\Lambda M) \, .
    \end{eqnarray*}
\end{definition}

In the following we shall now study how all these four structures interact with each other.
Recall that the Chas-Sullivan product defines a unital, associatve and graded commutative product on $\mathrm{H}_{\bullet}(\Lambda M)$, see \cite{chas:1999}.

\begin{prop}\label{prop_bimodule_structure}
    The pairings $*_l$ and $*_r$ endow $\mathrm{H}_{\bullet}(P_f M)$ with a bi-module structure over the Chas-Sullivan ring.
\end{prop}
This can be seen similarly to how one proves that the Chas-Sullivan product and intersection product of a closed manifold are associative, see \cite[Proposition B.1]{hingston:2017}.
A similar proof for a string topology type product can be found in \cite[Theorem 2.5]{stegemeyer2023string}.
Below we shall how the left- and the right-module structure are related to each other by comparing the pairing $A*_l X$ and $X*_r A$ for $A\in\mathrm{H}_{\bullet}(\Lambda M)$ and $X\in\mathrm{H}_{\bullet}(\Lambda M)$, see Proposition \ref{prop_commutativity_module}.

\begin{remark}\label{remark_id_just_cs}
    Recall that in the above definitions the involution $f\colon M\to M$ does not need to be fixed point free.
    In particular we can consider the pairing and the module structures in the case $f = \mathrm{id}_M$ and one sees directly that we get the Chas-Sullivan product in each case.
\end{remark}

We now want to study how the pairing $\wedge_f$ interacts with the other three structures.
The main goal will be to show that the direct sum $\mathrm{H}_{\bullet}(\Lambda M)\oplus \mathrm{H}_{\bullet}(P_f M)$ becomes an associative unital algebra.
In a brief interlude we turn to a purely algebraic situation.
Let $(A,\wedge_A)$ be a unital, associative graded algebra over a commutative ring $R$.
Moreover, let $B$ be an $A$-bimodule with module structures
$$   *_l \colon A\otimes B\to B \quad \text{and}\quad *_r \colon B\otimes A\to B .   $$
Assume we have also have a pairing $\wedge_B\colon B\otimes B\to A$.
We want to make the direct sum $A\oplus B$ into a unital, associative algebra.
In order to ensure that we end up with an associative algebra we make the following definition.
\begin{definition}\label{def_adapted}
    Let $(A,\wedge_A)$ be an associative algebra and $B$ be an $A$-bimodule.
    If $\wedge_B\colon B\otimes B\to A$ is a pairing, then we say that it is \textit{adapted} to $A$ if 
    \begin{eqnarray*}
           &  a \wedge_A (y\wedge_B z) = (a *_l y)\wedge_B z & \quad   (x\wedge_B y) \wedge_A a = x \wedge_B (y*_r a)        ,
           \\
           &  (x *_r a)\wedge_B z = x \wedge_B ( a *_l z) & \text{ and }  \quad 
           (x\wedge_B y)*_l z = x *_r (y\wedge_B z) 
    \end{eqnarray*}
    for all $a\in A$, $x,y,z\in B$.
\end{definition}
A direct check then shows the folllowing.
\begin{prop}
    Let $A$ be an associative and unital algebra and $B$ an $A$-bimodule.
    If $\wedge_B\colon B\otimes B\to A$ is an adapted pairing, then the direct sum $A\oplus B$ is an associative and unital algebra with multiplication given by
    $$    ( a, x ) \wedge (b, y) =   \big( a\wedge_A b + x\wedge_B y \,, \, a *_l y + x *_r b\big)     \quad \text{for}\,\,\,a,b\in A, x,y\in B . $$
    The unit is given by $(1,0)$ with $1\in A$ being the unit in $A$.
    Moreover, the algebra $(A,\wedge_A)$ is a subalgebra of $A\oplus B$.
\end{prop}
We call the algebra structure on $A\oplus B$ the \textit{extension algebra induced by} $\wedge_B$.

\begin{remark}\label{example_trivial_extension}
    Let $A$ be a unital associative algebra.
    Of course, $A$ is a bimodule over itself.
    We set $B = A$ and consider it as an $A$-bimodule.
    Moreover, the multiplication in $A$ defines a pairing $B\otimes B\to A$.
    This is adapted to the multiplication in $A$ in the sense of Definition \ref{def_adapted} and thus we obtain a unital associative algebra structure on $A\oplus B \cong A\oplus A$ given by
    $$     (a,x) \wedge (b,y) =  (a\cdot b + x\cdot y, a\cdot y + x\cdot b)     $$
    where $\cdot$ denotes the original product on $A$.
    We call this extension algebra the \textit{trivial extension}.
    In particular we note that if an algebra $C$ is isomorphic to the trivial extension algebra of an algebra $A$ then there is an element $c\in C$ with $c\neq \pm 1$ such that $c^2 = 1$.
    Hence, the non-existence of such an element can tell us that a certain extension algebra is not a trivial extension.
\end{remark}

\begin{prop}
    Let $M$ be a closed oriented manifold and $f\colon M\to M$ an involution.
    Take homology with coefficients in a commutative ring $R$.
    The pairing $\wedge_f\colon \mathrm{H}_{\bullet}(P_f M)\otimes \mathrm{H}_{\bullet}(P_f M)\to \mathrm{H}_{\bullet-n}(\Lambda M)$ is adapted to the Chas-Sullivan product.
\end{prop}
Like Proposition \ref{prop_bimodule_structure} this can be shown analogously to how one proves the associativity of the intersection product and of the Chas-Sullivan product.
Since the proof is quite lengthy but not very illuminating we shall omit it here.

\begin{theorem}\label{theorem_extended_cs_product}
    Let $M$ be a closed oriented manifold and $f\colon M\to M$ be an involution.
    Take homology with coefficients in a commutative unital ring $R$.
    The direct sum $\mathrm{H}_{\bullet}(\Lambda M)\oplus \mathrm{H}_{\bullet}(P_f M)$ equipped with the multiplication induced by the Chas-Sullivan product, the left and right module structure as well as the involution pairing $\wedge_f$ yields a unital and associative algebra.
\end{theorem}

We shall frequently consider the degree shifted homology
$$   \widehat{\mathbb{H}}_{\bullet}(M) : =   \mathrm{H}_{\bullet+n}(\Lambda M)\oplus \mathrm{H}_{\bullet+ n }(P_f M) .     $$
On this graded $R$-module the extended Chas-Sullivan product is a product of degree $0$ and the original Chas-Sullivan product is graded commutative in the classical sense.
If not specified otherwise however, the degree of a homology class always refers to the usual grading, i.e. the point class lives in degree $0$.
We now discuss the commutativity of the pairing $\wedge_f$.
\begin{prop}\label{prop_commutativity_loop_pairing}
    Let $M$ be a closed oriented manifold of dimension $n$ and $f\colon M\to M$ an involution.
    Take homology with coefficients in a commutative unital ring $R$.
    Then the pairing $\wedge_f$ is graded commutative in the sense that
    $$       X\wedge_f Y =   \mathrm{deg}(f) (-1)^{(|X|-n)(|Y|-n)} Y\wedge_f X   $$
    for $X,Y\in\mathrm{H}_{\bullet}(P_f M)$
\end{prop}
\begin{proof}
    Let $T\colon P_f M\times P_f M\to P_f M \times P_f M$ be the swapping map.
    Assume that $M$ is endowed with a $\mathbb{Z}_2$-invariant metric, with respect to the $\mathbb{Z}_2$-action induced by the involution $f$.
    Further let $\mathcal{T}\colon \mathrm{H}_i(P_f M)\otimes \mathrm{H}_j(P_f M) \to \mathrm{H}_j(P_fM)\otimes \mathrm{H}_i(P_f M)$ be the map swapping the tensor factors.
    We consider the following diagram.
    $$
    \begin{tikzcd}
        \mathrm{H}_i(P_f M )\otimes \mathrm{H}_j(P_f M) \arrow[]{r}{\mathcal{T}} \arrow[]{d}{\times}
        & 
        \mathrm{H}_j(P_f M )\otimes \mathrm{H}_i(P_f M) \arrow[]{d}{\times}
        \\
        \mathrm{H}_{i+j}(P_f M\times P_f M) \arrow[]{r}{T_*}
        \arrow[]{d}{} 
        & 
        \mathrm{H}_{i+j}(P_f M\times P_f M)\arrow[]{d}{}
        \\
        \mathrm{H}_{i+j}(U_C,U_C\setminus C_f) \arrow[]{r}{T_*} \arrow[]{d}{\tau_C \cap}
        & \mathrm{H}_{i+j}(U_C,U_C\setminus C_f) \arrow[]{d}{\tau_C \cap}
        \\
        \mathrm{H}_{i+j-n}(U_C) \arrow[]{d}{(\mathrm{concat}\circ R_C)_*} \arrow[]{r}{T_*} & 
        \mathrm{H}_{i+j-n}(U_C) \arrow[]{d}{(\mathrm{concat}\circ R_C)_*}
        \\
        \mathrm{H}_{i+j-n}(\Lambda ) \arrow[]{r}{ = } & \mathrm{H}_{i+j-n}(\Lambda ).
    \end{tikzcd}
    $$
    We claim that this diagram commutes up to sign.
    The first square commutes up to sign $(-1)^{ij}$.
    The second square commutes, since the swapping map maps the pair $(U_C,U_C\setminus C_f)$ to itself by our assumption on the $\mathbb{Z}_2$-invariance of the metric.
    For the third square, let $X \in \mathrm{H}_{\bullet}(U_C,U_C\setminus C_f)$, then by naturality we have
    $$   \tau_C \cap T_* X =  T_*( T^* \tau_C \cap X) .    $$
    Recall that we have $\tau_C = (\mathrm{ev}_1\times \mathrm{ev}_0)^* \tau_M$.
    Moreover, we have for $(\gamma,\sigma)\in U_C$ that
    \begin{eqnarray*}
        (\mathrm{ev}_1\times \mathrm{ev}_0)\circ T(\gamma,\sigma) & =& (\sigma(1),\gamma(0)) \\ 
        &=&
        (f(\sigma(0)),f(\gamma(1)) \\
        &=& t\circ (f\times f)  (\gamma(1),\sigma(0)).
        \end{eqnarray*}
        This shows that $\mathrm{ev}_1\times \mathrm{ev}_0\circ T =  t\circ (f\times f) \circ (\mathrm{ev}_1\times \mathrm{ev}_0)$, where $t\colon M\times M\to M\times M$ is the swapping map.
        Consequently we have that 
        $$   T^* \tau_C  =  (\mathrm{ev}_1\times \mathrm{ev}_0)^*  (f\times f)^* t^* \tau_M   .  $$   
        It is well-known that the map $t\colon M^2\to M^2$ is orientation-reversing if and only if $n$ is odd, hence $t^*\tau_M = (-1)^n\tau_M$.
        Now consider the map $f\times f\colon M\times M\to M\times M$.
        Clearly we have 
        $$  (f\times f)_* [M\times M] = f_* [M]\times f_*[M] =  \mathrm{deg}(f)^2 [M\times M]  =  [M\times M] .$$
        Note that the map $f\times f$ restricts to a homeomorphism of pairs $f\times f\colon (U_M, U_M\setminus \Delta M) \to (U_M, U_M\setminus \Delta M)$ and hence we have $(f\times f)^* \tau_M = (-1)^a \tau_M$.
        We compute that
        \begin{eqnarray*}
            \Delta_* [M]  & = & \tau_M \cap [M\times M] \\
            &=& \tau_M \cap \left(   (f\times f)_* ([M]\times [M] )   \right)  \\
            &=&
            (f\times f)_* \left(   (f\times f)^* \tau_M \cap [M\times M]   \right)
            \\ &=& (-1)^a (f\times f)_* \Delta_* [M] \\
            &=& (-1)^a \Delta_* f_*[M] =  (-1)^a \mathrm{deg}(f) \Delta_* [M] .
        \end{eqnarray*}
        In the last equality we used that $(f\times f)\circ \Delta = \Delta \circ  f $.
        Consequently, we have $(-1)^a = \mathrm{deg}(f)$.
        Therefore we see that the third square commutes up to sign $(-1)^n\mathrm{deg}(f)$.
        The last square commutes because the underlying diagram of maps commutes up to homotopy.
        More precisely, note that the maps $\mathrm{concat}\colon C_f \to \Lambda$ and $\mathrm{concat}\circ T\colon C_f \to \Lambda M$ are homotopic by the homotopy $H\colon C_f\times [0,1]\to \Lambda M$ given by 
        $$    H((\gamma,\sigma),s) = \begin{cases}
                \gamma(2t + (1-s)) , & 0\leq t \leq \tfrac{s}{2} \\
                \sigma(2t-s) , & \tfrac{s}{2} \leq t \leq \tfrac{1+s}{2} \\
                \gamma(2t-s-1), & \tfrac{1+s}{2}\leq t \leq 1.
        \end{cases}    $$
        Let $X\in\mathrm{H}_i(P_f M)$ and $Y\in\mathrm{H}_j(P_f M)$.
        The composition down the left hand side of the diagram gives $(-1)^{n-ni} X\wedge_f Y$ while the composition down the right hand side gives $(-1)^{n-nj}Y\wedge_f X$. Collecting all the signs yields the identity $$X\wedge_f Y = \mathrm{deg}(f) (-1)^{(|X|-n)(|Y|-n)} Y\wedge_f X . $$
        This completes the proof.
\end{proof}
We now prove a similar commutativity result for the left and the right-module structure.
\begin{prop}\label{prop_commutativity_module}
    Let $M$ be a closed oriented manifold and $f\colon M\to M$ an involution.
    Take homology with coefficients in a commutative unital ring $R$.
    Denote by $\Lambda f\colon \Lambda M\to \Lambda M$ map induced by $f$ on the free loop space.
    Then for $X\in\mathrm{H}_{\bullet}(\Lambda M)$ and $A\in\mathrm{H}_{\bullet}(P_f M)$ it holds that 
    $$    X*_l  A   =  \mathrm{deg}(f) (-1)^{(|X|-n)(|Y|-n)}   \, A*_r ((\Lambda f)_*X ).    $$
\end{prop}
\begin{proof}
      As in the proof of Proposition \ref{prop_commutativity_loop_pairing} let $T\colon P_f M\times P_f M\to P_f M \times P_f M$ be the swapping map and endow $M$ with a $\mathbb{Z}_2$-invariant metric with respect to $f$.
    Let $\mathcal{T}\colon \mathrm{H}_i(\Lambda M)\otimes \mathrm{H}_j(P_f M) \to \mathrm{H}_j(P_f M)\otimes \mathrm{H}_i(\Lambda M)$ be the swapping map.
    Consider the following diagram.
    $$
    \begin{tikzcd}
        \mathrm{H}_i(\Lambda M )\otimes \mathrm{H}_j(P_f M) \arrow[]{r}{\mathcal{T}\circ (\Lambda f_* \otimes \mathrm{id})} \arrow[]{d}{\times}
        & [2.5em]
        \mathrm{H}_j(P_f M )\otimes \mathrm{H}_i(\Lambda M) \arrow[]{d}{\times}
        \\
        \mathrm{H}_{i+j}(\Lambda M\times P_f M) \arrow[]{r}{(T\circ (\Lambda f\times \mathrm{id}))_*}
        \arrow[]{d}{} 
        & 
        \mathrm{H}_{i+j}(P_f M\times\Lambda M)\arrow[]{d}{}
        \\
        \mathrm{H}_{i+j}(U_{D_l},U_{D_l}\setminus D_l^f) \arrow[]{r}{(T\circ (\Lambda f\times \mathrm{id})_*} \arrow[]{d}{\tau_l \cap}
        & \mathrm{H}_{i+j}(U_{D_r},U_{D_r}\setminus D_r^f) \arrow[]{d}{\tau_r \cap}
        \\
        \mathrm{H}_{i+j-n}(U_{D_l}) \arrow[]{d}{( {R}_l)_*} \arrow[]{r}{(T\circ (\Lambda f\times \mathrm{id})_*} & 
        \mathrm{H}_{i+j-n}(U_{D_r}) \arrow[]{d}{( R_{r})_*}
        \\
        \mathrm{H}_{i+j-n}( D_l ) \arrow[]{r}{ (T\circ (\Lambda f\times \mathrm{id}))_* } \arrow[]{d}{\mathrm{concat}_*} & \mathrm{H}_{i+j-n}(D_r) \arrow[]{d}{\mathrm{concat}_*}
        \\
        \mathrm{H}_{i+j-n}(P_f M)\arrow[]{r}{=} & \mathrm{H}_{i+j-n}(P_f M) .
    \end{tikzcd}
    $$
    We claim that this diagram commutes up to sign.
    As in the proof of Proposition \ref{prop_commutativity_loop_pairing} we get a sign $(-1)^{ij}$ in the first square. 
    The second square commutes and we claim that the third square commutes up to the sign $(-1)^n\mathrm{deg}(f)$.
    Let $X\in\mathrm{H}_{\bullet}(U_{D_l},U_{D_l}\setminus D_l^f)$.
    Then by naturality we have
    $$   \tau_r \cap ((T\circ (\Lambda f\times \mathrm{id}))_* X) =   (T\circ (\Lambda f\times \mathrm{id}))\big(     (T\circ (\Lambda f\times \mathrm{id}))^* \tau_r \cap X  \big)  . $$
    Recall that $\tau_r = (\mathrm{ev}_1\times \mathrm{ev}_{\Lambda})^*\tau_M$.
    A direct computation shows that
    $$    (\mathrm{ev}_1\times \mathrm{ev}_{\Lambda})\circ T\circ (f\times \mathrm{id}) = (f\times f)\circ t\circ (\mathrm{ev}_{\Lambda}\times \mathrm{ev}_0)  .     $$
    Thus like in the proof of Proposition \ref{prop_commutativity_loop_pairing} we see that
    $$    (T\circ (\Lambda f\times \mathrm{id}))^* \tau_r = (-1)^n\mathrm{deg}(f) \tau_l .      $$
    This shows that the third square commutes up to the sign $(-1)^n\mathrm{deg}(f)$.
    The commutativity of the fourth square is easy to see so we turn to the commutativity of the last square.
    We claim that the underlying diagram of maps
    $$
    \begin{tikzcd}
        D_l \arrow[]{r}{T\circ (\Lambda f\times \mathrm{id})} \arrow[swap]{dr}{\mathrm{concat}}
        & [2.5em] D_r \arrow[]{d}{\mathrm{concat}}
        \\
        & P_f M
    \end{tikzcd}
    $$
    commutes up to homotopy.
    Indeed we compute that for $(\gamma,\sigma)\in D_l$ we have
    $$   \mathrm{concat}\circ T\circ (\Lambda f\times \mathrm{id})(\gamma,\sigma) =   \mathrm{concat}(\sigma,f\circ \gamma)  . $$
    Recall that we have a circle action $\chi_f\colon \mathbb{S}^1\times P_f M\to P_f M$.
    Let $\chi_{f,\!\frac{1}{4}}= \chi_f(\tfrac{1}{4},\cdot)\colon P_f M\to P_f M$ be the action by $\tfrac{1}{4}\in \mathbb{S}^1$.
    Then one checks that
    $   \chi_{f,{\!\frac{1}{4}}}\circ \mathrm{concat}(\gamma,\sigma) = \mathrm{concat}(\sigma, f\circ \gamma)       $ for $(\gamma,\sigma)\in D_l$
    and thus $\chi_{f,\!\frac{1}{4}}\circ \mathrm{concat} = \mathrm{concat}\circ T\circ (\Lambda f\times \mathrm{id})$.
    Clearly, $\chi_{f,\!\frac{1}{4}}  \simeq \chi_{f,0} = \mathrm{id}_{P_f M}$ and hence we see that $\mathrm{concat} \simeq \mathrm{concat}\circ T\circ (\Lambda f\times \mathrm{id})$.
    This shows the last square commutes.
    Collecting all the signs then yields the claimed identity.
\end{proof}

\begin{remark}\label{remark_bv_algebra}
    Let $M$ be a closed oriented manifold and $f\colon M\to M$ a smooth involution.
    We have now established an associative and unital product on the direct sum $\widehat{\mathbb{H}}_{\bullet}(M) := \mathrm{H}_{\bullet+n}(\Lambda M)\oplus \mathrm{H}_{\bullet+n}(P_f M)$.
    With the $\mathbb{S}^1$-action on $P_f M$ one can build an operator $B\colon {\mathrm{H}}_{\bullet}(M) \to {\mathrm{H}}_{\bullet+1}(M)$ by setting $B(X) = (\chi_f)_*([\mathbb{S}^1]\times X)$ for $X\in\mathrm{H}_{\bullet}(P_f M)$.
    It is easy to see that $B^2 = 0$.    
    On $\mathrm{H}_{\bullet}(\Lambda M)$ there is a circle action and a corresponding operator $B\colon \mathrm{H}_{\bullet}(\Lambda M)\to \mathrm{H}_{\bullet+1}(\Lambda M)$ as well.
    It is therefore an natural question what structure the product on $\widehat{\mathbb{H}}_{\bullet}(M)$ and the operator $B$ form together.
    In the string topology on the free loop space it is well-known that one obtains a BV-algebra, i.e. the bracket
    $$   \{X,Y \} =  (-1)^{|X|-n}\Delta_*(X\wedge_{\mathrm{CS}} Y) - (-1)^{|X|-n}\Delta_*(X)\wedge_{\mathrm{CS}} Y - X\wedge_{\mathrm{CS}} \Delta_* (Y)     $$
    is a derivation in each variable, see \cite{cohen2006string}.
    This bracket is a graded Lie bracket and hence the homology of the free loop space becomes a Gerstenhaber algebra.
    Note however that the graded commutativity of the product is a crucial ingredient here.
    We have seen above that the extended loop product is in general not graded commutative although the failure to be graded commutative can be expressed in terms of the map $f$.
    Therefore the author conjectures that there might be a variant of a BV algebra where one sprinkles in additional signs which is indeed realized on $\widehat{\mathbb{H}}_{\bullet}(M)$.
\end{remark}

\subsection{Invariance of the pairing and the module structure}\label{subsec_invariance_loop_pairings}

We begin by showing that two involution path spaces $P_f M$ and $P_g M$ are homotopy equivalent if the involutions $f$ and $g$ are homotopic.
\begin{prop}\label{prop_homotopy_invariance_involution_path_spaces}
    Let $M$ be a closed oriented manifold and $f$ and $g$ two involutions.
    If $f$ and $g$ are homotopic, then the involution path spaces $P_f M$ and $P_g M$ are homotopy equivalent.
\end{prop}
\begin{proof}
    Let $H\colon M\times I\to M$ be a homotopy between $f$ and $g$, i.e.
    $$   H(p,0) = f(p) \quad \text{and}\quad H(p,1) = g(p)     $$
    for $p\in M$.
    Define $\eta\colon M\to PM$ to be the map $\eta(p)(s) = H(p,s)$ for $p\in M$, $s\in I$.
    We define $\Phi\colon P_f M\to P_g M$ and $\Psi\colon P_g M\to P_f M$ by
    $$   \Phi(\gamma) = \mathrm{concat}(\gamma,\eta(\gamma(0))) \quad \text{and}\quad \Psi(\gamma) = \mathrm{concat}(\gamma,\overline{\eta(\gamma(0))}) .      $$
    Here, the path $\overline{\sigma}$ means the reversed path of the path $\sigma \in PM$, i.e. $\overline{\sigma}(t) = \sigma(1-t)$.
    One checks that these maps are homotopy inverse to each other.
\end{proof}
\begin{cor}\label{cor_involution_homotopic_to_id_just_free_loop_space}
    Let $f\colon M\to M$ be an involution which is homotopic to the identity, then $P_f M \simeq \Lambda M$.
\end{cor}

We shall now show that the involution pairing as well as the left- and right-module structures on $P_f M$ can be obtained from the respective structures on $P_g M$ if $f$ and $g$ are homotopic involutions.
Note that next to the homotopy equivalence $\Phi$ in the proof of Proposition \ref{prop_homotopy_invariance_involution_path_spaces} there is another homotopy equivalence $\Phi'\colon P_f M\to P_g M$ given by
$$  \Phi'(\gamma) = \mathrm{concat}(\overline{\eta(\gamma(1))}, \gamma) \quad \text{for}\,\,\,\gamma\in P_f M .      $$
Indeed, we have that $\overline{\eta(\gamma(1))}$ is a path from $g(\gamma(1)) = g(f(\gamma(0)))$ to $f(\gamma(1)) = \gamma(0)$, hence this is well-defined.
Note that the induced maps in homology
$   \Phi_*, \Phi'_* \colon \mathrm{H}_i(P_f M)\to \mathrm{H}_i(P_g M)      $
do not need to agree.

\begin{prop}\label{prop_pairing__for_homotopic_maps}
    Let $M$ be a closed oriented manifold and let $f,g\colon M\to M$ be smooth 
    involutions.
    Take homology with coefficients in a commutative, unital ring $R$.
    If $f$ and $g$ are homotopic then the pairings with respect to $f$ and $g$ are related to each other by the commutativity of the diagram
    $$ 
    \begin{tikzcd}
        \mathrm{H}_i(P_f M)\otimes \mathrm{H}_j(P_f M) \arrow[r, "\Phi'_*\otimes\, \Phi_*", "\cong"']
        \arrow[]{d}{\wedge_f}
        & [2.5em]
        \mathrm{H}_i(P_g M)\otimes \mathrm{H}_j(P_g M)
        \arrow[]{d}{\wedge_g}
        \\
        \mathrm{H}_{i+j-n}(\Lambda M) \arrow[]{r}{ = }
        &
        \mathrm{H}_{i+j-n}(\Lambda M) .
    \end{tikzcd}
    $$
\end{prop}
\begin{proof}
    Throughout this proof we shall write $U_{C,f}$ for the tubular neighborhood of $C_f$ inside $P_fM \times P_f M$ and similarly we shall write $U_{C,g}$ for the tubular neighborhood of $C_g$ inside $P_g M\times P_g M$.
    Similarly, we write $\tau_{C,f}$ and $\tau_{C,g}$, respectively for the Thom classes.
    We claim that the following diagram commutes.
    $$
    \begin{tikzcd}
        \mathrm{H}_i(P_f M\times P_f M) \arrow[]{r}{(\Phi'\times \Phi)_*}
        \arrow[]{d}{}
        & [2.5em]
        \mathrm{H}_i(P_g M \times P_g M)
        \arrow[]{d}{}
        \\
        \mathrm{H}_i(U_{C,f}, U_{C,f}\setminus C_f) \arrow[]{r}{(\Phi'\times \Phi)_*}
        \arrow[]{d}{\tau_{C,f}\cap}
        &
        \mathrm{H}_i(U_{C,g},U_{C,g}\setminus C^g) \arrow[]{d}{\tau_{C,g}\cap}
        \\
        \mathrm{H}_i(U_{C_f}) \arrow[]{r}{(\Phi'\times \Phi)_*}
        \arrow[]{d}{(\mathrm{concat}\circ R^f)_*}
        &
        \mathrm{H}_i(U_{C^g}) 
        \arrow[]{d}{(\mathrm{concat}\circ R^g)_*}
        \\
        \mathrm{H}_i(\Lambda M) \arrow[]{r}{=} &
        \mathrm{H}_i(\Lambda M) .
    \end{tikzcd}
    $$
    The commutativity of the first square is clear.
    For the middle square, let $X\in\mathrm{H}_i(U_{C,f},U_{C,f}\setminus C_f)$.
    We need to show that
    $$  (\Phi'\times \Phi)_* (\tau_{C,f}\cap X) =  \tau_{C,g}\cap ((\Phi'\times \Phi)_* X) .        $$
    By naturality the right hand side of the above equation equals
    $$ \tau_{C,g}\cap ((\Phi'\times \Phi)_* X)  =   (\Phi'\times \Phi)^*  \big(  (\Phi'\times \Phi)^*\tau_{C,g} \cap X     \big)  .    $$
    Recall that $\tau_{C,g} = (\mathrm{ev}_1\times \mathrm{ev}_0)^* \tau_M$.
    Since $(\mathrm{ev}_1\times \mathrm{ev}_0) \circ (\Phi'\times \Phi) = \mathrm{ev}_1\times \mathrm{ev}_0$ we have $ (\Phi'\times \Phi)^*\tau_{C,g} = \tau_{C,f}$ and thus the second square commutes.
    For the last square we note that the diagram of spaces and maps
    $$   
    \begin{tikzcd}
        U_{C_f} \arrow[]{r}{\Phi'\times \Phi}
        \arrow[swap]{d}{\mathrm{concat}\circ R^f}
        & [2.5em]
        U_{C^g} \arrow[]{d}{\mathrm{concat}\circ R^g}
        \\
        \Lambda M\arrow[]{r}{ \Theta}
        & \Lambda M
    \end{tikzcd}
    $$
    commutes up to homotopy, where $\Theta\colon \Lambda M\to \Lambda M$ is the map
    $$ \Theta(\gamma) =  \mathrm{concat}( \overline{\eta(f(\gamma(0)))} , \gamma, \eta(f(\gamma(0)))) \quad \text{for}\,\,\, \gamma\in \Lambda M.       $$
    By the spaghetti trick there is a homotopy $\Theta\simeq \mathrm{id}_{\Lambda M}$ and thus the third square commutes as well.
    One sees from the definition of the pairing $\wedge_f$ that the commutativity of the diagram shows the claim.
\end{proof}
Similarly, one can show the following statement about the module structure.
We shall omit the proof since it follows the same strategy as the proof of the previous proposition.
\begin{prop}\label{prop_module_for_homotopic_maps}
     Let $M$ be a closed oriented manifold and let $f,g\colon M\to M$ be smooth involutions.
    Take homology with coefficients in a commutative, unital ring $R$.
    If $f$ and $g$ are homotopic then the module structures $*_l$ and $*_r$ with respect to $f$ and $g$ are related by the commutativity of the diagrams
    $$ 
    \begin{tikzcd}
        \mathrm{H}_i(\Lambda M)\otimes \mathrm{H}_j(P_f M) \arrow[r, " \mathrm{id}\, \otimes \,\Phi_*", "\cong"']
        \arrow[]{d}{*_l}
        & [2.5em]
        \mathrm{H}_i(\Lambda M)\otimes \mathrm{H}_j(P_g M)
        \arrow[]{d}{*_l}
        \\
        \mathrm{H}_{i+j-n}(P_f M) \arrow[]{r}{\Phi_* }
        &
        \mathrm{H}_{i+j-n}(P_g M) 
    \end{tikzcd}
    $$
    and
    $$ 
    \begin{tikzcd}
        \mathrm{H}_i(P_f M)\otimes \mathrm{H}_j(\Lambda M) \arrow[r, "  \Phi'_*\otimes \, \mathrm{id} ", "\cong"']
        \arrow[]{d}{*_r}
        & [2.5em]
        \mathrm{H}_i(P_g M)\otimes \mathrm{H}_j(\Lambda M)
        \arrow[]{d}{*_r}
        \\
        \mathrm{H}_{i+j-n}(P_f M) \arrow[]{r}{\Phi'_* }
        &
        \mathrm{H}_{i+j-n}(P_g M) 
    \end{tikzcd}
    $$
\end{prop}

We now note what happens in the case that $f\colon M\to M$ is an involution which is homotopic to the identity.
Using Remark \ref{remark_id_just_cs} and Corollary \ref{cor_involution_homotopic_to_id_just_free_loop_space} we obtain the following corollary.

\begin{cor}\label{cor_operations_for_map_homotopic_to_id}
    Let $M$ be a closed oriented manifold and $f\colon M\to M$ an involution which is homotopic to the identity.
    Take homology with coefficients in a commutative, unital ring $R$.
    Then the pairing with respect to $f$ as well as the left- and the right-module structure of $\mathrm{H}_{\bullet}(P_f M)$ over the Chas-Sullivan ring of $M$ are completely determined by the Chas-Sullivan ring.
    More precisely, if $\Phi,\Phi'\colon P_f M\to \Lambda M$ are the homotopy equivalences as above then we have
    $$    A \wedge_f B =  \Phi_* (A) \wedge_{\mathrm{CS}} \Phi'_* (B) ,\qquad \Phi_*(X *_l A )  =  X \wedge_{\mathrm{CS}} \Phi_*(A) $$ 
    and 
    $$     \Phi'_*(A *_r X) =   \Phi'_* (A) \wedge_{\mathrm{CS}} X    \quad \text{for all}\,\,\, A,B\in\mathrm{H}_{\bullet}(P_f M) \,\,\text{and}\,\, X\in\mathrm{H}_{\bullet}(\Lambda M) .$$
\end{cor}
\begin{cor}\label{cor_trivial_extension}
    Let $M$ be a closed oriented manifold and $f\colon M\to M$ an involution which is homotopic to the identity.    
    Take homology with coefficients in a commutative, unital ring $R$.
    Assume that the maps $\Phi_*,\Phi'_*\colon \mathrm{H}_{\bullet}(P_f M)\to \mathrm{H}_{\bullet}(\Lambda M)$ agree in homology.
    Then the algebra induced by the extended Chas-Sullivan product is isomorphic to the trivial extension of the Chas-Sullivan product on $\mathrm{H}_{\bullet}(\Lambda M)$.
\end{cor}

\subsection{The extended loop product on odd-dimensional spheres}

We now consider the situation that $M= \mathbb{S}^n$ and $f = a \colon \mathbb{S}^n\to \mathbb{S}^n$ is the antipodal map. 

\begin{definition}
    Let $n\geq 1$ be an integer and consider the manifold $\mathbb{S}^n$.
    The space
    $$  P_a \mathbb{S}^n = \{   \gamma\in P\mathbb{S}^n\,|\, \gamma(1) = - \gamma(0) \}.      $$
    is called the \textit{space of antipodal paths} on $\mathbb{S}^n$.
\end{definition}
We want to compute the extended loop product on 
$$   \widehat{\mathrm{H}}_{\bullet}(\mathbb{S}^n)  = \mathrm{H}_{\bullet}(\Lambda \mathbb{S}^n) \oplus \mathrm{H}_{\bullet}(P_a \mathbb{S}^n) .  $$

\begin{prop}\label{prop_homotopy_equivalence_n_odd}
    Let $n$ be odd.
    The space of antipodal paths $P_a\mathbb{S}^n$ is homotopy equivalent to the free loop space $\Lambda \mathbb{S}^n$.
\end{prop}
\begin{proof}
   This follows from Corollary \ref{cor_involution_homotopic_to_id_just_free_loop_space} since the antipodal map on an odd-dimensional sphere is homotopic to the identity.
\end{proof}
\begin{remark}\label{remark_homotopy_antipodal_map_identity}
    Let $n$ be odd.
    Note that one can construct a homotopy between the antipodal map $a\colon \mathbb{S}^n\to\mathbb{S}^n$ and the identity explicitly as follows.
    Since on an odd-dimensional sphere there are nowhere vanishing vector fields we can take a unit vector field $X$ on $\mathbb{S}^n$.
    Let $\varphi\colon \mathbb{S}^n\to P_a\mathbb{S}^n$ be the map defined as $\varphi(p) = \gamma_{\pi \cdot X(p)}$ where $\gamma_w\colon I\to \mathbb{S}^n$, $w\in T\mathbb{S}^n$ is the geodesic starting in direction $w$ with respect to the standard metric on $\mathbb{S}^n$.
    Then a homotopy $H\colon \mathbb{S}^n\times I\to \mathbb{S}^n$ between the antipodal map and the identity is given by
    $    H(p,s) =    \varphi(p)(1-s)$ for $p\in\mathbb{S}^n, s\in I$.
\end{remark}

\begin{theorem}\label{theorem_extended_product_circle}
    Equip $\widehat{\mathbb{H}}_{\bullet}(\mathbb{S}^1)$ with the extended loop product where homology is taken with integer coefficients.
    There is an isomorphism of graded algebras
    $$      \widehat{\mathbb{H}}_{\bullet}(\mathbb{S}^1) \cong \Lambda[a] \otimes \mathbb{Z}[t,t^{-1}]      $$
    where $|a| = -1$ and $|t| = |t^{-1}|= 0$.
    The Chas-Sullivan algebra of $\mathbb{S}^1$ is the subalgebra $\Lambda[a]\otimes \mathbb{Z}[t^2,t^{-2}]$ and the extended Chas-Sullivan algebra is not a trivial extension.
\end{theorem}
\begin{proof}
    Note that we have
    $$   P_a\mathbb{S}^1 \simeq \Lambda \mathbb{S}^1 \cong \mathbb{S}^1 \times \Omega \mathbb{S}^1 \simeq \mathbb{S}^1 \times \mathbb{Z}         $$
    where the homeomorphism $\Lambda \mathbb{S}^1\cong \mathbb{S}^1\times \Omega \mathbb{S}^1$ holds since $\mathbb{S}^1$ is a Lie group.
    Moreover, it is well-known that $\Omega \mathbb{S}^1\simeq \mathbb{Z}$ since each connected component of $\Omega\mathbb{S}^1$ is contractible.
    Let $k =  k_0 + \tfrac{1}{2}$, $k_0\in\mathbb{Z}$.
    If $k> 0$, i.e. $k_0\geq 0$ denote by $\gamma_k\in P_a\mathbb{S}^1$ the antipodal path parametrized at constant speed starting at $p_0 = 1$ in the positively oriented direction wrapping $k_0$ times around the circle and the going to $-p_0$.
    For $k<0$ let $\gamma_k$ be the inverse path of $\gamma_{-k}$.
    In particular, $\gamma_{\tfrac{1}{2}}$ is the length-minimizing geodesic connecting $1$ to $-1$ which goes along the upper part of the circle and $\gamma_{-\tfrac{1}{2}}$ is the length-minimizing geodesic going along the lower part.
    Moreover, let $\varphi\colon \mathbb{S}^1\times P_a\mathbb{S}^1\to P_a\mathbb{S}^1$ be the action of $\mathbb{S}^1$ induced by the group structure of $\mathbb{S}^1$, i.e.
    $$   \varphi(s,\gamma)(t) =   \gamma(t) + s \quad \text{for}\,\,\,s\in\mathbb{S}^1,\,\,\gamma\in P_a\mathbb{S}^1,\,\,t\in [0,1].   $$
    A generating set for the homology of $P_a\mathbb{S}^1$ is given by the classes
    $$        \{ [\gamma_{k}] \in \mathrm{H}_0(P_a\mathbb{S}^1) \,|\, k\in \mathbb{Z} + \tfrac{1}{2} \} \cup  \{  \varphi_* ( [\mathbb{S}^1] \times [\gamma_k] )\in \mathrm{H}_1(P_a\mathbb{S}^1) \,|\,  k\in \mathbb{Z} + \tfrac{1}{2}   \}.           $$
    Denote the classes $  \varphi_* ( [\mathbb{S}^1] \times [\gamma_k] )$ by $\Gamma_k$.
    Similarly, for the free loop space of $\mathbb{S}^1$ we note that a generating set is given by 
    $$        \{ [\gamma_{k}']\in \mathrm{H}_0(\Lambda \mathbb{S}^1) \,|\, k\in \mathbb{Z}  \} \cup  \{  \varphi'_* ( [\mathbb{S}^1] \times [\gamma_k]' )\in\mathrm{H}_1(\Lambda \mathbb{S}^1) \,|\,  k\in \mathbb{Z}   \}.           $$
   Here, $\gamma_k'$, $k\in\mathbb{Z}$ is the loop winding $k$ times around the circle and $\varphi'\colon \mathbb{S}^1\times \Lambda \mathbb{S}^1\to \Lambda \mathbb{S}^1$ is the $\mathbb{S}^1$-action on the free loop space of the circle induced by the multiplication of $\mathbb{S}^1$.
    We define $\Gamma_k' = \varphi'_* ( [\mathbb{S}^1] \times [\gamma_k'] )\in\mathrm{H}_1(\Lambda\mathbb{S}^1)$ for $k\in\mathbb{Z}$.
    It is well-known, see \cite{cohen:2003}, that for the Chas-Sullivan product on $\mathbb{S}^1$ we have
    \begin{equation}\label{eq_cs_product_circle}
           [\gamma'_k]\wedge_{\mathrm{CS}}[\gamma'_l] = 0,\quad \Gamma'_k \wedge_{\mathrm{CS}} [\gamma'_l]  =  [\gamma'_k]\wedge_{\mathrm{CS}} \Gamma'_l  = \gamma'_{k+l} \quad \text{and}\quad \Gamma'_k \wedge_{\mathrm{CS}} \Gamma'_l = \Gamma'_{k+l}  
    \end{equation}
    for $k,l\in\mathbb{Z}$.
     Define $\eta\colon \mathbb{S}^1\to P_a\mathbb{S}^1$ by taking $\eta(p)$ to be the injective antipodal path going in the positive direction connecting $p$ and $-p$ for $p\in\mathbb{S}^1$.
    This induces the homotopy equivalences $\Phi\colon P_a\mathbb{S}^1\to \Lambda\mathbb{S}^1$ and $\Phi'\colon P_a\mathbb{S}^1\to \Lambda \mathbb{S}^1$
    as in Subsection \ref{subsec_invariance_loop_pairings}.    
    One checks that we have
    $$  \Phi_* [\gamma_k] =  [\gamma'_{k+\tfrac{1}{2}}] \quad \text{while}\quad \Phi'_* [\gamma_k] =  [\gamma'_{k-\tfrac{1}{2}}]      $$
    for $k\in\mathbb{Z}+\tfrac{1}{2}$.
    The analogous expressions hold for the classes of the form $\Gamma_k$, $k\in\mathbb{Z} + \tfrac{1}{2}$.    
   Then we see by the identities in \eqref{eq_cs_product_circle} and Corollary \ref{cor_operations_for_map_homotopic_to_id} that
    $$      [\gamma_k]\wedge_a [\gamma_l] = 0,\quad \Gamma_k\wedge_a [\gamma_l] =  [\gamma'_{k+l}] , \quad [\gamma_k] \wedge_a [\Gamma_l] = [\gamma'_{k+l}]      \quad \text{and}\quad  \Gamma_k \wedge_a \Gamma_l =  \Gamma'_{k+l}      $$
    for $k,l\in \mathbb{Z} + \tfrac{1}{2}$.
     One can also compute the module structure of $\mathrm{H}_{\bullet}(P_a\mathbb{S}^1)$ without problems and finds that
     $$     [\gamma'_k]*_l [\gamma_l]= [\gamma_l]*_r [\gamma'_k] =   0 , \quad \Gamma'_k *_l  [\gamma_l] = [\gamma_l]*_r \Gamma'_k  = [\gamma_{k+l}]     $$
    as well as
    $$   [\gamma'_k] *_l  \Gamma_l =  \Gamma_l *_r [\gamma'_k]   =[\gamma_{k+l}]    \quad \text{and}\quad \Gamma'_k *_l \Gamma_l = \Gamma_l *_r \Gamma'_k =  \Gamma_{k+l}  $$
    for $k\in\mathbb{Z}$, $l\in\mathbb{Z}+\tfrac{1}{2}$.
     We see that the algebra structure on the direct sum $\mathrm{H}_{\bullet}(\Lambda \mathbb{S}^1) \oplus \mathrm{H}_{\bullet}(P_a\mathbb{S}^1)$ behaves as follows.
     If we set $t = \Gamma_{\tfrac{1}{2}}$ and $t^{-1}= \Gamma_{-\tfrac{1}{2}}$ as well as $a = [\gamma_0]$ then there is an isomorphism of graded algebras
    $$      \widehat{\mathbb{H}}_{\bullet}(\mathbb{S}^1) \cong \Lambda[a] \otimes \mathbb{Z}[t,t^{-1}]      $$
    where $|a| = -1$ and $|t| = |t^{-1}|= 0$.
    Moreover, the Chas-Sullivan algebra of $\mathbb{S}^1$ is the subalgebra $\Lambda[a]\otimes \mathbb{Z}[t^2,t^{-2}]$.
    Finally, we see that this is extension algebra not a trivial extension since there is no element $X\in\mathrm{H}_{\bullet}(P_a\mathbb{S}^1)$ such that $X^2 = 1 = \Gamma_0'$, see Remark \ref{example_trivial_extension}.  
\end{proof}

We end this section by showing that for some odd-dimensional spheres the maps $\Phi$ and $\Phi'$ can be chosen in such a way that they induce the same map in homology.
This implies that the extension of the Chas-Sullivan algebra induced by the antipodal pairing is a trivial extension by Corollary \ref{cor_trivial_extension}.
\begin{lemma}
    Let $\mathbb{S}^n$ be an odd-dimensional sphere admitting at least two linearly independent nowhere vanishing vector fields. Then the maps $\Phi\colon P_a\mathbb{S}^n\to \Lambda \mathbb{S}^n$ and $\Phi'\colon P_a\mathbb{S}^n\to \Lambda \mathbb{S}^n$ are homotopic.
\end{lemma}
\begin{proof}
    We recall the construction of $\Phi$ and $\Phi'$ as in the proof of Proposition \ref{prop_homotopy_invariance_involution_path_spaces} and in the discussion before Proposition \ref{prop_pairing__for_homotopic_maps}.
    Let $X$ be a nowhere vanishing vector field on $\mathbb{S}^n$ such that $X(-p) = - X(p)$ under the identification $T_p\mathbb{S}^n \cong T_{-p}\mathbb{S}^n$.
    Such a vector field exists since if we write $n = 2m-1$ we can take e.g. 
    $$   X(x_1,x_2,\ldots, x_{2m}) =  (x_2,-x_1,\ldots, x_{2m},-x_{2m-1}) \quad \text{for}\,\,\, (x_1,\ldots, x_{2m}) \in \mathbb{S}^{2m-1}   .  $$
    By rescaling we can achieve that $|X(p)| = \pi$ for all $p\in\mathbb{S}^n$.
    Let $\gamma_v\colon I\to \mathbb{S}^n$ be the geodesic with initial condition $\Dot{\gamma}_v(0) = v$.
    Define the map $\eta\colon \mathbb{S}^n\to P\mathbb{S}^n$ by setting $\eta(p) = \gamma_{X(-p)}$.
    We have
    $$  \Phi(\sigma) = \mathrm{concat}(\sigma, \gamma_{X(-\sigma(0))}) \quad \text{for}\,\,\, \sigma\in P_a\mathbb{S}^n .   $$
    Furthermore, we have
    $$   \Phi'(\sigma) =   \mathrm{concat}(\overline{\eta(\sigma(1))},\sigma)   = \mathrm{concat}(\overline{\gamma_{X(\sigma(0))}},\sigma)    \quad \text{for}\,\,\,\sigma\in P_a\mathbb{S}^n     $$
    and it is clear that $\Phi'$ is homotopic to the map $\Phi''\colon P_a\mathbb{S}^n\to \Lambda \mathbb{S}^n$ given by
    $$   \Phi'' (\sigma)  =    \mathrm{concat}(\sigma,\overline{\eta(\sigma(1))}) =  \mathrm{concat}(\sigma,\overline{\gamma_{X(\sigma(0))}}) \quad \text{for}\,\,\,\sigma\in P_a\mathbb{S}^n  .       $$
    Note that the geodesic $\gamma_{X(p)}$ satisfies
    $   \Dot{\gamma}_{X(p)}(1) = - X(p)       $
    under the canonical identification of $T_p \mathbb{S}^n \cong (p)^{\perp}$ with $T_{-p}\mathbb{S}^n \cong (-p)^{\perp} \cong (p)^{\perp}$.
    Hence, the path $\overline{\gamma_{X(p)}}$ is equal to the geodesic $\gamma_{X(p)|_{-p}}$ where $X(p)|_{-p}$ means that we consider the vector $X(p)$ as a tangent vector to $-p$.
    Consequently, we have
    $$   \Phi''(\sigma) =  \mathrm{concat}(\sigma, \gamma_{X(\sigma(0))|_{-\sigma(0)}})  \quad \text{for}\,\,\,\sigma\in P_a\mathbb{S}^n .      $$
    By using a second linearly independent vector field we can now rotate the antipodal path $\gamma_{X(\sigma(0))|_{-\sigma(0)}}$ into the path $\gamma_{X(-\sigma(0))}$.
    More precisely, assume that $Y$ is a second vector field on $\mathbb{S}^n$ which is nowhere vanishing and linearly independent to $X$.
    We can assume without loss of generality that $Y$ is orthogonal to $X$ at each point and that $|Y| = \pi$.
    Then a homotopy between the maps $\Phi$ and $\Phi''$ is given by defining $H\colon P_a\mathbb{S}^n\times I \to \Lambda \mathbb{S}^n$ as
    $$    H(\sigma,s) =    \mathrm{concat}(\sigma, \gamma_{\mathrm{cos}(\pi s) X(-\sigma(0)) + \mathrm{sin}(\pi s)Y(-\sigma(0))})  \quad \text{for}\,\,\,\sigma\in P_a\mathbb{S}^n, \,\, s\in I.       $$
    One checks that $H(\sigma,0) = \Phi''(\sigma)$ and $H(\sigma,1) = \Phi(\sigma)$.
    Hence $\Phi\simeq \Phi''\simeq \Phi'$ as claimed.
\end{proof}

The next Corollary is then immediate from Proposition \ref{prop_pairing__for_homotopic_maps}.
\begin{cor}
    Let $\mathbb{S}^n$ be an odd-dimensional sphere admitting at least two linearly independent vector fields.
    Then we have the identity
    $$    A\wedge_a B =  \Phi_*(A) \wedge_{\mathrm{CS}}  \Phi_*(B)   \quad \text{for}\,\,\, A,B\in\mathrm{H}_{\bullet}(P_a\mathbb{S}^n)  .$$
    Moreover, for the module multiplications we have
    $$    \Phi_*( X*_l A) = X\wedge_{\mathrm{CS}} \Phi_*(A) \quad \text{and} \quad \Phi_*( A*_r X) =  \Phi_*(A) \wedge_{\mathrm{CS}} X         $$
    for $A\in\mathrm{H}_{\bullet}(P_a\mathbb{S}^n)$ and $X\in\mathrm{H}_{\bullet}(\Lambda \mathbb{S}^n)$.
    Hence, the extension algebra on $\mathrm{H}_{\bullet}(\Lambda \mathbb{S}^n) \oplus \mathrm{H}_{\bullet}(P_a\mathbb{S}^n)$ induced by the antipodal pairing is a trivial extension in the sense of Remark \ref{example_trivial_extension}.
\end{cor}
We recall that the Chas-Sullivan algebra of an odd-dimensional sphere satisfies 
$$  \big(\mathbb{H}_{\bullet}(\mathbb{S}^n),\wedge_{\mathrm{CS}}\big) \cong \frac{\mathbb{Z}[A,U]}{(A^2)} \quad \text{with}\,\,\, |A| = -n,\,\,|U| = n-1 , $$
see \cite{cohen:2003}.
A straight-forward computation then shows the following.
\begin{cor}\label{cor_extendedn_product_sphere}
    Let $\mathbb{S}^n$ be an odd-dimensional sphere admitting at least two linearly independent vector fields.
    Take homology with integer coefficients.
    Then the extended Chas-Sullivan algebra is isomorphic to the polynomial algebra
    $$   \widehat{\mathbb{H}}_{\bullet}(\mathbb{S}^n) \cong \frac{\mathbb{Z}[A,U,{E},B,V]}{\big(  A^2, B^2, BV - AU, V^2 - U^2, {E}B = A, {E}V = U, {E}^2 -1    \big)}     $$
    where $|A|= |B| = -n$, $|{E}| = 0$ and $|U| = |V| = n-1$. 
    The Chas-Sullivan algebra of $\mathbb{s}^n$ is the subalgebra generated by the classes $A$ and $U$.
\end{cor}
\begin{remark}\label{remark_non-trivial_subalgebra}
    Consider the situation of Corollary \ref{cor_extendedn_product_sphere}.
    \begin{enumerate}
        \item 
    One checks that the subalgebra generated by the classes $B$ and $V$ is also isomorphic to the Chas-Sullivan algebra of an odd-dimensional sphere, i.e.
    $$   ( B,V )  \cong \frac{\mathbb{Z}[A',U']}{((A')^2)}   \quad \text{with}\,\,\, |A'| = -n,\,\, |U'| = n-1 .  $$
    \item     It is well-known that an odd-dimensional sphere $\mathbb{S}^n$ admits at least two linearly independent vector fields if and only if $n = 3$ mod $4$.
    In the case $n =1$ we have seen above that the extended Chas-Sullivan algebra is not a trivial extension of the Chas-Sullivan algebra.
    It would be interesting to understand the cases $n = 1$ mod $4$ for $n\geq 5$.
    \end{enumerate}
\end{remark}

\section{The energy functional on the space of antipodal paths}\label{sec_antipodal_geodesics}

In this section we consider the space of antipodal paths on the sphere.
On this space we study the critical points of the energy functional determined by a Riemannian metric. 
In particular we compute the Morse indices of the critical points induced by the standard metric.

If $M$ is a closed manifold and $g$ a Riemannian metric on $M$ then there is an energy functional $E\colon PM\to [0,\infty)$ which is defined as
$$   E(\gamma) = \int_0^1 g_{\gamma(t)}(\Dot{\gamma }(t), \Dot{\gamma}(t))\,\mathrm{d}t \quad \text{for}\,\,\,\gamma\in PM .  $$
We will consider the restriction of $E$ to submanifolds of $PM$ but shall abuse notation and always write $E$ for these restrictions as well.
It turns out that the critical points of the energy functional on the free loop space $\Lambda M$ are the \textit{closed geodesics} in $M$, i.e. geodesics $\gamma\colon I\to M$ with
$$   \gamma(0) = \gamma(1) \quad \text{and}\quad \Dot{\gamma}(0) = \Dot{\gamma}(1)  .    $$
We shall now follow the well-known story of the closed geodesics on the free loop space,, see \cite{klingenberg:78} and \cite{klingenberg:1995}, and consider the critical points of the energy functional on the space of antipodal paths on the sphere.
For the following definition recall that there is a canonical identification of the tangent spaces 
$    T_p \mathbb{S}^n  \cong T_{-p}\mathbb{S}^n$ for all $p\in\mathbb{S}^n$.
This identification is given by the standard embedding of the sphere $\mathbb{S}^n$ into $\mathbb{R}^{n+1}$ and the observation that $T_p \mathbb{S}^n $ and $T_{-p}\mathbb{S}^n $ are the same subspaces of $\mathbb{R}^{n+1}$.

\begin{definition}\label{def_antipodal_geodesics}
    Let $n\geq 1 $ be an integer and consider the $n$-sphere $\mathbb{S}^n$.
    An \textit{antipodal geodesic} is a geodesic $\gamma\colon I\to \mathbb{S}^n$ such that
    $$   \gamma(1) = -\gamma(0) \quad \text{and}\quad \Dot{\gamma}(1)  = -\Dot{\gamma}(0)     $$
    under the identification of $T_p \mathbb{S}^n $ and $T_{-p}\mathbb{S}^n $ explained above.
\end{definition}

In the following we will show that the antipodal geodesics are precisely the critical points of the energy functional on $P_a\mathbb{S}^n$.
First, we need to study the tangent spaces of $P_a\mathbb{S}^n$.
We now recall from \cite[p. 167]{klingenberg:78} that the tangent space to an element $\gamma\in PM$ is the space of $H^1$-vector fields along $\gamma$, see \cite{klingenberg:78} for details.
Moreover, by \cite[Proposition 2.4.1]{klingenberg:78} we see directly that the tangent space of $P_a\mathbb{S}^n$ behaves as follows.
For $\gamma\in P_a\mathbb{S}^n$ we have
$$    T_{\gamma} P_a\mathbb{S}^n  \cong \{ \xi \,\,\,H^1\text{vector field along} \,\,   \gamma\,|\,   \xi(0) = -\xi(1) \}  $$
where we again use the identification $T_p\mathbb{S}^n \cong T_{-p}\mathbb{S}^n$.
As in \cite[Lemma 2.4.3]{klingenberg:78} we get the following characterization of the critical points of the energy functional on $P_a\mathbb{S}^n$.
\begin{lemma}
    Let $n\geq 1$ and consider the space of antipodal paths $P_a\mathbb{S}^n$.
    The critical points of the energy functional on $P_a\mathbb{S}^n$ are precisely the antipodal geodesics.
\end{lemma}
\begin{remark}
    Consider the space of antipodal paths on the sphere $\mathbb{S}^n$.
    \begin{enumerate}
        \item Note that not every geodesic $\gamma\colon I\to \mathbb{S}^n$ with $\gamma(1) = -\gamma(0)$ is a critical point of the energy functional.
        This is analogous to the fact that not every geodesic $\gamma\colon I\to \mathbb{S}^n$ with $\gamma(1) = \gamma(0)$ is a closed geodesic.
        \item On the free loop space $\Lambda M$ of a closed manifold we have that the trivial loops are critical points of the energy functional.
        Clearly these are precisely the loops with minimal energy.
        Hence, in this case we understand the global minima of the energy functional very well.
        In contrast, on the space of antipodal paths $P_a\mathbb{S}^n$, it depends on the particular metric how the antipodal paths of minimal energy behave.
    \end{enumerate}
\end{remark}

We now turn to the standard metric on the sphere.
It is well-known how the geodesics behave for the standard metric.
In particular we observe the following.
Let $\gamma\colon [0,\infty)\to \mathbb{S}^n$ be a unit speed geodesic starting at the point $\gamma(0) = p$.
Then at times $t =   (2k+1) \pi, \quad k\in \mathbb{N}_0   $ we have $\gamma(t) = -p$ and $\Dot{\gamma}(t) = -\Dot{\gamma}(0)$.
The antipodal geodesics thus do not come as isolated critical points, but the critical manifolds are diffeomorphic to the unit tangent bundle $U\mathbb{S}^n$.
We say that an antipodal geodesic $\gamma\in P_a\mathbb{S}^n$ is \textit{prime} if the map $\gamma\colon I\to \mathbb{S}^n$ is injective.
Clearly, the prime antipodal geodesics are the global minimum of the energy functional.
We say that an antipodal geodesic has multiplicity $i\geq 1$ if it is the concatenation of a closed geodesic of multiplicity $i$ with a prime antipodal geodesic.
The prime antipodal geodesics are said to have multiplicity $0$.
We summarize the behavior of the antipodal geodesics in the following Lemma.
\begin{prop}\label{prop_cr_submflds}
    Consider the $n$-sphere with the standard metric.
    The critical sets of the energy functional on $P_a\mathbb{S}^n$ are the manifolds
    $$      \Sigma_i  =  \{   \gamma \text{ is an antipodal geodesic of multiplicity }i \} ,\quad i\geq 0.         $$
    Moreover, the critical values are
    $    E(\Sigma_i ) =  {((2i+1)\pi)^2}, \quad i\geq 0     $
    with ${\pi^2}$ being the global minimum.
    Each critical manifold $\Sigma_i$ is diffeomorphic to the unit tangent bundle with the diffeomorphism given by
    $$   U\mathbb{S}^n \to \Sigma_i,\quad v \mapsto \gamma_{(2i+1)\pi v} ,\quad i\geq 0   $$
    where $\gamma_w\colon [0,1]\to\mathbb{S}^n$ is the geodesic on $\mathbb{S}^n$ starting in direction $w\in T\mathbb{S}^n$. 
\end{prop}
Grove shows in \cite{grove1973condition} that the Palais-Smale condition, also called condition (C), holds for the energy functional on certain path spaces including the space of antipodal paths on the sphere.

We shall now study how the Morse index behaves for the energy functional of the standard metric. 
We will eventually see that the energy functional is a Morse-Bott function.

Let $X$ be a Hilbert manifold and $f\in C^2(X)$ a twice continuously differentiable function satisfying condition (C).
If $x\in X$ is a critical point, i.e. $df_x = 0$, then we consider the Hessian $D^2f_x \colon T_x X\times T_x X\to \mathbb{R}$.
The \textit{index} of $x$, denoted by $\mathrm{ind}(x)$ is defined as the dimension of a maximal subspace of $T_x X$ on which the Hessian is negative definite.
Similarly, the \textit{nullity} of a critical point is defined as the dimension of the kernel of the Hessian.
If the set of critical points is a disjoint union which consists of connected finite-dimensional submanifolds such that for each critical point $x\in X$ that lies in the critical submanifold $\Sigma$ we have $\mathrm{null}(x) = \mathrm{dim}(\Sigma)$, then $f$ is called a \textit{Morse-Bott function}.

Assume that the Riemannian metric on $\mathbb{S}^n$ is invariant under the antipodal map. In this case we can related the index and nullity of the energy functional on $P_a\mathbb{S}^n$ to the corresponding quantities on $\Lambda \mathbb{R}P^n$ as we show in the following Lemma.
\begin{lemma}\label{lemma_index_and_null_preserved_under_covering}
    Let $g$ be a $\mathbb{Z}_2$-invariant metric on $\mathbb{S}^n$ and let $g'$ be the induced Riemannian metric on $\mathbb{R}P^n$.
    Consider the energy functionals $E\colon P_a\mathbb{S}^n\to [0,\infty)$ induced by $g$ and $E'\colon \Lambda \mathbb{R}P^n\to [0,\infty)$ induced by $g'$.
    If $\gamma\in P_a\mathbb{S}^n$ is a critical point of $E$, then $\gamma' = \pi\circ \gamma$ is a critical point of $E'$ where $\pi\colon \mathbb{S}^n\to \mathbb{R}P^n$ is the quotient map.
    Then we have $\mathrm{ind}(\gamma) =\mathrm{ind}(\gamma') $ as well as $\mathrm{null}(\gamma) = \mathrm{null}(\gamma')$.
\end{lemma}
\begin{proof}
    Recall from Proposition \ref{prop_involution_path_space_is_universal_covering} that the map $p\colon P_a\mathbb{S}^n\to \Lambda_1\mathbb{R}P^n, \gamma\mapsto \pi\circ \gamma$ is a covering, where $\Lambda_1\mathbb{R}P^n \subseteq \Lambda \mathbb{R}P^n$ is the component of non-contractible loops.
    Moreover by construction we have $E = E'\circ p$ and it is clear that critical points of $E$ are mapped to critical points of $E'$.
    Since $p$ is a covering, it is locally trivial.
    The index and the nullity of a critical point are determined by the local behavior of the function and hence if $\gamma\in P_a\mathbb{S}^n$ is a critical point of $E$, then its index and nullity agree with the index and nullity of $p(\gamma) = \pi\circ \gamma\in \Lambda \mathbb{R}P^n$.
\end{proof}
Since the standard metric on $\mathbb{S}^n$ is $\mathbb{Z}_2$-invariant, we can use the above Lemma and the results on index and nullity of the compact rank one symmetric spaces in \cite{ziller:1977} to see that the following holds.
\begin{lemma}\label{lemma_index-and-nullity}
    Consider the standard metric on $\mathbb{S}^n$.
    Then for an antipodal geodesic $\gamma$ of multiplicity $k$ we have
    $$   \mathrm{ind}(\gamma) = 2k(n-1) \quad\text{and}\quad \mathrm{null}(\gamma) = 2n-1 .      $$
\end{lemma}

\begin{theorem}
    Consider the sphere $\mathbb{S}^n$ with the standard metric. 
    The induced energy functional on the space of antipodal paths $P_a\mathbb{S}^n$ is a Morse-Bott function.
\end{theorem}
\begin{proof}
    As mentioned above, the energy functional on the space of antipodal paths satisfies condition (C) of Palais-Smale, see \cite{grove1973condition}.
    We have seen in Proposition \ref{prop_cr_submflds} that the critical submanifolds are each of dimension $2n-1$ which agrees with the nullity of each critical point by Lemma \ref{lemma_index-and-nullity}.
    Consequently, the energy functional on $P_a\mathbb{S}^n$ is a Morse-Bott function.
\end{proof}

    In the next section we shall even see that the energy functional for the standard metric is a \textit{perfect} Morse-Bott function.
\begin{remark}
    Consider the space of antipodal paths $P_a\mathbb{S}^n$ equipped with the energy functional induced by the standard metric.
   Note that by Lemma \ref{lemma_index-and-nullity} the index of an antipodal geodesic $\gamma$ equals the number of conjugate points in its interior counted with multiplicity.
        This is analogous to the situation of the energy functional on the free loop space on the standard sphere, see \cite{oancea:2015}.
        Moreover, it shows that the Morse index of an antipodal geodesic $\gamma\in P_a\mathbb{S}^n$ as a critical point in $P_a\mathbb{S}^n$ agrees with the Morse index on the space $\Omega_{p,-p}\mathbb{S}^n$ of paths
        $$  \Omega_{p,-p}\mathbb{S}^n= \{\gamma\in P\mathbb{S}^n\,|\, \gamma(0) = p = -\gamma(1) \}   .   $$
        This is because the Morse index in a path space of the form $\Omega_{pq}\mathbb{S}^n$ equals the sum of the conjugate points along the geodesic counted with multiplicity by \cite[Theorem 2.5.9]{klingenberg:1995}.
\end{remark}

\section{Morse-Bott theory on the space of antipodal paths}\label{sec_completing_manifolds}

In this section we will show that the energy functional on the space of antipodal paths $P_a\mathbb{S}^n$ is a perfect Morse-Bott function.
We use the notion of completing manifolds which we will introduce in the first subsection.
We then define explicit completing manifolds on the space of antipodal paths and on the free loop space of the sphere and determine the homology of $P_a\mathbb{S}^n$ explicitly.

\subsection{Completing manifolds and perfect Morse-Bott functions}

Consider a Hilbert manifold $X$ and let $f \colon X\to \mathbb{R}$ be a smooth function on $X$ satisfying the Palais-Smale condition (C).
For $a\in \mathbb{R}$ we shall use the notation
$$     X^{\leq a} = \{ x\in X \,|\, f(x) \leq a \} \quad \text{and} \quad X^{<a} = \{ x\in X \,|\, f(x) < a\} .      $$
Let $a$ be a critical value of $f$ and assume that the set of critical points $\Sigma$ at level $a$ is a non-degenerate finite-dimensional critical submanifold of finite index $k$ with orientable negative bundle.
The behavior of the relative homology $\mathrm{H}_{\bullet}(X^{\leq a},X^{<a})$ is then well known, i.e. there is an isomorphism
$$   \mathrm{H}_{\bullet}(X^{\leq a},X^{<a}) \cong \mathrm{H}_{\bullet -k}(\Sigma)        $$
where we take coefficients in a commutative unital ring $R$.
In applications the homology of these critical submanifolds is usually easier to understand than the homology of $X$.
Therefore, it is desirable to find conditions which imply that the homology of $X$ is built up by respective homology groups of the critical manifolds.
\begin{definition}[\cite{oancea:2015}, Definition 6.1] \label{def_completing_mfld}
Let $X$ be a Hilbert manifold and let $f$ be a smooth real-valued function on $X$ satisfying condition (C).
Let $a$ be a critical value of $f$ and assume that $\Sigma$ is a non-degenerate connected oriented critical submanifold at level $a$ of index $k$ and of dimension $l = \mathrm{dim}(\Sigma)$.
Assume that $k$ and $l$ are both finite.
A \textit{completing manifold} for $\Sigma$ is a closed, orientable manifold $\Gamma$ of dimension $k + l$ with an embedding $\varphi: \Gamma \to X^{\leq a}$ such that the following holds.
There is an $l$-dimensional submanifold $S\subseteq \Gamma$ such that $\varphi|_S$ maps $S$ homeomorphically onto $\Sigma$.
Moreover, there is a retraction map $p:\Gamma\to S$.
Furthermore, the embedding $\varphi$ induces a map of pairs
$  \varphi \colon (\Gamma,\Gamma\setminus S) \to (X^{\leq a}, X^{<a})    $.
\end{definition}

We now recall some facts about \textit{Gysin maps} for finite-dimensional manifolds.
Assume that $f\colon M\to N$ is a map between oriented manifolds $M$ and $N$.
Then the Gysin map
$$f_!\colon \mathrm{H}_j(N)\to \mathrm{H}_{j+\mathrm{dim}(M)-\mathrm{dim}(N)}(M)$$ is defined by
$$f_! \colon \mathrm{H }_j(N) \xrightarrow[]{(PD_{N})^{-1}}  \mathrm{H}^{\mathrm{dim}(N) - j}(N) \xrightarrow[]{f^*} \mathrm{H}^{\mathrm{dim}(N)-j}(M) \xrightarrow[]{PD_{M}} \mathrm{H}_{\mathrm{dim}(M) - (\mathrm{dim}(N)-j)}(M) . $$
Here $PD_B\colon \mathrm{H}^{i}(B) \to \mathrm{H}_{\mathrm{dim}(B)-i}(B)$ stands for the Poincaré duality isomorphism on the oriented manifold $B$. 
In the above situation of a completing manifold the existence of a retraction $p\colon \Gamma \to S$ implies the following.
The Gysin map $p_!\colon \mathrm{H}_i(S)\to \mathrm{H}_{i+k}(\Gamma)$ is a right inverse to the Gysin map $s_!\colon \mathrm{H}_i(\Gamma)\to \mathrm{H}_{i-k}(S)$ where $s\colon S\hookrightarrow \Gamma$ is the inclusion.
Up to sign one can see that the Gysin map $s_!$ equals the composition
$$    \mathrm{H}_i(\Gamma) \to \mathrm{H}_i(\Gamma,\Gamma\setminus S) \xrightarrow[]{\cong}  \mathrm{H}_{i-k}(S)   $$
where the first map is induced by the inclusion of pairs and the second one being the Thom isomorphism, see \cite[Theorem VI.11.3]{bredon:2013}.
Consequently, the map $\mathrm{H}_{\bullet}(\Gamma)\to \mathrm{H}_{\bullet}(\Gamma,\Gamma\setminus S)$ is surjective.
In particular we obtain the following.

\begin{prop}[\cite{oancea:2015}, Lemma 6.2] \label{prop_completing_mfld}
Let $X$ be a Hilbert manifold, $f$ a smooth real-valued function on $X$ satisfying condition (C) and $a$ be a critical value of $f$.
Assume that the set of critical points at level $a$ is a non-degenerate critical submanifold $\Sigma$ of index $k$.
If there is a completing manifold $\Gamma$ for $\Sigma$ then
$$   \mathrm{H}_{\bullet}(X^{\leq a}) \cong \mathrm{H}_{\bullet}(X^{<a}) \oplus \mathrm{H}_{\bullet }(X^{\leq a},X^{<a}) \cong \mathrm{H}_{\bullet}(X^{<a}) \oplus \mathrm{H}_{\bullet -k}(\Sigma)    .      $$
The injection $\mathrm{H}_{\bullet -k}(\Sigma) \hookrightarrow \mathrm{H}_{\bullet}(X^{\leq a})$ is given by the composition $\varphi_* \circ p_!$ where $\varphi\colon \Gamma\to X^{\leq a}$ is the embedding of the completing manifold and $p\colon \Gamma\to S\cong\Sigma$ is the retraction.
\end{prop}

Now, take homology with coefficients in a commutative ring $R$.
If the long exact sequence of the pair $(X^{\leq a}, X^{<a})$ has vanishing connecting homomorphisms for each $i\geq 0$, then we obtain short exact sequences
$$   0\xrightarrow[]{} \mathrm{H}_i(X^{<a}) \xrightarrow[]{} \mathrm{H}_i(X^{\leq a},X^{<a}) \xrightarrow[]{} \mathrm{H}_i(X^{\leq a}, X^{<a}) \xrightarrow[]{} 0 .   $$
If these are split for all critical values $a$ and for all $i\geq 0$, we say that the function $f$ is a \textit{perfect} Morse-Bott function for $R$-coefficients.
If every critical submanifold has a completing manifold, we have seen above that the function $f$ is perfect with respect to arbitrary coefficients.
Using completing manifolds Ziller shows in \cite{ziller:1977} that the energy function on the free loop space of a compact symmetric space is a perfect Morse-Bott function for $\mathbb{Z}_2$-coefficients.

Assume that $f\colon X\to \mathbb{R}$ is a Morse-Bott function on a Hilbert manifold $X$ which is bounded from below such that each critical manifold admits a completing manifold.
Then we have an isomorphism
$$    \mathrm{H}_{\bullet}(X) \cong   \bigoplus_{\Sigma \text{ cr. mfld.}} \mathrm{H}_{\bullet - \lambda_{\Sigma}} (\Sigma)     $$
where the direct sum runs over all critical manifolds $\Sigma$ and where $\lambda_{\Sigma}$ is the index of $\Sigma$.
Moreover, as seen above we do not just get this abstract isomorphism but we get explicit injections
$   \mathrm{H}_i(\Sigma) \hookrightarrow \mathrm{H}_{i+\lambda_{\Sigma}}(X)    $
given by the composition
$$     \mathrm{H}_i(\Sigma) \cong \mathrm{H}_i(S) \xrightarrow[]{p_!} \mathrm{H}_{i+\lambda_{\Sigma}} (\Gamma) \xrightarrow[]{\varphi_*} \mathrm{H}_{i+\lambda_{\Sigma}} (X)  .   $$
Thus, we obtain a very explicit way of writing down generators of the homology of $X$.

\subsection{Completing manifolds on the space of antipodal paths}

We shall now see that we can indeed find completing manifolds for the critical manifolds in the space of antipodal paths as well as for the free loop space of a sphere.
On the free loop space these completing manifolds agree with Ziller's completing manifolds in \cite{ziller:1977}.
Note that Ziller treats general compact symmetric spaces, therefore he defines these manifolds as certain quotients of a product of groups.
In our case it is easier to define the completing manifolds as fiber products of the unit tangent bundle of the sphere.

For $l\geq 0$ define $\Gamma_l$ to be the fiber product
$$    \underbrace{U\mathbb{S}^n\times_{\mathbb{S}^n} U\mathbb{S}^n\times_{\mathbb{S}^n} \times \ldots \times_{\mathbb{S}^n} U\mathbb{S}^n}_{l+1 \text{ times}} .    $$
In particular we have that $\Gamma_0 = U\mathbb{S}^n$.
Note that the projection onto the first factor yields submersions $p_l\colon \Gamma_l \to U\mathbb{S}^n$.
Moreover, there is a section of $p_l$ given by the diagonal, i.e. the map $s_l\colon U\mathbb{S}^n\to \Gamma_l$ given by
$$   s_l( x) =   (x,\ldots, x) \quad \text{for}\,\,\,x\in U\mathbb{S}^n    $$
is a section of $p_l$.
Note that for a point $x\in U\mathbb{S}^n$ we shall also write $x = (p,v)$ with $p\in\mathbb{S}^n$ the basepoint and $v\in T_p\mathbb{S}^n$.
An element of $\Gamma_l$ can thus also be written as 
$    (p,v_0,\ldots, v_l) \in \Gamma_l   $ with $(p,v_i)\in U\mathbb{S}^n$ for each $i\in\{0,\ldots,l\}$.
Note that for a point $(p,v)\in U\mathbb{S}^n$ there is a unique length-minimizing antipodal geodesic $\gamma\colon [0,1]\to \mathbb{S}^n$ going from $p$ to $-p$ such that $\Dot{\gamma}(0) = \pi v$.
We denote this geodesic by $\gamma_{(p,v)}\colon [0,1]\to \mathbb{S}^n$.
For $m\geq 0$ we shall now describe embeddings $\Gamma_{2m+1}\to \Lambda \mathbb{S}^n$ and $\Gamma_{2m}\to P_a\mathbb{S}^n$, respectively.
First, let $l = 2m+1$ be odd with $m\geq 0$.
Then define $f_{2m+1}\colon \Gamma_{2m+1}\to \Lambda \mathbb{S}^n$ by setting
$$      f_{2m+1}(p,v_0,\ldots, v_{2m+1})(t) =  \begin{cases}
     \gamma_{(p,v_0)}((2m+2) t) , & 0\leq t \leq \tfrac{1}{2m+2}
     \\
     \gamma_{(-p,-v_1)}((2m+2) t - 1) , & \tfrac{1}{2m+2}\leq t \leq \tfrac{2}{2m+2} \\
     \vdots & \vdots \\
     \gamma_{(p,v_{2m})}((2m+2)t-2m), & \tfrac{2m}{2m+2} \leq t \leq \tfrac{2m+1}{2m+2} \\
     \gamma_{(-p,-v_{2m+1})}((2m+2)t- (2m+1)), & \tfrac{2m+1}{2m+2} \leq t\leq 1 .
\end{cases}    $$
Note in particular that the composition 
$$     f_{2m+1}\circ s_{2m+1} \colon U\mathbb{S}^n\to \Lambda \mathbb{S}^n    $$
is the embedding of the unit tangent bundle $U\mathbb{S}^n$ as the critical manifold of closed geodesics of multiplicity $m+1$.
Similarly, for $l = 2m$ let $f_{2m}\colon \Gamma_{2m}\to P_a\mathbb{S}^n$ be the map defined by
$$     f_{2m}(p,v_0,v_1,\ldots, v_{2m})(t) = \begin{cases}
     \gamma_{(p,v_0)}((2m+1) t) , & 0\leq t \leq \tfrac{1}{2m+1}
     \\
     \gamma_{(-p,-v_1)}((2m+1) t - 1) , & \tfrac{1}{2m+1}\leq t \leq \tfrac{2}{2m+1} \\
     \vdots & \vdots \\
     \gamma_{(-p,-v_{2m-1})}((2m+1)t- (2m-1)), & \tfrac{2m-1}{2m+1} \leq t \leq \tfrac{2m}{2m+1} \\
     \gamma_{(p,v_{2m})}((2m+1)t- 2m), & \tfrac{2m}{2m+1} \leq t\leq 1 .
\end{cases}    $$
Again, the composition $f_{2m}\circ s_{2m}\colon U\mathbb{S}^n\to P_a\mathbb{S}^n$ gives the critical manifold of antipodal geodesics of multiplicity $m$.
For both even and odd $l$ one checks that the only critical points in the image of $f_l$ are the points that come from the image of $s_l\colon U\mathbb{S}^n\to \Gamma_l$.
Hence, if we compose the embedding $f_l$ with an arbitrarily short flow of the gradient flow of the energy functional then we obtain an embedding $\widetilde{f}_l\colon \Gamma_l\to X$ such that the conditions of Definition \ref{def_completing_mfld} are satisfied, where $X$ is either $\Lambda \mathbb{S}^n$ or $P_a\mathbb{S}^n$.
Moreover, one checks from the computation of the index in Lemma \ref{lemma_index-and-nullity} that $\mathrm{dim}(\Gamma_l) = \mathrm{dim}(U\mathbb{S}^n) + \mathrm{ind}(\gamma) $ where $\gamma$ is a critical point of the energy functional in either the free loop space or the space of antipodal paths. 
Since for all $l\geq 0$ the manifold $\Gamma_l$ is orientable, we have thus established the following.
\begin{prop}\label{prop_existence_completing_manifolds}
    Consider the sphere $\mathbb{S}^n$ with the standard metric and consider the energy functional on the space of antipodal paths on the sphere.
    Then each critical manifold has a completing manifold.
\end{prop}
In Appendix \ref{sec_appendix_orientations} we shall specify explicitly how we orient the critical manifolds as well as the completing manifolds in the case that $n$ is even.
This will be important in the computation of the extended loop product in Section \ref{sec_computation_spheres}.

Note that we also found completing manifolds for the critical manifold of positive energy in the free loop space $\Lambda \mathbb{S}^n$.
In this case also the constant loops are a critical manifold $\Sigma_0\cong \mathbb{S}^n$.
This critical manifold is a completing manifold for itself since the constant loops have Morse index $0$.
Hence, the statement of Proposition \ref{prop_existence_completing_manifolds} also holds for the free loop space of $\mathbb{S}^n$ with the standard metric, see also \cite{ziller:1977} and \cite{oancea:2015}.

In the next subsection we shall use the completing manifolds to compute the homology of $P_a\mathbb{S}^n$.
Before we do so we shall note some more properties of the completing manifolds $\Gamma_l$ which we will later use for the computation of the string topology operations.
Consider the fiber bundle $\pi\colon U\mathbb{S}^n\to \mathbb{S}^n$.
For $l\geq 0$ the composition $\pi_l := \pi\circ p_l\colon \Gamma_l\to \mathbb{S}^n$ is a fiber bundle as well.
For $l,m\geq 0$ we consider the fiber product $\Gamma_l\times_{\mathbb{S}^n}\Gamma_m$.
Note that we have
$$    \Gamma_l\times_{\mathbb{S}^n}\Gamma_m = \{ (p,v_0,\ldots, v_l),(q,w_0,\ldots, w_m)\in \Gamma_l\times \Gamma_m \,|\,  q = p \}   $$
and hence there is a diffeomorphism $   \Phi_0\colon \Gamma_l\times_{\mathbb{S}^n}\Gamma_m \to \Gamma_{l+m+1}      $ given by
$$    \Phi_0((p,v_0,\ldots,v_l),(p,w_0,\ldots, w_m)) =  (p,v_0,\ldots, v_l,w_0,\ldots, w_m) .     $$
Moreover, we look at the fiber product
$$   \Gamma_l\times_{\Delta_a}\Gamma_m =  \{  (p,v_0,\ldots, v_l),(q,w_0,\ldots, w_m)\in \Gamma_l\times \Gamma_m \,|\,  q = -p \}    . $$
Again there is a diffeomorphism
$$   \Phi_1\colon \Gamma_l \times_{\Delta_a}\Gamma_m\to \Gamma_{l+m+1}    $$
given by
$$     \Phi_1((p,v_0,\ldots,v_l),(-p,w_0,\ldots,w_m)) = (p,v_0,\ldots, v_l, - w_0,\ldots, -w_m) .   $$
We shall now see that these diffeomorphisms as well as the embeddings of the completing manifolds are compatible with concatenation of paths.
Let $\phi_0\colon U\mathbb{S}^n\times_{\mathbb{S}^n}U\mathbb{S}^n\to U\mathbb{S}^n$ be the map which is induced by the projection on the first factor.
Similarly, let $\phi_1\colon U\mathbb{S}^n\times_{\Delta_a}U\mathbb{S}^n\to U\mathbb{S}^n$ be the map induced by the projection on the first factor.
Moreover, for $l,m\geq 0$ let $\tau= \tfrac{l+1}{l+m+2}$ and for $\gamma,\sigma\in P\mathbb{S}^n$ such that $\gamma(1) = \sigma(0)$ let $\mathrm{concat}^{\tau}(\gamma,\sigma)$ be the concatenated path
$$    \mathrm{concat}^{\tau}(\gamma,\sigma) (t) =\begin{cases}
    \gamma( \tfrac{t}{\tau}), & 0\leq t\leq \tau \\
    \sigma(\tfrac{t-\tau}{1-\tau}), & \tau\leq t\leq 1 .
\end{cases}     $$
\begin{lemma} \label{completing_manifolds_diffeo}
The diffeomorphisms $\Phi_0$ and $\Phi_1$ are compatible with the completing manifold structure in the following sense.
\begin{enumerate}
    \item 
If $l$ and $m$ are even then the following diagram commutes.
$$
\begin{tikzcd}
\Gamma_l\times_{\Delta_a} \Gamma_m \arrow[]{r}{f_l\times f_m} 
\arrow[d , "\Phi_1", "\cong"']  
& [2em]  P_a\mathbb{S}^n \hphantom{i}_{\mathrm{ev}_1}\times_{\mathrm{ev}_0} P_a \mathbb{S}^n \arrow[d, "\mathrm{concat}^{\tau}" ]
\\
\Gamma_{l+m+1} \arrow[]{r}{f_{l+m+1}} & \Lambda \mathbb{S}^n
\end{tikzcd}
$$
\item 
If $l$ is even and $m$ is odd then the following diagram commutes.
$$
\begin{tikzcd}
\Gamma_l\times_{\Delta_a} \Gamma_m \arrow[]{r}{f_l\times f_m} 
\arrow[d , "\Phi_1", "\cong"']  
& [2em]  P_a\mathbb{S}^n \hphantom{i}_{\mathrm{ev}_1}\times_{\mathrm{ev}_{\Lambda}} \Lambda \mathbb{S}^n \arrow[d, "\mathrm{concat}^{\tau}" ]
\\
\Gamma_{l+m+1} \arrow[]{r}{f_{l+m+1}} & P_a \mathbb{S}^n
\end{tikzcd}
$$
\item 
If $l$ is odd and $m$ is even then the following diagram commutes.
$$
\begin{tikzcd}
\Gamma_l\times_{\mathbb{S}^n} \Gamma_m \arrow[]{r}{f_l\times f_m} 
\arrow[d , "\Phi_0", "\cong"']  
& [2em]  \Lambda\mathbb{S}^n \hphantom{i}_{\mathrm{ev}_{\Lambda}}\times_{\mathrm{ev}_0} P_a \mathbb{S}^n \arrow[d, "\mathrm{concat}^{\tau}" ]
\\
\Gamma_{l+m+1} \arrow[]{r}{f_{l+m+1}} & P_a \mathbb{S}^n
\end{tikzcd}
$$
\item
If $l$ and $m$ are odd then the following diagram commutes.
$$
\begin{tikzcd}
\Gamma_l\times_{\mathbb{S}^n} \Gamma_m \arrow[]{r}{f_l\times f_m} 
\arrow[d , "\Phi_0", "\cong"']  
& [2em]  \Lambda\mathbb{S}^n \hphantom{i}_{\mathrm{ev}_{\Lambda}}\times_{\mathrm{ev}_{\Lambda}} \Lambda \mathbb{S}^n \arrow[d, "\mathrm{concat}^{\tau}" ]
\\
\Gamma_{l+m+1} \arrow[]{r}{f_{l+m+1}} & \Lambda \mathbb{S}^n
\end{tikzcd}
$$
\end{enumerate}

Moreover, the diffeomorphisms $\Phi_0$ and $\Phi_1$ are compatible with the submersions $p_l\colon \Gamma_l\to U\mathbb{S}^n$ in the sense that the diagrams
$$
\begin{tikzcd}
\Gamma_l\times_{\mathbb{S}^n} \Gamma_m \arrow[]{r}{p_l\times p_m} 
\arrow[d , "\Phi_0", "\cong"']  
&  [2em]  U\mathbb{S}^n\times_{\mathbb{S}^n}U\mathbb{S}^n \arrow[d, "\phi_0" ]
\\
\Gamma_{l+m+1} \arrow[]{r}{p_{l+m+1}} & U\mathbb{S}^n 
\end{tikzcd}
$$
and
$$
\begin{tikzcd}
\Gamma_l\times_{\Delta_a} \Gamma_m \arrow[]{r}{p_l\times p_m} 
\arrow[d , "\Phi_1", "\cong"']  
&  [2em]  U\mathbb{S}^n\times_{\Delta_a}U\mathbb{S}^n \arrow[d, "\phi_1" ]
\\
\Gamma_{l+m+1} \arrow[]{r}{p_{l+m+1}} & U\mathbb{S}^n 
\end{tikzcd}
$$
commute.
Finally, for $n$ even the diffeomorphisms $\Phi_0$ and $\Phi_1$ are orientation-preserving.  
\end{lemma}
\begin{proof}
By going through the respective definitions one can check that the diagrams in the statement of the lemma commute.
We prove that $\Phi_0$ and $\Phi_1$ are orientation-preserving in Appendix \ref{sec_appendix_orientations}.
\end{proof}

\subsection{The homology of the space of antipodal paths}

In this subsection we use the completing manifolds to compute the homology of the space of antipodal paths on the sphere.

\begin{theorem}\label{theorem_homology_additively}
    Consider the space of antipodal paths $P_a\mathbb{S}^n$ and take homology with integer coefficients.
    There is an isomorphism
    $$    \mathrm{H}_i(P_a\mathbb{S}^n) \cong  \bigoplus_{m=0}^{\infty}  \mathrm{H}_{i- m(2n-2)} (U\mathbb{S}^n) .     $$
    Explicitly, we obtain the following:
    \begin{enumerate}
        \item For $n\geq 3$ odd, we have
        $$      \mathrm{H}_i(P_a\mathbb{S}^n)  \cong \begin{cases}
            \mathbb{Z} , & i = m(n-1) , m\geq 0 \,\,\,\text{or}\,\,\, i = m(n-1) + n, m\geq 0
            \\
            0 & \text{else} .
        \end{cases}         $$
        \item For $n\geq 4$ even, we have
        $$      \mathrm{H}_i(P_a\mathbb{S}^n)  \cong \begin{cases}
            \mathbb{Z} , & i = m(2n-2) , m\geq 0 \,\,\,\text{or}\,\,\, i = m(2n-2) + 2n-1, m\geq 0
            \\
            \mathbb{Z}_2, & i = m(2n-2) + n-1, m\geq 0
            \\
            0 & \text{else} .
        \end{cases}         $$
        \item For $n = 2$ we have
        $$
            \mathrm{H}_i(P_a\mathbb{S}^2) \cong \begin{cases}
                \mathbb{Z}, & i = 2m, \,\,\, m\geq 0 \\
                \mathbb{Z}_2, & i = 1 \\
                \mathbb{Z}\oplus \mathbb{Z}_2, & i= 2m+1, \,\,\, m\geq 1 . 
            \end{cases}
        $$
    \end{enumerate}
\end{theorem}
\begin{proof}
    Since every critical manifold has a completing manifold the first isomorphism follows from the behavior of the index, see Lemma \ref{lemma_index-and-nullity} and the existence of the completing manifolds.
    For the explicit description of the homology we use the well-known homology of the unit tangent bundle of the sphere, see e.g. \cite[p. 20]{ziller:1977}.
\end{proof}

\begin{remark}
    The computation of the homology of $P_a\mathbb{S}^n$ for $n$ odd does not require completing manifolds of course, since by Proposition \ref{prop_homotopy_equivalence_n_odd} that $P_a\mathbb{S}^n\simeq \Lambda \mathbb{S}^n$.
\end{remark}

The computation of the homology in Theorem \ref{theorem_homology_additively} shows the following.
\begin{cor}
    Let $n\geq 2$ be even. 
    Then the space of antipodal paths $P_a\mathbb{S}^n$ is not homotopy equivalent to the free loop space $\Lambda \mathbb{S}^n$. 
    Moreover, for $n\geq 4$ the fibration $\mathrm{ev}_0\colon P_a\mathbb{S}^n\to \mathbb{S}^n$, $\gamma\mapsto \gamma(0)$ does not admit a section.
\end{cor}
\begin{proof}
     Comparing the homology of $P_a\mathbb{S}^n$ for $n$ even to the homology of the free loop space $\Lambda \mathbb{S}^n$, see e.g. \cite[p. 21]{ziller:1977}, we see that they are not isomorphic.
  
    For the second statement assume that there is a section $s\colon \mathbb{S}^n\to P_a\mathbb{S}^n$.
    In particular $s$ is then an embedding of $\mathbb{S}^n$ and we can consider the long exact sequence of the pair $(P_a\mathbb{S}^n,\mathbb{S}^n)$.
    Since the evaluation $\mathrm{ev}_0$ is a left-inverse to the inclusion $s\colon \mathbb{S}^n\to P_a\mathbb{S}^n$ we see that for each $i\geq 0$ we obtain a split short exact sequence
    $$   0\xrightarrow[]{} \mathrm{H}_i(\mathbb{S}^n) \xrightarrow[]{s_*} \mathrm{H}_i(P_a\mathbb{S}^n) \xrightarrow[]{} \mathrm{H}_i(P_a\mathbb{S}^n,\mathbb{S}^n) \xrightarrow[]{} 0 .   $$
    In particular for $i=n$ we obtain an injection $   \mathbb{Z}\hookrightarrow \mathrm{H}_n(P_a\mathbb{S}^n)   $.
    However, by Theorem \ref{theorem_homology_additively} we have that $\mathrm{H}_n(P_a\mathbb{S}^n) = 0$, thus yielding a contradiction.
\end{proof}
\begin{remark}
       For $n$ odd the fibration $P_a\mathbb{S}^n\to \mathbb{S}^n$ does have sections.
       Note however, that there is no canonical section while for the free loop space of a closed manifold $M$ there is always a canonical section of the free loop fibration given by the embedding of $M$ as the constant loops.
\end{remark}

Following \cite{oancea:2015} we now graphically represent the homology of the space of antipodal paths on the sphere using the perfectness of the energy functional.
The following diagram is an intuitive way of understanding how the shifted homologies of the critical manifolds make up the homology of $P_a\mathbb{S}^n$ for both $n$ even and odd.
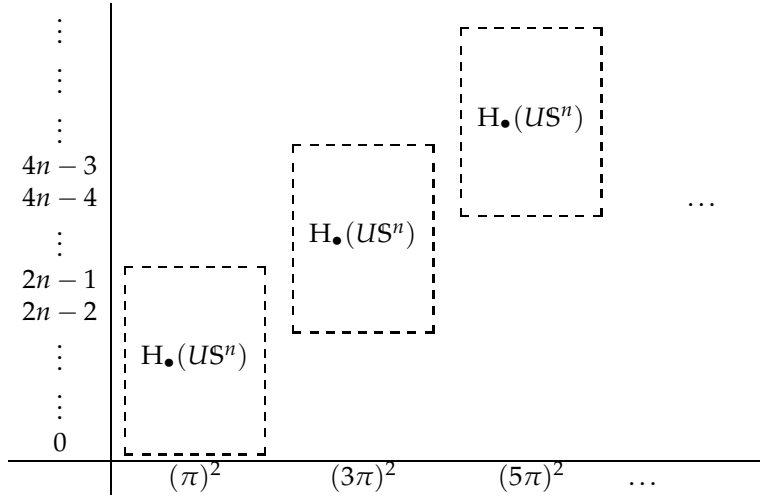
\begin{center}
\begin{tabular}[c]{c|cccccc}
$ \vdots $ &    &  &  &  &    &   \\
$ \vdots  $    &   &   &   &  & &    \\
$ \vdots  $ &   &   &  &  & & 
\\
$ 4n- 3$  &  &    &  &   & & \\
$ 4n -4   $ &  &    &   &    &   \dots & \\
$ \vdots   $ &  &  &  &   &   & \\
$ 2n-1 $  
&  &   &  &      &    & \\
$ 2n-2 $    &     &     &  &  & & \\
 $\vdots $ &   &  &    &   
 & & \\
 $\vdots $   &    &   &  &  &  &  \\
$ 0 $    &
\raisebox{0pt}[0pt][0pt]{\raisebox{62pt}{\dbox{$\xymatrix@R=18pt{ \\ \mathrm{H}_{\bullet}(U \mathbb{S}^n) \\ \quad }$}}} &
\raisebox{0pt}[0pt][0pt]{\raisebox{108pt}{\dbox{$\xymatrix@R=18pt{ \\  \mathrm{H}_{\bullet}(U \mathbb{S}^n) \\ \quad }$}}} &
\raisebox{0pt}[0pt][0pt]{\raisebox{152pt}{\dbox{$\xymatrix@R=18pt{ \\ \mathrm{H}_{\bullet}(U\mathbb{S}^n) \\ \quad }$}}} &
&
\\\hline
& $(\pi)^2$ & $(3 \pi)^2$ & $(5\pi)^2$ & \dots
\end{tabular}
\end{center} 
On the horizontal axis we denote the energy of the critical manifolds and on the vertical axis we denote the homological degree.
In contrast, the picture for the free loop space looks as follows.
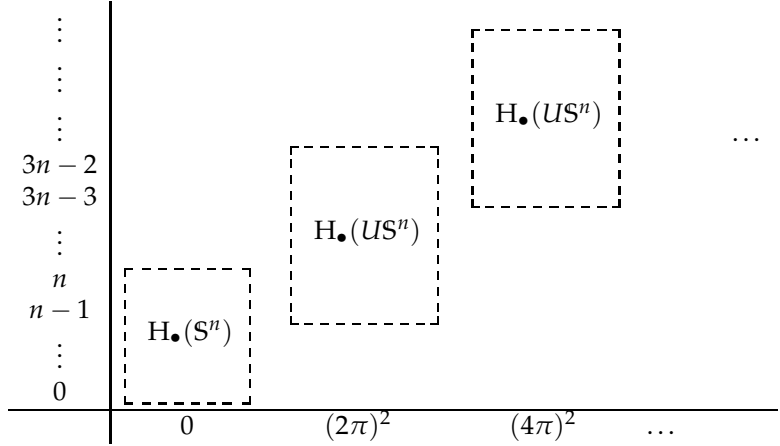
\begin{center}
\begin{tabular}[c]{c|cccccc}
 $ \vdots $  &   &    &  &  & & \\
 $ \vdots $  &   &    &    &  & & \\
 $\vdots$ &  &  &  &  & &  \dots\\
 $ 3n-2 $ &  &  &  &  & & \\
$3n-3$ &  &  &  &  & & \\
$\vdots$  &   &  &  &  &   & \\
$ n$ &  &  &  &       &    &    \\
$n-1$ &  &     &    &    &   & \\
$ \vdots$ &  &     &     &     &    &    \\
$ 0 $   &
\raisebox{0pt}[0pt][ 0pt]{\raisebox{42pt}{\dbox{ $\xymatrix@R=8pt{ \\ \mathrm{H}_{\bullet}( \mathbb{S}^n) \\ \quad }$}}} &
 \raisebox{ 0pt}[ 0pt][ 0pt]{ \raisebox{ 88pt} { \dbox{ $\xymatrix@R=16pt{ \\  \mathrm{H}_{\bullet}(U \mathbb{S}^n ) \\ \quad }$}}} &
 \raisebox{ 0pt }[ 0pt ][0pt ]{ \raisebox{132pt}{\dbox{ $\xymatrix@R=16pt{  \\ \mathrm{H}_{\bullet}(U\mathbb{S}^n) \\ \quad }$}}} &
&
\\\hline
& $0$ & $(2 \pi)^2$ & $(4\pi)^2$ & \dots
\end{tabular}
\end{center} 
Recall that $P_a\mathbb{S}^n \simeq \Lambda \mathbb{S}^n$ for $n$ odd and consequently $\mathrm{H}_{\bullet}(P_a\mathbb{S}^n)\cong \mathrm{H}_{\bullet}(\Lambda \mathbb{S}^n)$ in this case.
We find it remarkable that for odd $n$ we have thus found two different ways of building up the homology of the free loop space by stacking together copies of the homology of $U\mathbb{S}^n$.

Finally, we want to remark another interesting feature of the homology of the space of antipodal paths of an even-dimensional sphere.
From Theorem \ref{theorem_homology_additively} it follows easily that if we take $n$ to be even the direct sum $\widetilde{\mathrm{H}}_{\bullet}(\Lambda \mathbb{S}^n;\mathbb{Q})\oplus \mathrm{H}_{\bullet}(P_a\mathbb{S}^n;\mathbb{Q})$ behaves as follows.
$$ \widetilde{\mathrm{H}}_{\bullet}(\Lambda \mathbb{S}^n;\mathbb{Q})\oplus \mathrm{H}_{\bullet}(P_a\mathbb{S}^n;\mathbb{Q})  \cong 
\begin{cases}
            \mathbb{Q} , & i = m(n-1) , m\geq 0 \,\,\,\text{or}\,\,\, i = m(n-1) + n, m\geq 0
            \\
            0 & \text{else} .
        \end{cases}
$$
By $\widetilde{\mathrm{H}}_{\bullet}(X)$ we mean the reduced homology of a topological space $X$.
Thus comparing this with the homology of the free loop space of an odd-dimensional sphere, we see that \textit{formally} the direct sum $\widetilde{\mathrm{H}}_{\bullet}(\Lambda \mathbb{S}^n;\mathbb{Q})\oplus \mathrm{H}_{\bullet}(P_a\mathbb{S}^n;\mathbb{Q})$ has the homology of the free loop space of an odd-dimensional sphere.
We shall see in Section \ref{sec_computation_spheres} that - up to signs - this formal behavior even holds multiplicatively for the extended Chas-Sullivan product.


\section{The extended loop product on even-dimensional spheres}\label{sec_computation_spheres}

We use structure of the completing manifolds in the previous subsection to compute the extended loop product for even-dimensional spheres.
Throughout this section we always assume that $n\geq 2$ is even and we only consider homology with rational coefficients.
In this section we call the pairing $\wedge_a\colon \mathrm{H}_{\bullet}(P_a\mathbb{S}^n)\otimes \mathrm{H}_{\bullet}(P_a\mathbb{S}^n)\to \mathrm{H}_{\bullet}(\Lambda \mathbb{S}^n)$ the \textit{antipodal pairing}.

Let $k\in\mathbb{N}_0$ be even.
Recall from Section \ref{sec_completing_manifolds} that we had the completing manifolds $\Gamma_k$ which came with an embedding $f_k\colon \Gamma_k\to P_a\mathbb{S}^n$.
Moreover, there is a fiber bundle $p_k\colon \Gamma_k\to U\mathbb{S}^n$ and the composition
$$   (f_k)_*\circ (p_k)_! \colon \mathrm{H}_i(U\mathbb{S}^n) \to \mathrm{H}_{i+\lambda_k}(P_a \mathbb{S}^n)       $$
is injective where $\lambda_k = k(n-1)$.
Moreover, since $n$ is even we know that
$$   \mathrm{H}_i(U\mathbb{S}^n) \cong \begin{cases}
    \mathbb{Q},\quad i=0,2n-1\\
     \{0\} \quad \text{else} .
\end{cases}      $$
Hence, we have generators $[p_0]\in \mathrm{H}_0(U\mathbb{S}^n)$ where $p_0\in U\mathbb{S}^n$ is a point and the fundamental class $[U\mathbb{S}^n]\in \mathrm{H}_{2n-1}(U\mathbb{S}^n)$.
We thus get a set of generators of $\mathrm{H}_{\bullet}(P_a\mathbb{S}^n)$ by defining
$$   A_k =  (f_k)_* \circ (p_k)_! ([p_0])\in \mathrm{H}_{\lambda_k}(P_a\mathbb{S}^n)      \quad \text{and}\quad   B_k =  (f_k)_* \circ (p_k)_! ([U\mathbb{S}^n])\in \mathrm{H}_{\lambda_k+2n-1}(P_a\mathbb{S}^n)   .   $$
By Theorem \ref{theorem_homology_additively} the classes $A_k$ and $B_k$, $k\geq\mathbb{N}_0$ even generate all of $\mathrm{H}_{\bullet}(P_a\mathbb{S}^n)$.
Recall that in Section \ref{sec_completing_manifolds} we also constructed manifolds $\Gamma_k$ for odd $k$ and embeddings $f_k\colon \Gamma_k\to \Lambda \mathbb{S}^n$.
Thus, we can use the same construction to obtain a set of generators for the homology of the free loop space.
We define for odd $k\in\mathbb{N}$ the classes
$$    A_k' = (f_k)_*\circ (p_k)_!([p_0])\in \mathrm{H}_{\lambda_k}(\Lambda \mathbb{S}^n)\quad \text{and}\quad B_k' = (f_k)_*\circ (p_k)_![U\mathbb{S}^n]\in\mathrm{H}_{\lambda_k+2n-1}(\Lambda \mathbb{S}^n)        $$
again with $\lambda_k = k(n-1)$.
Let $q_0\in\Lambda \mathbb{S}^n$ be a constant loop.
Further we denote the image of the fundamental class in the free loop space by $[\mathbb{S}^n] \in\mathrm{H}_n(\Lambda \mathbb{S}^n)$.
Then the classes $A_k'$, $B_k'$ for odd $k\in\mathbb{N}$ together with the classes $[q_0]\in\mathrm{H}_0(\Lambda\mathbb{S}^n)$ and $[\mathbb{S}^n]\in\mathrm{H}_n(\Lambda \mathbb{S}^n)$ generate the homology of $\Lambda \mathbb{S}^n$ additively.
We now compute the extended loop product of $\mathbb{S}^n$ by explicitly computing all the occurring antipodal pairings and the module structure.
We begin with some preparations.
Consider the fiber product
$$   U\mathbb{S}^n\times_{\Delta_a}U\mathbb{S}^n = \{ (v,w)\in U\mathbb{S}^n\times U\mathbb{S}^n \,|\, \pi(v) = -\pi(w) \}    $$
where $\pi\colon U\mathbb{S}^n\to \mathbb{S}^n$ is the projection.
We have an inclusion $j\colon U\mathbb{S}^n\times_{\Delta_a}U\mathbb{S}^n\hookrightarrow U\mathbb{S}^n\times U\mathbb{S}^n$ as well as a projection $\phi_1\colon U\mathbb{S}^n\times_{\Delta_a}U\mathbb{S}^n\to U\mathbb{S}^n $ given by $\phi_1(v,w) = v$.
The map $\phi_1$ is a fiber bundle with fiber $\mathbb{S}^{n-1}$.
This fiber bundle has a section $s\colon U\mathbb{S}^n\to U\mathbb{S}^n\times_{\Delta_a}U\mathbb{S}^n$ given by $s(v) = (v,-v)$ where $-v\in T_{-p}\mathbb{S}^n$ for $v\in T_p\mathbb{S}^n$.
We consider the composition $\sigma = j\circ s\colon U\mathbb{S}^n\to U\mathbb{S}^n\times U\mathbb{S}^n$.
The Gysin map $\sigma_!\colon \mathrm{H}_i(U\mathbb{S}^n\times U \mathbb{S}^n)\to \mathrm{H}_{i-2n+1}(U\mathbb{S}^n)$ induces a product 
$$   \cap_a \colon \mathrm{H}_{i}(U\mathbb{S}^n) \otimes \mathrm{H}_j(U\mathbb{S}^n) \to \mathrm{H}_{i+j-2n+1}(U\mathbb{S}^n)       $$
by defining $X\cap_a Y := \sigma_! (X\times Y)$.
\begin{lemma}\label{lemma_twisted_intersection_product}
    The product $\cap_a$ satisfies
    $$   [p_0]\cap_a [p_0 ] = 0,\quad [p_0]\cap_a [U\mathbb{S}^n] = -[U\mathbb{S}^n]\cap_a [p_0] = [p_0] \quad \text{and}\quad [U\mathbb{S}^n]\cap_a[U\mathbb{S}^n] = [U\mathbb{S}^n] .   $$
\end{lemma}
\begin{proof}
    Since the product $\cap_a$ is given as the composition
    $$    \cap_a \colon \mathrm{H}_i(U\mathbb{S}^n)\otimes \mathrm{H}_j(U\mathbb{S}^n) \xrightarrow[]{\times}\mathrm{H}_{i+j}(U\mathbb{S}^n\times U\mathbb{S}^n) \xrightarrow[]{\sigma_!} \mathrm{H}_{i+j-2n+1}(U\mathbb{S}^n)    $$
    we need to understand how the Gysin map $\sigma_!$ behaves.
    
       On the class $[p_0] \times [p_0]\in\mathrm{H}_0(U\mathbb{S}^n\times U\mathbb{S}^n)$ the Gysin maps yield zero for degree reasons.
    Consider the class $[p_0]\times [U\mathbb{S}^n]$.
    Let $\varphi\in\mathrm{H}^{2n-1}(U\mathbb{S}^n)$ be the cohomology class dual to the fundamental class, i.e. the Kronecker pairing satisfies $\langle \varphi , [U\mathbb{S}^n]\rangle = 1$.
    Then we have $\mathrm{PD}^{-1}([p_0]\times [U\mathbb{S}^n]) = \varphi\times 1 = \mathrm{pr}_1^*\varphi$.
    Thus,
    $$   \sigma^*  \mathrm{PD}^{-1}([p_0]\times [U\mathbb{S}^n])  =  (\mathrm{pr}_1\circ \sigma)^* \varphi =  \varphi   $$
    since $\mathrm{pr}_1\circ \sigma = \mathrm{id}_{U\mathbb{S}^n}$.
    Hence we see that $\sigma_! ( [p_0]\times [U\mathbb{S}^n])  = [p_0]$.
    For the class $[U\mathbb{S}^n]\times [p_0]$ we have $\mathrm{PD}^{-1}([U\mathbb{S}^n]\times [p_0])= -1\times \varphi = - \mathrm{pr}_2^*\varphi$.
    Therefore, 
    \begin{equation}\label{eq_computation_with_degree_differential_antipodal_map}
           \sigma^*  (\mathrm{PD}^{-1}([U\mathbb{S}^n]\times [p_0]) ) =   -\sigma^* (\mathrm{pr}_2^* \varphi)  =  - (\mathrm{pr}_2\circ \sigma)^*\varphi .   
    \end{equation}
    The composition $\mathrm{pr}_2\circ \sigma \colon U\mathbb{S}^n\to U\mathbb{S}^n$ is the map induced by the differential of the antipodal map $Da \colon U\mathbb{S}^n\to U\mathbb{S}^n$.
    We claim that the degree of $Da$ is $1$.
    Indeed by the Lefschetz fixed point theorem the Lefschetz number 
    $$   \Lambda_{Da}  =  \mathrm{tr}(Da_*|_{\mathrm{H}_0(U\mathbb{S}^n)}) + (-1)^{2n-1} \mathrm{tr}(Da_*|_{\mathrm{H}_{2n-1}}(U\mathbb{S}^n))  $$
    must equal $0$ since $Da$ does not have fixed points.
    Since the first summand is obviously $1$ we must have that $Da_*|_{\mathrm{H}_{2n-1}(U\mathbb{S}^n)}$ is the identity.
    Consequently by equation \eqref{eq_computation_with_degree_differential_antipodal_map} we obtain that $$  \sigma_! ([U\mathbb{S}^n] \times [p_0]) = \mathrm{PD}(-(Da)^*\varphi) =  - [p_0] .  $$
    Finally, for $[U\mathbb{S}^n]\times [U\mathbb{S}^n]$ one sees with a similar computation that $\sigma_!([U\mathbb{S}^n]\times [U\mathbb{S}^n]) = [U\mathbb{S}^n]$.
\end{proof}
Note that the Gysin map $\sigma_! \colon \mathrm{H}_i(U\mathbb{S}^n\times U\mathbb{S}^n)\to \mathrm{H}_{i-2n+1}(U\mathbb{S}^n)$ factors through the maps
 $$   \mathrm{H}_{i}(U\mathbb{S}^n\times U\mathbb{S}^n) \xrightarrow[]{j_!} \mathrm{H}_{i-n}(U\mathbb{S}^n\times_{\Delta_a} U\mathbb{S}^n) \xrightarrow[]{s_!} \mathrm{H}_{i-2n+1}(U\mathbb{S}^n) .  $$
We also have a Gysin map $(\phi_1)_! \colon \mathrm{H}_i(U\mathbb{S}^n) \to \mathrm{H}_{i+n-1}(U\mathbb{S}^n\times_{\Delta_a}U\mathbb{S}^n)$.
\begin{lemma}\label{lemma_subspace_j_and_p_agree}
    The subspaces $\mathrm{im}(j_!)$ and $\mathrm{im}((\phi_1)_!)$ of $\mathrm{H}_{\bullet}(U\mathbb{S}^n\times_{\Delta_a}U\mathbb{S}^n)$ agree and there is a splitting
    $$  \mathrm{H}_{\bullet}(U\mathbb{S}^n\times_{\Delta_a}U\mathbb{S}^n) = \mathrm{im}((\phi_1)_!) \oplus \mathrm{ker}(s_!) .      $$
\end{lemma}
\begin{proof}
    The splitting follows from the fact that $s_!$ is surjective since $s_!\circ (\phi_1)_! = (\phi_1\circ s)_! = \mathrm{id}_{\mathrm{H}_{\bullet}(U\mathbb{S^n})}$.
    We show the equality $\mathrm{im}((\phi_1)_!) =  \mathrm{im}(j_!)$ explicitly.
    Since $E= U\mathbb{S}^n\times_{\Delta_a}U\mathbb{S}^n$ is the total space of an $\mathbb{S}^{n-1}$-bundle over $U\mathbb{S}^n$ with a section the Gysin sequence shows that the rational homology of $E$ is 
    $$   \mathrm{H}_i(E;\mathbb{Q}) \cong \begin{cases}
        \mathbb{Q}, \quad i = 0,n-1,2n-1,3n-2 \\ \{0\} & \text{else} .
    \end{cases}   $$
    We denote a generator in degree $n-1$ by $x$ and the generator in degree $3n-2$ is given by the fundamental class $[E]$.
    We claim that 
    $$     \mathrm{im}(j_!) = \mathrm{span}\{x,[E]\} = \mathrm{im}((\phi_1)_!) .     $$
    It is clear that $(\phi_1)_! = \mathrm{span}\{x,[E]\}$ by degree reasons since $(\phi_1)_!$ is injective.
    We explicitly compute the image of $j_!$.
    Note that the class $[p_0]\times [p_0]\in\mathrm{H}_0(U\mathbb{S}^n\times U\mathbb{S}^n)$ is mapped to zero for degree reasons.
    The classes $[p_0]\times [U\mathbb{S}^n]\in \mathrm{H}_{2n-1}(U\mathbb{S}^n)$ and $[U\mathbb{S}^n]\times [p_0]\in\mathrm{H}_{2n-1}(U\mathbb{S}^n)$ as well as the class $[U\mathbb{S}^n]\times [U\mathbb{S}^n]\in \mathrm{H}_{4n-2}(U\mathbb{S}^n \times U\mathbb{S}^n)$ have non-trivial image under $j_!$ since by Lemma \ref{lemma_twisted_intersection_product} these classes have non-trivial image under $\sigma_! =  s_!\circ j_!$.
    For degree reasons we therefore must have that $\mathrm{im}(j_!) = \mathrm{span}\{x,[E]\}$.
\end{proof}

\begin{theorem}
    Let $n\in\mathbb{N}$ be even.
    Consider the antipodal pairing $\wedge_a\colon \mathrm{H}_{\bullet}(P_a\mathbb{S}^n)^{\otimes2}\to \mathrm{H}_{\bullet}(\Lambda \mathbb{S}^n)$.
    For $k,l\in\mathbb{N}_0$ even we have
    $$     A_k \wedge_a A_l =  0,\quad A_k \wedge_a B_l = - B_k \wedge_a A_l =  A'_{k+l+1} \quad \text{and}\quad B_k \wedge B_l = B_{k+l+1}'        $$
\end{theorem}
\begin{proof}
Let $l,m\in\mathbb{N}$ be even numbers and set
$$   \lambda_l  =   l(n-1) \quad \text{and} \quad \lambda_m =  m(n-1) .      $$
We consider the following diagram where we write $P$ for the space of antipodal paths $P_a \mathbb{S}^n$.
$$
\begin{tikzcd}
    \mathrm{H}_{i-\lambda_l}(U\mathbb{S}^n) \otimes \mathrm{H}_{j-\lambda_m}(U\mathbb{S}^n) 
    \arrow[]{r}{(p_l)_!\otimes (p_m)_!}
    \arrow[]{d}{\times}
     \arrow[dr, phantom,shift left=0.8ex, "\mathrm{(1)}"]
    & [2.4em]
    \mathrm{H}_i(\Gamma_l)\otimes \mathrm{H}_j(\Gamma_m) 
    \arrow[]{r}{(f_l)_*\otimes (f_m)_*}
    \arrow[]{d}{\times}
     \arrow[dr, phantom,shift left=0.8ex, "\mathrm{(2)}"]
    & [2.4em]
    \mathrm{H}_i(P)\otimes \mathrm{H}_j(P)
    \arrow[]{d}{\times}
    \\
    \mathrm{H}_{i+j-\lambda_m-\lambda_l}(U\mathbb{S}^n\times U\mathbb{S}^n) 
    \arrow[]{r}{(p_l\times p_m)_!}
    \arrow[]{d}{j_!}
     \arrow[dr, phantom,shift left=0.8ex, "\mathrm{(3)}"]
    & 
    \mathrm{H}_{i+j}(\Gamma_l\times \Gamma_m)
    \arrow[]{r}{(f_l\times f_m)_*}
    \arrow[]{d}{J_!}
     \arrow[dr, phantom,shift left=0.8ex, "\mathrm{(4)}"]
    &
    \mathrm{H}_{i+j}(P\times P)
    \arrow[]{d}{}
    \\
    \mathrm{im}(j_!)
    \arrow[]{r}{(p_l\times p_m)_!}
    \arrow[]{d}{s_!}
     \arrow[dr, phantom,shift left=0.8ex, "\mathrm{(5)}"]
    &
    \mathrm{H}_{i+j-n}(\Gamma_l\times_{\Delta_a}\Gamma_m)
    \arrow[]{r}{(f_l\times f_m)_*}
    \arrow[]{d}{(\Phi_1)_*}
     \arrow[dr, phantom,shift left=0.8ex, "\mathrm{(6)}"]
    &
    \mathrm{H}_{i+j-n}(P  \hphantom{i}_{\mathrm{ev}_1}\times_{\mathrm{ev}_0} P)
    \arrow[]{d}{(\mathrm{concat}_{\tau})_*}
    \\
    \mathrm{H}_{i+j-\lambda_{m+l+1}-n}(U\mathbb{S}^n)
    \arrow[]{r}{(p_{l+m+1})_!} 
    &
    \mathrm{H}_{i+j-n}(\Gamma_{l+m+1})
    \arrow[]{r}{(f_{l+m+1})_*}
    &
    \mathrm{H}_{i+j-n}(\Lambda \mathbb{S}^n) 
\end{tikzcd}
$$

Here, $J\colon \Gamma_l\times_{\Delta_a}\Gamma_m\hookrightarrow \Gamma_l\times \Gamma_m$ is the inclusion.
We claim that the diagram commutes.
The first square commutes by \cite[Proposition VI.14.3]{bredon:2013} which says that we have 
$$   (p_l\times p_m)_! (a\times b) =  (-1)^{(\mathrm{dim}(U\mathbb{S}^n) + \mathrm{dim}(\Gamma_l)) (\mathrm{dim}(\Gamma_m) - |b|)} (p_l)_!(a) \times (p_m)_! (b)         $$
for $a\in\mathrm{H}_{\bullet}(U\mathbb{S}^n)$ and $b\in\mathrm{H}_{\bullet}(U\mathbb{S}^n)$.
Note that $\mathrm{dim}(\Gamma_l) = \mathrm{dim}(U\mathbb{S}^n) + l(n-1)$ and since $l$ is even we see that the sum of the dimensions of $U\mathbb{S}^n$ and $\Gamma_l$ is even.
We obtain
$$   (p_l\times p_m)_! (a\times b) =  (p_l)_!(a) \times (p_m)_! (b)         $$
for $a,b\in\mathrm{H}_{\bullet}(U\mathbb{S}^n)$.
The commutativity of the second square is clear.
For the third square we note that for maps $f\colon N\to M$ and $g\colon M\to K$ we have $(f\circ g)_! = g_! \circ f_!$ see \cite[Proposition VI.14.1]{bredon:2013}.
The commutativity of the third square thus follows from the fact that $(p_l\times p_m) \circ J = j\circ (p_l\times p_m)$.
For the fourth square we recall from \cite[Theorem VI.11.3]{bredon:2013} that for a Gysin map for an embedding $i\colon N\to M$ between finite-dimensional manifolds we have
$$     i_! x =   (-1)^{(\mathrm{dim}(M)-\mathrm{dim}(N))(\mathrm{dim}(M)-|x|)} R_* (\tau\cap x)      $$
for $x\in \mathrm{H}_{\bullet}(M)$ where $\tau$ is the Thom class of the embedding $N\hookrightarrow M$ and $R$ the retraction of a tubular neighborhood $U_N\to N$.
Since in our situation $\mathrm{dim}(\Gamma_l\times \Gamma_m) - \mathrm{dim}(\Gamma_l\times_{\Delta_a} \Gamma_m) = n$ is even the sign is trivial.
We thus have
$$   J_! a =  (R')_*(\tau'\cap a)     $$
where $\tau' = ((a\circ \pi\circ p_l)\times (\pi\circ p_m ))^*\tau_{\mathbb{S}^n} \in \mathrm{H}_n(U', U'\setminus \Gamma_l\times_{\Delta_a}\Gamma_m)$ is the Thom class and $U'$ is a tubular neighborhood of $\Gamma_l\times_{\Delta_a}\Gamma_m$ inside $\Gamma_l\times \Gamma_m$.
Note that we orient $\Gamma_l\times_{\Delta_a}\Gamma_m$ such that for a point $u\in \Gamma_l\times_{\Delta_a}\Gamma_m$ the orientations of $ N_u(\Gamma_l\times_{\Delta_a}\Gamma_m\hookrightarrow \Gamma_l\times \Gamma_m)\oplus  T_u \Gamma_l\times_{\Delta_a}\Gamma_m $ and $ T_u \Gamma_l\times \Gamma_m$ agree.
Here, the product $\Gamma_l\times \Gamma_m $ is equipped with the usual product orientation and the normal space inherits an orientation from $\mathbb{S}^n$ since the normal bundle is the pullback of the tangent bundle of the sphere. 
See Appendix \ref{sec_appendix_orientations} for details.
Recall that the map on the right hand side of square (4) is given by
$  (R_a)_*(\tau_{C} \cap x)      $
for $x\in\mathrm{H}_{i+j}(P\times P)$.
In order to show the commutativity of the fourth square we thus need to show that
$$     (R_a)_* (\tau_C\cap (f_l\times f_m)_* X)  =   (f_l\times f_m)_* R'_*(\tau'\cap X)        $$
for $X\in\mathrm{H}_{\bullet}(\Gamma_l\times \Gamma_m)$.
By naturality the left hand side of this equation yields
$$     (R_a)_* (\tau_C\cap (f_l\times f_m)_* X)  =  (R_a)_*(f_l\times f_m)_* ((f_l\times f_m)^*\tau_C\cap  X)  .    $$
Recall that $\tau_C = (\mathrm{ev}_1\times \mathrm{ev}_0)^* \tau_{\mathbb{S}^n}$ and hence
$   (f_l\times f_m)^* \tau_C = ( (\mathrm{ev}_1\times \mathrm{ev}_0) \circ (f_l\times f_m))^* \tau_{\mathbb{S}^n} $.
A direct computation shows that
$$   (\mathrm{ev}_1\times \mathrm{ev}_0)\circ (f_l\times f_m) (x,y) = ( (a\circ \pi\circ p_l)(x) ,   \pi\circ p_m(y)  ) \quad \text{for}\,\,\,x\in \Gamma_l,y\in\Gamma_m  $$
and therefore $(f_l\times f_m)^*\tau_C = \tau'$.
Since there is a homotopy $(f_l\times f_m)\circ R' \simeq R_a \circ (f_l\times f_m)$ it follows the fourth square commutes.
Note that if for the fifth square we had drawn the diagram
$$
\begin{tikzcd}
    \mathrm{H}_{i+j-\lambda_m-\lambda_l}(U\mathbb{S}^n\times_{\Delta_a}U\mathbb{S}^n) \arrow[]{r}{(p_l\times p_m)_!} \arrow[]{d}{s_!} & \mathrm{H}_{i+j-n}(\Gamma_l\times_{\Delta_a}\Gamma_m) \arrow[]{d}{(\Phi_1)_*} 
    \\
    \mathrm{H}_{i+j-\lambda_{m+l+1}-n}(U\mathbb{S}^n)\arrow[]{r}{(p_{l+m+1})_!} 
    &\mathrm{H}-{i+j-n}(\Gamma_{l+m+1})
\end{tikzcd}
$$
then this diagram can not commute since the composition along the top and the right vertical arrow is injective while the map $s_!$ has a non-trivial kernel.
We claim that this diagram however does commute if we restrict to classes $x\in\mathrm{H}_{\bullet}(U\mathbb{S}^n\times_{\Delta_a}U\mathbb{S}^n)$ such that $x\in\mathrm{im}(j_!)$.
By Lemma \ref{lemma_subspace_j_and_p_agree} we have $\mathrm{im}(j_!) =\mathrm{im}((\phi_1)_!)$ for the map $\phi_1\colon U\mathbb{S}^n\times_{\Delta_a}U\mathbb{S}^n\to U\mathbb{S}^n$ given by $\phi_1(p,q)  = p$.
One checks from the definitions that the diagram
$$
\begin{tikzcd}
    U\mathbb{S}^n\times_{\Delta_a}U\mathbb{S}^n   \arrow[]{d}{\phi_1}  & \Gamma_l\times_{\Delta_a}\Gamma_m \arrow[swap]{l}{p_l\times p_m} \arrow[]{d}{\Phi_1}
    \\
    U\mathbb{S}^n  & \Gamma_{l+m+1} \arrow[swap]{l}{p_{l+m+1}} 
\end{tikzcd}
$$
commutes and hence for $x = (\phi_1)_! y\in \mathrm{H}_{\bullet}(U\mathbb{S}^n\times_{\Delta_a}U\mathbb{S}^n)$ we get
$$   (\Phi_1)_!(p_{l+m+1})_!s_! x = (\Phi_1)_! (p_{l+m+1})_! y =  (p_l\times p_m)_! (\phi_1)_! y = (p_l\times p_m)_! x .     $$
Since $\Phi_1$ is an orientation-preserving diffeomorphism, see Appendix \ref{sec_appendix_orientations} we have that $(\Phi_1)_* = (\Phi_1)_!^{-1}$ and it follows that the fifth square commutes.
The commutativity of the sixth square is clear by Lemma \ref{completing_manifolds_diffeo}.
Using the explicit expressions in Lemma \ref{lemma_twisted_intersection_product} and the definitions of the classes $A_l$ and $B_l$, $l\in\mathbb{N}_0$ yields the result.
\end{proof}

In a similar manner one can compute the module structure and finds the following.
For $k\in\mathbb{N}$ odd and $m\in\mathbb{N}_0$ even we have
$$    A_k'  *_l A_m = 0, \quad A_k' *_l B_m =  B_k' *_l A_m =   A_{k+m+1}\quad \text{and}\quad B_k'  *_l B_m =  B_{k+m+1} .  $$
Note that there is indeed no minus sign in the relation $A_k'*_l B_m = B_k'*_l A_m$.
Similarly, for the right-module structure we get that for $k\in \mathbb{N}_0$ even and $m\in\mathbb{N}$ odd we have
$$    A_k *_r A_m' = 0, \quad A_k *_r B_m' = - B_k *_r A_m' =  A_{k+m+1} \quad \text{and}\quad B_k *_r B_m' =  B_{k+m+1} .      $$
This is almost a full computation of the antipodal pairing and the module structure.
However, there is still the critical manifold of constant loops on the free loop space which we need to take into account.
Recall that it induces classes $[p_0]\in\mathrm{H}_0(\Lambda \mathbb{S}^n)$ and $[\mathbb{S}^n]\in\mathrm{H}_n(\Lambda \mathbb{S}^n)$.
Since $[\mathbb{S}^n]$ is the unit of the Chas-Sullivan ring, we have $[\mathbb{S}^n] *_l X = X *_r [\mathbb{S}^n] = X$ for all $X\in\mathrm{H}_{\bullet}(P_a\mathbb{S}^n)$.
The multiplication with $[p_0]$ however is trivial for degree reasons, i.e. $[p_0] *_l X = X *_r [p_0] = 0$ for all $X\in\mathrm{H}_{\bullet}(P_a\mathbb{S}^n)$.
We have thus completely computed the extended loop product for even-dimensional spheres.

\begin{theorem}\label{theorem_full_extended_loop_product_even_dim_spheres}
    Consider the extended loop product for $\mathbb{S}^n$ with $n$ even and for the antipodal map $a\colon \mathbb{S}^n\to \mathbb{S}^n$.
    Take homology with rational coefficients.
    There is an isomorphism of graded algebras
    $$    \widehat{\mathbb{H}}_{\bullet}(\mathbb{S}^n) \cong \frac{\mathbb{Q}\langle A,B,U\rangle}{(A^2,  AU + UA ,B^2, BA, AB, BU, UB)} \quad \text{with} \,\,\, |A| = |B| =  -n,\,\,|U | = n-1 .        $$
\end{theorem}
\begin{proof}
    One checks with the above computations that the classes $$[p_0]\in\mathrm{H}_0(\Lambda \mathbb{S}^n),\quad [\mathbb{S}^n]\in\mathrm{H}_n(\Lambda \mathbb{S}^n),\quad A_0 \in \mathrm{H}_0(P_a\mathbb{S}^n) \quad \text{and}\quad B_0 \in\mathrm{H}_{2n-1}(P_a\mathbb{S}^n)  $$
    generate $\widehat{\mathbb{H}}_{\bullet}(\mathbb{S}^n)$ multiplicatively with $[\mathbb{S}^n]$ being the unit.
    We define $A := A_0$, $B:= [p_0]$ and $U:= B_0$ and check that if one makes the identifications
    $$   A_k' = AU^k,\quad B_k' = U^{k+1}, \quad A_m = AU^m\quad \text{and}\quad B_m = U^{m+1}    $$
    for $k$ odd and $m$ even then one gets the claimed isomorphism of algebras.
\end{proof}

Take the reduced homology $\widetilde{\mathrm{H}}_{\bullet}(\Lambda \mathbb{S}^n)$ and consider the direct sum $\widetilde{\mathrm{H}}_{\bullet}(\Lambda \mathbb{S}^n)\colon = \widetilde{\mathrm{H}}_{\bullet}(\Lambda \mathbb{S}^n)\oplus \mathrm{H}_{\bullet}(P_a\mathbb{S}^n)$.
One checks that the extended loop product can be restricted to $\widetilde{\mathrm{H}}_{\bullet}(\mathbb{S}^n)$ and we get the following result.
\begin{cor}\label{cor_extended_product_on_reduced_homology}
    The extended loop product induces an algebra on the direct sum
    $$    \widetilde{\mathbb{H}}_{\bullet}(\mathbb{S}^n) = \widetilde{\mathrm{H}}_{\bullet+n}(\Lambda \mathbb{S}^n) \oplus  \mathrm{H}_{\bullet+n}(P_a\mathbb{S}^n)     $$
    and this algebra is isomorphic to
     $$    \widetilde{\mathbb{H}}_{\bullet}(\mathbb{S}^n) \cong \frac{\mathbb{Q}\langle A,U\rangle}{(A^2,  AU + UA )} \quad \text{with} \,\,\, |A| = -n,\,\,|U | = n-1 .        $$
\end{cor}

\begin{remark}  As already mentioned in the introduction, it is remarkable that in Corollary \ref{cor_extended_product_on_reduced_homology} we formally almost recover the Chas-Sullivan algebra of an odd-dimensional sphere.
        In fact for $n$ odd there is an isomorphism of algebras $$(\mathrm{H}_{\bullet}( \Lambda \mathbb{S}^n),\wedge_{\mathrm{CS}})\cong \frac{\mathbb{Q}[A,U]}{(A^2)} \quad \text{with}\,\,\, |A| = -n,\,\, |U| = n-1 . $$
      \end{remark}
      \begin{remark}
    We note at this point that the map $B\colon \mathrm{H}_{\bullet}(P_a\mathbb{S}^n)\to \mathrm{H}_{\bullet+1}(P_a\mathbb{S}^n)$ induced by the circle action, see Remark \ref{remark_bv_algebra}, can also partially be computed using the completing manifolds.
    For degree reasons it is clear that $B(B_{2m}) = 0$ and we claim that $B(A_{2m} )=\lambda B_{2m-2}$ for some rational number $\lambda \neq 0$.
    The non-triviality of $B(A_{2m})$ can be seen with the completing manifolds and therefore we must have that $B(A_{2m})$ is a multiple of $B_{2m-2}$.
    However, it turns out that this computation becomes very technical and since we do not need the circle action on $P_a\mathbb{S}^n$ in the rest of this article we chose to omit this computation.
\end{remark}

\section{Extensions of the string topology coproduct}\label{sec_extensions_of_coproduct}

While we focussed on extensions of the loop product in the previous section we now consider extensions of the string topology coproduct.
In the setting of a closed oriented manifold with an involution we define a copairing and comodule structures.
We then show how these structures are related to the string topology coproduct in case that $f$ is homotopic to the identity.
In the last part of this section we compute these co-operations for even-dimensional spheres using again the completing manifolds that we studied above.

\subsection{Definition of the co-operations}

Let $M$ be a closed oriented manifold and let $f\colon M\to M$ be a fixed-point free involution.
Note that if we have a loop $\gamma\in \Lambda M$ and a time $s\in (0,1)$ such that
$   \gamma(s) = f(\gamma(0))   $
then the restrictions $\gamma|_{[0,s]}$ and $\gamma|_{[s,1]}$ are both in the path space $P_f M$.
We want to define a copairing in homology of the form $\mathrm{H}_{\bullet}(\Lambda M)\to \mathrm{H}_{\bullet+1-n}(P_f M\times P_f M)$.
Fix a Riemannian metric $g$ on $M$.
Since $M$ is compact there is a positive number $\rho_0 >0$ such that $\mathrm{d}_M(p,f(p)) \geq \rho_0$ for all $p\in M$, where $\mathrm{d}_M$ is the distance function induced by $g$.

Consider the map $\mathrm{ev}_f\colon \Lambda M\times I\to M\times M$ given by
$   \mathrm{ev}_f (\gamma,s) = (f(\gamma(0)),\gamma(s)) $ for  $\gamma\in \Lambda M$.
We consider the pullback of the diagonal 
$$   F_f = (\mathrm{ev}_f)^{-1}(\Delta M)  = \{ (\gamma,s)\in\Lambda M\times I\,|\, \gamma(s)=f(\gamma(0)) \} .    $$
Consider the oben neighborhood 
$$   U_f =    (\mathrm{ev}_f)^{-1}(U_M) = \{ (\gamma,s)\in\Lambda M\times I\,|\, \mathrm{d}(f(\gamma(0)),\gamma(s))< \epsilon\} .  $$
Moreover, for an $\epsilon_0> 0$ with $\epsilon_0<\epsilon$ define
$    U_{M,\geq\epsilon_0} = \{(p,q)\in U_M \,|\,  \mathrm{d}_M(p,q) \geq \epsilon_0\}       $.
The Thom class $\tau_M\in\mathrm{H}^n(U_M, U_M\setminus M)   $ 
induces a class in 
$   \tau'_M \in \mathrm{H}^n(U_M, U_{M,\geq\epsilon_0})     $.
We consider the space
$   U_{f,\geq\epsilon_0} = (\mathrm{ev}_f)^{-1}(U_{M,\geq\epsilon_0})      $
and pull back the Thom class to obtain a class
$$   \tau_f = (\mathrm{ev}_f)^*\tau'_M \in \mathrm{H}^n(U_f, U_{f,\geq\epsilon_0}) .     $$
We shall cap with this Thom class and understand this cap product as map
$$   \mathrm{H}^n(U_f,U_{f,\geq\epsilon_0})\otimes  \mathrm{H}_k(\Lambda M\times I,\Lambda M\times \partial I) \to \mathrm{H}_{k-n}(U_f)   $$
where we use a relative version of the cap product, see \cite[Appendix A]{hingston:2017}.
Furthermore we define a map $R_f\colon U_f\to F_f$ as follows.
Let $(\gamma,s)\in U_f$ and define the paths 
$$\gamma_1 = \mathrm{concat}(\gamma|_{[0,s]}, \overline{\gamma(s)f(\gamma(0))}) \quad \text{and}\quad \gamma_2 = \mathrm{concat}(\overline{f(\gamma(0))\gamma(s)},\gamma_{[s,1]}) . $$
Define a loop $\widetilde{\gamma}\in \Lambda M$ as the concatenation of $\gamma_1$ and $\gamma_2$ parametrized affine linearly in such a way that $\gamma_1$ is run through in the interval $[0,s]$ and $\gamma_2$ is run through in the interval $[s,1]$.
We then set $R_f (\gamma,s) =(\widetilde{\gamma},s)$.
The map $R_f$ is homotopic to a retraction $U_f\to F_f$ which can be seen similarly to \cite[Lemma 2.11]{hingston:2017}.
Define a map $\mathrm{cut}\colon F_f\to P_f M\times P_f M$ by setting
$$   \mathrm{cut}(\gamma,s) =  ( \gamma|_{[0,s]},\gamma|_{[s,1]})     $$
for $(\gamma,s)\in F_f$ where we reparametrize the restrictions of $\gamma$ by a linear map such that they are again parametrized on the unit interval.
Finally, let $[I]\in\mathrm{H}_1(I,\partial I)$ be the positively oriented generator of the homology of the unit interval relative to its boundary.
We now define the copairing.

\begin{definition}\label{definition_copairing}
    Consider a closed oriented $n$-dimensional manifold $M$ with a fix-point free involution $f\colon M\to M$.
    Take homology with coefficients in a commutative unital ring $R$.
    The \textit{loop copairing with respect to} $f$ is defined as the composition
\begin{eqnarray*}
     \vee_f  \colon  \mathrm{H}_{i}(\Lambda M) &\xrightarrow[]{\hphantom{coi}\times [I]\hphantom{coi}}& \mathrm{H}_{i+1}(\Lambda M\times I, \Lambda M\times \partial I)
     \\
     &\xrightarrow{ \hphantom{coi} \tau_f\cap \hphantom{coii} }& \mathrm{H}_{i+1-n}(U_f)
    \\
    &\xrightarrow[]{\hphantom{ci}(\mathrm{R}_f)* \hphantom{coi}} & \mathrm{H}_{i+1-n}(F_f) 
    \\ & \xrightarrow[]{\hphantom{icii}\mathrm{cut}_*\hphantom{ici}} & \mathrm{H}_{i+1-n}(P_f M \times P_f M) \, .
    \end{eqnarray*}
\end{definition}
\begin{remark}
    Consider the loop copairing $\vee_f$.
    \begin{enumerate}
        \item Note that if we take homology with coefficients in a field we can use the Künneth isomorphism to obtain a map
        $   \mathrm{H}_i(\Lambda M) \to \big(\mathrm{H}_{\bullet}(P_f M )\otimes \mathrm{H}_{\bullet}(P_f M)\big)_{i+1-n}    
        $.
        This justifies the name copairing.
        \item The definition of the copairing stays very close to the definition of the string topology coproduct in \cite{hingston:2017}.
        However, note that since $U_f\cap (\Lambda M\times \partial I) = \emptyset$, we do not have any relative terms whereas the string topology coproduct is naturally defined as a coproduct $\mathrm{H}_{\bullet}(\Lambda M,M)\to \mathrm{H}_{\bullet}(\Lambda M,M)^{\otimes 2}$.
        \end{enumerate}
\end{remark}

Next, we recall the definition of the string topology coproduct.
Consider the product of the free loop space with the unit interval $\Lambda M\times I$.
Define the map
$  \mathrm{ev}_I\colon \Lambda M\times I\to M\times M     $ by $   \mathrm{ev}_I(\gamma,s)  = (\gamma(0),\gamma(s)) \,    $.
Consider the spaces
\begin{equation*}
  U_{\mathrm{GH}} = (\mathrm{ev}_I)^{-1}(U_M) = \{ (\gamma,s) \in \Lambda M\times I\,|\, \mathrm{d}(\gamma(0),\gamma(s))< \epsilon\} 
  \end{equation*}
  and
  \begin{equation*}
  U_{\mathrm{GH},\geq\epsilon_0} =  (\mathrm{ev}_I)^{-1}(U_{M,\geq\epsilon_0}) =
  \{ (\gamma,s)\in U_{\mathrm{GH}}\,|\, \mathrm{d}(\gamma(0),\gamma(s))\geq \epsilon_0\}\,  .
\end{equation*}
We have that $U_{\mathrm{GH}}$ is an open neighborhood of the space  
$$ F_I =  (\mathrm{ev}_I)^{-1}(\Delta M)  =   \{ (\gamma,s)\in \Lambda\times I\,|\, \gamma(0)=\gamma(s)\}\,   . $$
The map $\mathrm{ev}_I$ induces a map of pairs $ \mathrm{ev}_I\colon    (U_{\mathrm{GH}} ,U_{\mathrm{GH},\geq\epsilon_0}) \to (U_M,U_{M,\geq\epsilon_0})     $
and we can pull back the class $\tau_M'$ to obtain a class
$$     \tau_{\mathrm{GH}} = (\mathrm{ev}_I)^* \tau_M'\in \mathrm{H}^n( U_{\mathrm{GH}} ,U_{\mathrm{GH},\geq\epsilon_0})\,  .       $$
Moreover, there is a retraction map $R_{\mathrm{GH}}\colon U_{\mathrm{GH}}\to F_I$ and we have a cutting map $\mathrm{cut}\colon F_I\to \Lambda M\times \Lambda M$ given by
$$   \mathrm{cut}(\gamma,s) = (\gamma|_{[0,s]},\gamma|_{[s,1]})  \quad \text{for}\,\,\,(\gamma,s)\in F_{I} .$$

\begin{definition} Let $M$ be a closed oriented manifold of dimension $n$.
Take homology with coefficients in a commutative unital ring $R$.
The \textit{string topology coproduct} is defined as the map
\begin{eqnarray*}
 \vee: \mathrm{H}_{\bullet}(\Lambda M,M) &\xrightarrow{\hphantom{coi}\times [I]\hphantom{coi}}& \mathrm{H}_{\bullet +1}(\Lambda M \times I, \Lambda M\times\partial I\cup M\times I)
     \\ &\xrightarrow{ \hphantom{oc} \tau_{\mathrm{GH}}\cap \hphantom{co} }& \mathrm{H}_{\bullet +1-n}(U_{\mathrm{GH}}, \Lambda M \times\partial I\cup M\times I) 
     \\
    &\xrightarrow{ \hphantom{c}({R}_{\mathrm{GH}})_* \hphantom{ci} }&
    \mathrm{H}_{\bullet +1-n}(F_{I}, \Lambda M \times\partial I\cup M\times I) \\
    &\xrightarrow{ \hphantom{cii} (\mathrm{cut})_* \hphantom{ic} }& \mathrm{H}_{\bullet +1-n}(\Lambda M\times\Lambda M, \Lambda M\times M\cup M\times \Lambda M)\, . 
    \end{eqnarray*}
\end{definition}

\begin{remark}
    The above definition of the string topology coproduct is called the \textit{Thom-signed} coproduct by Hingston and Wahl in \cite{hingston:2017}.
        In order to ensure good algebraic properties one needs to introduce additions signs if the manifold is odd-dimensional, see \cite{hingston:2017}.
\end{remark}

We define two more pairings as follows.
Consider the maps ${\mathrm{ev}}_{f,l},\mathrm{ev}_{f,r}\colon P_f M\times I\to M\times M$ given by
$$     \mathrm{ev}_{f,l}(\gamma,s) =  (\gamma(0),\gamma(s)) \quad \text{and}\quad \mathrm{ev}_{f,r}(\gamma,s)  =   (\gamma(s),\gamma(1))     \quad \text{for}\,\,\,   (\gamma,s)\in P_f M\times I  .   $$
Consider the sets
$$      F_{f,l} = (\mathrm{ev}_{f,l})^{-1}(\Delta M)  = \{(\gamma,s)\in P_f M\times I\,|\, \gamma(0) = \gamma(s) \}        $$
as well as 
$$     F_{f,r} = (\mathrm{ev}_{f,r})^{-1}(\Delta M) = \{(\gamma,s)\in P_f M\times I\,|\, \gamma(s) = \gamma(1) \} .         $$
We have tubular neighborhoods
$$   U_{f,l} = (\mathrm{ev}_{f,l})^{-1}(U_M) \quad \text{and}\quad U_{f,r} = (\mathrm{ev}_{f,r})^{-1}(U_M )    $$
as well as the sets
$$    U_{f,l,\geq\epsilon_0} = (\mathrm{ev}_{f,l})^{-1}(U_{M,\geq\epsilon_0}) \quad \text{and}\quad U_{f,r,\geq\epsilon_0} = (\mathrm{ev}_{f,r})^{-1}(U_{M,\geq\epsilon_0})  .     $$
We define cohomology classes 
$$   \tau_{f,l} = (\mathrm{ev}_{f,l})^*\tau_M'\in\mathrm{H}^n(U_{f,l}, U_{f,l,\geq\epsilon_0}) \quad \text{and}\quad \tau_{f,r} = (\mathrm{ev}_{f,r})^*\tau_M' \in \mathrm{H}^n(U_{f,r}, U_{f,r,\geq\epsilon_0}) .      $$
Furthermore, there are retractions $R_{f,l}\colon U_{f,l}\to F_{f,l}$ and $R_{f,r}\colon U_{f,r}\to F_{f,r}$.
Note that the cutting map induces maps
$$     \mathrm{cut}\colon F_{f,l} \to \Lambda M\times P_f M\quad \text{and}\quad \mathrm{cut}\colon F_{f,r} \to P_f M\times \Lambda M.      $$
With this input we can define the following operations.

\begin{definition}
    Let $M$ be a closed oriented $n$-dimensional manifold and $f$ a fixed point free involution.
    Take homology with coefficients in a commutative unital ring $R$.
    \begin{enumerate}
        \item We define the \textit{left comodule pairing} map $\vee_{f,l}\colon \mathrm{H}_i(P_f M)\to \mathrm{H}_{i+1-n}(\Lambda M\times P_f M,M\times P_f M)$ to be the composition
        \begin{eqnarray*}
     \vee_{f,l}  \colon  \mathrm{H}_{i}(P_f M) &\xrightarrow[]{\hphantom{coi}\times [I]\hphantom{coi}}& \mathrm{H}_{i+1}(P_f M\times I, P_f M\times \partial I)
     \\
     &\xrightarrow{ \hphantom{co} \tau_{f,l}\cap \hphantom{coi} }& \mathrm{H}_{i+1-n}(U_{f,l}, P_f M\times \{0\})
    \\
    &\xrightarrow[]{\hphantom{ci}(\mathrm{R}_{f,l})_* \hphantom{c}} & \mathrm{H}_{i+1-n}(F_{f,l}, P_f M\times \{0\}) 
    \\ & \xrightarrow[]{\hphantom{icii}\mathrm{cut}_*\hphantom{ici}} & \mathrm{H}_{i+1-n}(\Lambda M \times P_f M, M\times P_f M) \, .
    \end{eqnarray*}
    \item Moreover, we define the \textit{right comodule pairing} map $\vee_{f,r}\colon \mathrm{H}_i(P_f M)\to \mathrm{H}_{i+1-n}(P_f M\times \Lambda M,P_f M\times M)$ to be the composition
    \begin{eqnarray*}
     \vee_{f,r}  \colon  \mathrm{H}_{i}(P_f M) &\xrightarrow[]{\hphantom{coi}\times [I]\hphantom{coi}}& \mathrm{H}_{i+1}(P_f M\times I, P_f M\times \partial I)
     \\
     &\xrightarrow{ \hphantom{coi} \tau_{f,r}\cap \hphantom{co} }& \mathrm{H}_{i+1-n}(U_{f,r}, P_f M\times \{1\})
    \\
    &\xrightarrow[]{\hphantom{ci}(\mathrm{R}_{f,r})* \hphantom{c}} & \mathrm{H}_{i+1-n}(F_{f,r}, P_f M\times \{1\}) 
    \\ & \xrightarrow[]{\hphantom{icii}\mathrm{cut}_*\hphantom{ici}} & \mathrm{H}_{i+1-n}(P_f M \times \Lambda M, P_f M\times M) \, .
    \end{eqnarray*}
       \end{enumerate}
\end{definition}
\begin{remark}
    If we take field coefficients then $\vee_{f,l}$ and $\vee_{f,r}$ induce maps
    $$    \mathrm{H}_{\bullet}(P_f M)\to \mathrm{H}_{\bullet}(\Lambda M,M)\otimes \mathrm{H}_{\bullet}(P_f M)\quad \text{and}\quad \mathrm{H}_{\bullet}(P_f M)\to \mathrm{H}_{\bullet}(\Lambda M,M)\otimes \mathrm{H}_{\bullet}(P_f M) .   $$
    The author believes that the left and right comodule pairing map is indeed a comodule in the algebraic sense over the string topology coproduct coalgebra.
    The proof of this fact should stay relatively close to the proof of the graded coassociativity of the string topology coproduct in \cite[Theorem 2.14]{hingston:2017}.
    Note that the string topology coproduct is only graded coassociative so also in our setting we should expect some signs to appear in the comodule structure.
    We do not carry this out in detail in this manuscript since our focus in the rest of the paper will be on understanding these operations in the context of even-dimensional spheres and the applications in the next section.
\end{remark}

For the rest of this subsection we take homology with field coefficients.
We define an \textit{extended coproduct} on the graded vector space $\overline{\mathbb{H}}_{\bullet}(M) := \mathrm{H}_{\bullet}(\Lambda M,M)\oplus \mathrm{H}_{\bullet}(P_f M)$ by defining $\overline{\vee}_f \colon \overline{\mathbb{H}}_{\bullet}(M)\to \overline{\mathbb{H}}_{\bullet}(M)^{\otimes 2}$ as follows. 
Recall that the map $\vee_f$ is defined on the absolute homology $\mathrm{H}_{\bullet}(\Lambda M)$.
We can use it to define a map on $\mathrm{H}_{\bullet}(\Lambda M,M)$ by noting that the long exact sequence of the pair $(\Lambda M,M)$ splits because of the evaluation map $\mathrm{ev}\colon \Lambda M\to M$.
Hence, there is a canonical map $\iota \colon \mathrm{H}_{\bullet}(\Lambda M,M)\to \mathrm{H}_{\bullet}(\Lambda M)$ and hence we can precompose $\vee_f$ with this map to obtain a well-defined map on relative homology $\mathrm{H}_{\bullet}(\Lambda M,M)$.
For classes $X\in\mathrm{H}_{\bullet}(\Lambda M,M)$ and $Y\in\mathrm{H}_{\bullet}(P_f M)$ define the extended coproduct to be
$$    \overline{\vee}_f(X+Y) =   \vee (X) + \vee_f( \iota(X))  + \vee_{f,l} Y + \vee_{f,r} Y      .   $$

Finally, we define dual operations in cohomology.
Recall from \cite{hingston:2017} that for cohomology classes $\varphi,\psi\in \mathrm{H}^{\bullet}(\Lambda M,M)$ the \textit{Goresky-Hingston product} of $\varphi$ and $\psi$ is the cohomology class $\varphi\ostar_{\mathrm{GH}}\psi$ which is determined by demanding that
$$   \langle \varphi\ostar_{\mathrm{GH}} \psi, X \rangle =  \langle \varphi \times \psi, \vee X  \rangle  $$
for all $X\in\mathrm{H}_{\bullet}(\Lambda M,M)$.
This uniquely determines the class $\varphi\ostar_{\mathrm{GH}} \psi$ and thus defines a map $\ostar_{\mathrm{GH}} \colon \mathrm{H}^i(\Lambda M,M)\otimes \mathrm{H}^j(\Lambda M,M) \to \mathrm{H}^{i+j+n-1}(\Lambda M,M)$.
Similarly we define operations
\begin{eqnarray*}
   &   \ostar_f\colon \mathrm{H}^i(P_f M)\otimes \mathrm{H}^j(P_f M)\to \mathrm{H}^{i+j+n-1}(\Lambda M)\\ & \ostar_{f,l}\colon \mathrm{H}^i(\Lambda M,M)\otimes \mathrm{H}^j(P_f M)\to \mathrm{H}^{i+j+n-1}(P_f M)     
\end{eqnarray*}
as well as
 $$  \ostar_{f,r}\colon \mathrm{H}^i(P_f M)\otimes \mathrm{H}^j(\Lambda M,M)\to \mathrm{H}^{i+j+n-1}(P_f M)   $$
by defining them to be the dual of $\vee_f$, $\vee_{f,l}$ and $\vee_{f,r}$, respectively.
Note that the product $\ostar_f$ goes to the cohomology $\mathrm{H}^{\bullet}(\Lambda M)$.
Again, since we have split exact sequences
$$     0\to \mathrm{H}^{\bullet}(\Lambda M,M)\to \mathrm{H}^{\bullet}(\Lambda M) \to \mathrm{H}^{\bullet}(M)\to 0     $$
we get a map $\kappa \colon \mathrm{H}^{\bullet}(\Lambda M)\to \mathrm{H}^{\bullet}(\Lambda M,M)$.
Hence, the product $\ostar_f$ induces a map $\mathrm{H}^{\bullet}(P_f M)^{\otimes 2}\to \mathrm{H}^{\bullet+n-1}(\Lambda M,M)$.
Overall, we obtain a product on the direct sum $\overline{\mathbb{H}}^{\bullet}(M) = \mathrm{H}^{\bullet}(\Lambda M,M)\oplus \mathrm{H}^{\bullet}(P_f M)$ by setting
$$    (\varphi, \alpha)\,  \overline{\ostar}\, (\psi,\beta)  =   (\varphi\ostar_{\mathrm{GH}} \psi + \kappa (\alpha\ostar_f \beta ) \,, \, \varphi\ostar_r \beta + \alpha \ostar_l \psi)     $$
for $\varphi,\psi\in\mathrm{H}^{\bullet}(\Lambda M,M)$ and $\alpha,\beta\in\mathrm{H}^{\bullet}(P_f M)$.

\subsection{The copairing as a lift of the string topology coproduct}

We now study the relation between the copairing $\vee_f$ and the string topology coproduct in the case that the involution $f\colon M\to M$ is homotopic to the identity.
Recall from Section \ref{sec_involution_path_space_string_operations} that the pairing $\wedge_f\colon \mathrm{H}_i(P_f M)\otimes \mathrm{H}_j(P_f M)\to \mathrm{H}_{i+j-n}(\Lambda M)$ can be expressed by the Chas-Sullivan product in case that $f$ is homotopic to the identity.
An analogous statement for the copairing cannot hold, simply because the copairing takes value in the absolute homology $\mathrm{H}_{\bullet}(P_f M)^{\otimes 2}$ while the string topology coproduct maps to the relative homology $\mathrm{H}_{\bullet}(\Lambda M,M)^{\otimes 2}$.
However, we will show that the copairing can be understood as a lift of the string topology coproduct along the natural map $\mathrm{H}_{\bullet}(\Lambda M)\to \mathrm{H}_{\bullet}(\Lambda M,M)$.
Similar identities hold for the comodule structure as we will sketch in the end of this subsection.

Let $M$ be a closed oriented manifold and assume that $f\colon M\to M$ is a fixed-point free involution which is homotopic to the identity.
Let $H\colon M\times I\to M$ be a homotopy between $f$ and the identity of $M$, i.e.
$   H(\cdot, 0) =   f$ and $ H(\cdot, 1) =\mathrm{id}_M   $.
Define $\eta \colon M\to PM$ to be the map
$$  \eta (p ) (t) =  H(p,t) \quad \text{for}\,\,\,p\in M,\,t\in I .    $$
Furthermore we define a map $\Theta\colon \Lambda M\to \Lambda M$ by
$$   \Theta(\gamma) =  \mathrm{concat} (  \eta(\gamma(0)),\gamma,\overline{\eta(\gamma(0))})        $$
for $\gamma\in \Lambda M$.
Note that $\Theta$ is homotopic to the identity by the Spaghetti trick.
In order to compare the tubular neighborhoods $U_f$ and $U_{\mathrm{GH}}$ we need the following lemmas.

\begin{lemma}\label{lemma_rho_0}
Let $M$ be a closed manifold with an involution $f\colon M\to M$ without fixed points.
Let $\rho_0 >0$ be such that $\mathrm{d}_M(p,f(p)) \geq \rho_0$ for all $p\in M$.
Consider the energy functional on $\Lambda M$ with respect to the chosen Riemannian metric.
    Let $\gamma\in\Lambda M$ be a loop with $E(\gamma) = a$.
    Then for every $t\in [0,\tfrac{\rho_0^2}{8a}) \cup (1- \tfrac{\rho_0^2}{8a},1]$ we have
    $\mathrm{d}_M(\gamma(0),\gamma(t)) < \tfrac{\rho_0}{2}$.
\end{lemma}
\begin{proof}
    Assume that the statement is false, i.e. there is a $t\in [0,\tfrac{\rho_0^2}{8a})$ such that $\mathrm{d}_M(\gamma(0),\gamma(t))\geq \tfrac{\rho_0}{2}$.
    By the inequality relating the distance function and the energy functional \cite[Proposition 2.4.4]{klingenberg:1995} we have
    $$     \frac{\rho_0^2}{4}  =  \mathrm{d}_M^2 (\gamma(0),\gamma(t_0)) \leq L^2 \big(  \gamma|_{[0,t_0]}   \big) \leq 2 t_0 E\big(c|_{[0,t_0]}\big) \leq 2 t_0 a .    $$
    Hence, $t_0 \geq \frac{\rho^2}{8a}$ which yields a contradiction.
    The argument for $t\in (1- \tfrac{\rho_0^2}{8a},1] $ is analogous.
\end{proof}
In the following we will always assume that the number $\epsilon$ which we chose for the definition of the tubular neighborhood $U_M$ satisfies $\epsilon < \tfrac{\rho_0}{2}$.
For $\gamma\in \Lambda M$ we set 
$$  t_{\gamma}^0 = \begin{cases}
    \tfrac{1}{4}, & E(\gamma) \leq \tfrac{\rho_0^2}{2} \\
    \tfrac{\rho_0^2}{8 E(\gamma)}, & E(\gamma) \geq \tfrac{\rho_0^2}{2}  . 
\end{cases}    $$

The map $t^0\colon\Lambda M\to \mathbb{R}$, $\gamma\mapsto t_{\gamma}^0$ is continuous.
We define a map $i\colon \Lambda \mathbb{S}^n\times I\to I$ by setting
$$     i(\gamma,s) = \begin{cases}
    \tfrac{1+t_{\gamma}^0}{3 t_{\gamma}^0} s , & 0 \leq s \leq t_{\gamma}^0 \\
    \tfrac{1}{3}s  + \tfrac{1}{3}, & t_{\gamma}^0 \leq s \leq 1 -t_{\gamma}^0 \\
    \tfrac{1 + t_{\gamma}^0}{3 t_{\gamma}^0} s +  \tfrac{2 t_{\gamma}^0 - 1}{3 t_{\gamma}^0} , &   1 - t_{\gamma}^0 \leq s \leq 1 .
\end{cases}   $$
Define $\widetilde{\Theta}\colon \Lambda M\times I\to \Lambda M\times I$ by
$   \widetilde{\Theta} (\gamma,s) =  (\Theta(\gamma),i(\gamma,s)) 
       $.

\begin{lemma}\label{lemma_phi_does_what_it_should}
    The map $\widetilde{\Theta}$ induces a map of quadruples 
     $$    (\Lambda M \times I  , U_f ,U_{f,\geq\epsilon_0} , \Lambda M \times \partial I)\to  (\Lambda M \times I, U_{\mathrm{GH}}, U_{\mathrm{GH},\geq\epsilon_0}, \Lambda M \times \partial I) .   $$
     Moreover, $\widetilde{\Theta}$ is a homotopic to the identity map of the pair $(\Lambda M\times I,\Lambda M\times \partial I)$ and as maps $U_f \to U_{\mathrm{GH}}$ it holds that $\mathrm{ev}_I\circ \widetilde{\Theta} = \mathrm{ev}_f$.
\end{lemma}
\begin{proof}
    It is clear that $\widetilde{\Theta}$ maps $\Lambda M \times \partial I$ to itself.
    Let $(\gamma,s)\in U_f$.
    We need to show that $(\Theta(\gamma),i(\gamma,s))\in U_{\mathrm{GH}}$, i.e. that
    $    \mathrm{d}_M(\Theta(\gamma)(0),\Theta(\gamma)(i(\gamma,s))) < \epsilon      $.
    By definition we have
    $$   \mathrm{d}_M(f(\gamma(0)),\gamma(s)) < \epsilon < \tfrac{\rho_0}{2}        $$ and by Lemma \ref{lemma_rho_0} this implies that $s\in [t_{\gamma}^0,1-t_{\gamma}^0]$.
    In fact, if we had $s < t_{\gamma}^0$ then $s < \tfrac{\rho_0^2}{8 E(\gamma)}$ and by Lemma \ref{lemma_rho_0} it follows that $\mathrm{d}_M(\gamma(0),\gamma(s)) < \tfrac{\rho_0}{2}$.
    The triangle inequality then shows that $\mathrm{d}_M(\gamma(0),f(\gamma(0))) < \rho_0$ which is a contradiction.
    Similarly one argues that $s> 1-t_{\gamma}^0$ cannot hold.
    Consequently, $s\in[t_{\gamma}^0, 1 - t_{\gamma}^0]$ and hence $i(\gamma,s) = \tfrac{1}{3}s+\tfrac{1}{3}$.
    By definition of $\Theta$ we therefore get $      \Theta(\gamma)(i(\gamma,s)) = \gamma(s)     $.
    Moreover, since $\Theta(\gamma)(0) = f(\gamma(0))$ we see that
    $$      \mathrm{d}_M(\Theta(\gamma)(0),\Theta(\gamma)(i(\gamma,s))) = \mathrm{d}_M(f(\gamma(0)),\gamma(s)) < \epsilon .          $$
    The above argument also shows that $U_{a,\geq\epsilon_0}$ is mapped into $U_{\mathrm{GH},\geq\epsilon_0}$.
    The rest of the lemma can be seen easily.
\end{proof}

Let $\Phi_0\colon P_f M \to \Lambda M$ be the homotopy equivalence
$$    \Phi_0(\gamma) = \mathrm{concat}(\eta(\gamma(0),\gamma),\quad \gamma\in P_f M       $$
and let $\Phi_1\colon P_f M\to \Lambda M$ be the homotopy equivalence
$$   \Phi_1(\gamma) = \mathrm{concat}(\gamma,\overline{\eta(\gamma(1))}) ,\quad \gamma\in P_f M . $$
Define an operation 
$$   \vee'_f =  (\Phi_0\times \Phi_1)\circ \vee_f \colon \mathrm{H}_i(\Lambda M)\to \mathrm{H}_{i+1-n}(\Lambda M\times \Lambda M) .    $$
In the next theorem we shall see that this map is compatible with the string topology coproduct.

\begin{theorem}\label{theorem_copairing_as_lift}
    Let $M$ be a closed oriented manifold and let $f\colon M\to M$ be an involution which is homotopic to the identity.
    Take homology with coefficients in a commutative unital ring $R$.
    The map $\vee'_f$ is a lift of the string topology coproduct, i.e. the diagram
    $$
    \begin{tikzcd}
        \mathrm{H}_i(\Lambda M) \arrow[]{r}{\vee_f'} \arrow[]{d}{}
        &
        \mathrm{H}_{i+1-n}(\Lambda M\times \Lambda M) \arrow[]{d}{}
        \\
        \mathrm{H}_i(\Lambda M,M) \arrow[]{r}{\vee} 
        &
    \mathrm{H}_{i+1-n}(\Lambda M \times \Lambda M,\Lambda M\times M\cup M \times \Lambda M) 
    \end{tikzcd}
    $$
    commutes where the horizontal arrows are induced by the inclusion of pairs.
\end{theorem}
\begin{proof}
   We claim that the following diagram commutes. 
   $$    
   \begin{tikzcd}
       \mathrm{H}_i(\Lambda M )
       \arrow[]{r}{}
       \arrow[]{d}{\times [I]}
       & [2.5em]
       \mathrm{H}_i(\Lambda  M ,M)
       \arrow[]{d}{\times [I]}
        \\
        \mathrm{H}_{i+1}(\Lambda M \times I,\Lambda M \times \partial I)
        \arrow[]{r}{\widetilde{\Theta}_*}
        \arrow[]{d}{\tau_f\cap} 
        &
        \mathrm{H}_{i+1}(\Lambda M \times I,\Lambda M \times \partial I\cup M\times I)
        \arrow[]{d}{\tau_{\mathrm{GH}}\cap}
        \\
        \mathrm{H}_{i+1-n}(U_f) 
        \arrow[]{r}{\widetilde{\Theta}_*}
        \arrow[]{d}{(\mathrm{cut}\circ {R}_f)_*}
        &
        \mathrm{H}_{i+1-n}(U_{\mathrm{GH}}, M\times I\cup\Lambda M \times \partial I)
        \arrow[]{dd}{(\mathrm{cut}\circ {R}_{\mathrm{GH}})_*}
        \\
        \mathrm{H}_{i+1-n}(P_f M\times P_f M) 
        \arrow[]{d}{(\Phi_0\times \Phi_1)_*}
        &
        \\
        \mathrm{H}_{i+1-n}(\Lambda M \times \Lambda  M) \arrow[]{r}{}
        & \mathrm{H}_{i+1-n}(\Lambda M \times \Lambda M, \Lambda  M\times M\cup M\times \Lambda M)
   \end{tikzcd}
   $$
   The commutativity of the first square is clear since by Lemma \ref{lemma_phi_does_what_it_should} we have that
   $$  \widetilde{\Theta}_* \colon \mathrm{H}_{\bullet}(\Lambda M \times I,\Lambda M \times \partial I)\to \mathrm{H}_{i+1}(\Lambda M \times I,\Lambda M \times \partial I\cup M\times I)      $$
   is equal to the map which is induced by the inclusion of pairs $(\Lambda M \times I,\Lambda M \times \partial I)\to (\Lambda M \times I,\Lambda M \times \partial I\cup M\times I)$.
    For the second square let $X\in\mathrm{H}_{\bullet}(\Lambda M\times I,\Lambda M \times \partial I)$ then we have
$$         \tau_{\mathrm{GH}} \cap (\widetilde{\Theta}_* X) =   \widetilde{\Theta}_*  ( \widetilde{\Theta}^* \tau_{\mathrm{GH}} \cap X)      $$
by naturality.
We want to show that $\widetilde{\Theta}^* \tau_{\mathrm{GH}} = \tau_f$.
Recall that $\tau_{\mathrm{GH}} = \mathrm{ev}_I^* \tau_{f}$ and therefore
$   \widetilde{\Theta}^* \tau_{\mathrm{GH}} =   \big(  \mathrm{ev}_I\circ \widetilde{\Theta} \big)^*\tau_{f}    $.
By Lemma \ref{lemma_phi_does_what_it_should} we have
$     \mathrm{ev}_I \circ \widetilde{\Theta} = \mathrm{ev}_f     $ which shows the commutativity of the second square.
The commutativity of the last square follows from the commutativity of the underlying diagram of maps.
The statement of the theorem follows directly from the commutativity of the diagram.
\end{proof}

\begin{remark}
    The fact that $\vee_f'$ is a lift of the string topology coproduct seems to be closely related to the construction of lifts of the coproduct in \cite{cieliebak2023loop}.
    Moreover, note that the lift $\vee_f'$ depends on the choices of $\Phi_0$ and $\Phi_1$.
    This can be seen by explicit computations on the circle.
\end{remark}

One can derive similar statements for the comodule structure.
Define
$$    \vee_{f,l}'\colon \mathrm{H}_i(\Lambda M)\to \mathrm{H}_{i+1-n}(\Lambda M\times \Lambda M ,M\times \Lambda M )   $$ and $$ \vee_{f,r}'\colon \mathrm{H}_i(\Lambda  M)\to \mathrm{H}_{i+1-n}(\Lambda M\times \Lambda M,\Lambda  M\times M)  $$
as the compositions
$$   \vee_{f,l}' =     (\mathrm{id}\times \Phi_1)_*\circ    \vee_{f,l} \circ (\Psi_1)_*  \quad \text{and}\quad \vee_{f,r}' = (\Phi_0\times \mathrm{id})_* \circ \vee_{f,r} \circ (\Psi_0)_*  $$
where $\Psi_0,\Psi_1 \colon \Lambda M\to P_f M$ are homotopy inverses to $\Phi_0,\Phi_1$, respectively.
One can show that the diagrams
$$  
\begin{tikzcd}
    \mathrm{H}_i(\Lambda M) \arrow[]{r}{\vee_{f,l}'}  \arrow[]{d}{} & \mathrm{H}_{i+1-n}(\Lambda M\times \Lambda M, M\times \Lambda M) \arrow[]{d}{}
    \\
      \mathrm{H}_i(\Lambda M,M) \arrow[]{r}{\vee} 
        &
    \mathrm{H}_{i+1-n}(\Lambda M \times \Lambda M,\Lambda M\times M\cup M \times \Lambda M) 
\end{tikzcd}
$$
and 
$$  
\begin{tikzcd}
    \mathrm{H}_i(\Lambda M) \arrow[]{r}{\vee_{f,r}'}  \arrow[]{d}{} & \mathrm{H}_{i+1-n}(\Lambda M\times \Lambda M,\Lambda M\times M) \arrow[]{d}{}
    \\
      \mathrm{H}_i(\Lambda M,M) \arrow[]{r}{\vee} 
        &
    \mathrm{H}_{i+1-n}(\Lambda M \times \Lambda M,\Lambda M\times M\cup M \times \Lambda M)  
\end{tikzcd}
$$
commute where the vertical arrows are induced by the respective inclusions of pairs.
The proof follows along the same lines as the proof of Theorem \ref{theorem_copairing_as_lift}.

\subsection{The extended coproduct for even-dimensional spheres}

In this final subsection we explicitly compute the copairing for even-dimensional spheres, taking homology with rational coefficients.
We begin by considering the cohomology ring of the completing manifolds.
Recall that for $l\in\mathbb{N}$ the manifold $\Gamma_l$ is defined as the $l+1$-fold fibered product over $\mathbb{S}^n$ of the unit tangent bundle $U\mathbb{S}^n$ with itself.
We have a sequence of maps $$\Gamma_l\xrightarrow[]{\rho_l}\Gamma_{l-1}\xrightarrow[]{\rho_{l-1}} \ldots \xrightarrow[]{\rho_{2}} \Gamma_1 \xrightarrow[]{\rho_{1}} \Gamma_0 = U\mathbb{S}^n$$ where each map is given by dropping the last factor.
These maps are all sphere bundles with a section and hence $\Gamma_l$ is an iterated sphere bundle over $U\mathbb{S}^n$.
A section to the map $\rho_i$ is given by the map $\sigma_i\colon \Gamma_{i-1}\to \Gamma_i$, $i\in\{1,\ldots,l\}$
with
$$ \sigma_i(p,v_0,\ldots, v_{i-1}) =  (p,v_0,\ldots, v_{i-1},v_{i-1})  .    $$
Since each fiber is an odd-dimensional sphere $\mathbb{S}^{n-1}$ we can read off the cohomology ring of $\Gamma_l$ by inductively considering the Gysin sequences of the respective steps of the iterated sphere bundle.
We obtain an isomorphism
$$   \mathrm{H}^{\bullet}(U\mathbb{S}^n) \cong \frac{\mathbb{Q}[\alpha,\xi_1,\ldots, \xi_l]}{(\alpha^2,\xi_1^2, \ldots, \xi_l^2)}     $$
with $|\alpha| = 2n-1$, $|\xi_1| = \ldots |\xi_l| = n-1$.
Recall that for $l = 2m+1$ odd we have an embedding $f_l\colon \Gamma_l\to \Lambda \mathbb{S}^n$.
In order to compute the antipodal copairing we need to consider the space $U_{\Gamma_l} = (f_l\times \mathrm{id}_I)^{-1}(U_a)$ where $U_a$ is the open set defined before Definition \ref{definition_copairing}.
Since all paths in the image of $f_l$ are broken geodesics one sees that there is a $\delta>0$ such that
$$   U_{\Gamma_l} =  \Gamma_l \times \bigsqcup_{i=1}^{m+1} (\tfrac{2i-1}{2m+2}- \delta, \tfrac{2i-1}{2m+2} + \delta ) .     $$
Moreover, we define $   U_{\Gamma_l,\geq\epsilon_0}  := (f\times \mathrm{id}_I)^{-1}(U_{f,\geq\epsilon_0})   $ and we find that there is a $\delta_0>0$ with $\delta_0 < \delta$ such that
$$   U_{\Gamma_l,\geq\epsilon_0} =    \Gamma_l \times \bigsqcup_{i=1}^{m+1} (\tfrac{2i-1}{2m+2}- \delta, \tfrac{2i-1}{2m+2} - \delta_0 ] \cup [  \tfrac{2i-1}{2m+2} + \delta_0, \tfrac{2i-1}{2m+2} + \delta ) .     $$
For an index $i\in\{1,\ldots,m+1\}$ we introduce the notation $I_i =  (\tfrac{2i-1}{2m+2}-\delta, \tfrac{2i-1}{2m+2}+\delta)$ as well as $J_i = (\tfrac{2i-1}{2m+2}-\delta,\tfrac{2i-1}{2m+2}-\delta_0]\cup [\tfrac{2i-1}{2m+2}+\delta_0, \tfrac{2i-1}{2m+2}+\delta) $.
The cohomology of the pair $(U_{\Gamma_l},U_{\Gamma_l,\geq\epsilon_0})$ can be expressed as 
$$    \mathrm{H}^{\bullet}(U_{\Gamma_l},U_{\Gamma_l,\geq\epsilon_0}) \cong   \mathrm{H}^{\bullet}(\Gamma_l) \otimes \bigoplus_{i=1}^{m+1} \mathrm{H}^{\bullet}(I_i, J_i) .  $$
Since $(I_i,J_i)\simeq (I,\partial I)$ we have that
$$  \mathrm{H}^k (I_i,J_i) \cong \begin{cases}
    \mathbb{Q}, & k= 1\\
    \{0\}, & \text{else} .
\end{cases}   $$
Let $\chi_i\in\mathrm{H}^1(I_i,J_i)$ be the class such that the Kronecker pairing $\langle \chi_i, [I_i]\rangle = 1$ where $[I_i]\in\mathrm{H}_1(I_i, J_i)$ which is represented by the singular chain $I_i\colon [0,1]\to I_i$ defined by $I_i(t) = 2\delta_0 t + \tfrac{2i-1}{2m+2} -\delta_0$.
We see that the cohomology $\mathrm{H}^n(U_{\Gamma_l},U_{\Gamma_l,\geq\epsilon_0})$ is generated by the classes $\xi_r\times \chi_i$ for $r \in\{1,\ldots, 2m+1\}$ and $i\in \{1,\ldots, m\}$.
The map $f_l\times \mathrm{id}_I$ induces a map of pairs $(f_l\times \mathrm{id}_I)\colon (U_{\Gamma_l},U_{\Gamma_l,\geq\epsilon_0})  \to (U_a, U_{a,\geq\epsilon_0})$.
We want to express the pullback of the Thom class $\tau_a\in\mathrm{H}^n(U_a,U_{a,\geq\epsilon_0})$ by this map explicitly through the classes $\xi_r$ and $\chi_i$.
\begin{lemma}\label{lemma_pullback_thom_class}
    The pullback of the Thom class $\tau_a$ under the map $f_l\times\mathrm{id}_I$ satisfies
    $$  (f_l\times \mathrm{id}_I)^* \tau_a =   \sum_{i=1}^{m+1}  \xi_{2i-1}\times \chi_i .     $$
\end{lemma}
The proof follows the idea of the proof of \cite[Lemma 7.1]{stegemeyer2024string}.
\begin{proof}
    Since $(f_l\times \mathrm{id}_I)^*\tau_a \in \mathrm{H}^{n-1}(\Gamma_l)\otimes \bigoplus_{i=1}^{m+1} \mathrm{H}^1(I_i, J_i)$ there are numbers $c_{j,i}\in \mathbb{Q}$ such that 
    $$   (f_l\times \mathrm{id}_I)^*\tau_a =   \sum_{j=1}^l \sum_{i=1}^{m+1} c_{j,i} \xi_j \times \chi_i .       $$
    We want to show that $c_{j,i} = \delta_{ji}$, i.e. the Kronecker delta.
    First, we describe homology classes $X_k\in \mathrm{H}_{n-1}(\Gamma_l)$ which form a dual basis to the $\xi_j$.
    Fix a basepoint $(p_0,u_0)\in U\mathbb{S}^n$ and take $\mathbb{S}^{n-1}$ to be the equator in $\mathbb{S}^n$ orthogonal to $p_0$.
    For $j\in \{1,\ldots, l\}$ consider the map $s_j\colon \mathbb{S}^{n-1}\to \Gamma_l$ given by
    $$   s_j(v) = (p,u_0,\ldots, u_0, v , v,\ldots, v) \quad \text{for} \,\,\,v\in\mathbb{S}^{n-1} .       $$
    Here, the first entry $v$ comes at the $j$'th position where we recall that the first $u_0$ is at position $0$.
    We define $X_j = (s_j)_* [\mathbb{S}^{n-1}]\in \mathrm{H}_{n-1}(\Gamma_l)$.
    Then from the Gysin sequence of the iterated sphere bundle $\Gamma_l\to U\mathbb{S}^n$ we see that $\langle X_k,\xi_j\rangle = \delta_{kj}$.
    Furthermore, let $[I_i]\in \mathrm{H}_1(I_i,J_i)$ for $i\in\{1,\ldots, m+1\}$ be the generator which we defined above.
    By our choice of the class $X_j$ we obtain
    $$    c_{j,i} =  \langle X_j\times [I_i], (f_l\times\mathrm{id}_I)^* \tau_a  \rangle  =  \langle (f_l\times \mathrm{id}_I)_* (X\times [I_i]), \tau_a \rangle . $$
    Moreover, recall that $\tau_a = \mathrm{ev}_a^*\tau_{\mathbb{S}^n}$ and hence we get
    $$    c_{j,i} =   \langle   (\mathrm{ev}_a \circ ((f_l \circ s_j) \times \mathrm{id}_I))_* ([\mathbb{S}^n] \times [I_i]) , \tau_{\mathbb{S}^n}\rangle .    $$
    Denote the composition $\mathrm{ev}_a \circ ((f_l \circ s_j) \times \mathrm{id}_I)\colon \mathbb{S}^{n-1}\times I_i \to U_{\mathbb{S}^n}$ by $h'_{j,i}$.
    Then we see that $h'_{j,i}$ factors through a map
    $$     
    \begin{tikzcd}
        \mathbb{S}^{n-1}\times (I_i,J_i) \arrow[]{rr}{h'_{j,i}}\arrow[swap]{rd}{h_{j,i}} & & (U_{\mathbb{S}^n},U_{\mathbb{S}^n,\geq\epsilon_0}) 
        \\
        & (B_{-p_0}, B_{-p_0,\geq\epsilon_0}) \arrow[swap]{ru}{\iota} & 
    \end{tikzcd}
    $$
    where $B_{-p_0} =\{ q\in\mathbb{S}^n\,|\, \mathrm{d}(q,-p_0)< \epsilon\}$ and $B_{-p_0,\geq\epsilon_0} =  \{q\in B_{-p_0} \,|\, \mathrm{d}(q,-p_0)\geq \epsilon_0\}  $
    and $\iota$ is the obvious inclusion of pairs.
    Consequently we can express the number $c_{j,i}$ as 
    $$     c_{j,i} = \langle (h_{j,i})_* ([\mathbb{S}^{n-1}] \times [I_i] , \iota^* \tau_{\mathbb{S}^n} \rangle .    $$
    The pair $(B_{-p_0},B_{-p_0,\geq\epsilon_0})$ is homotopy equivalent to $(\mathbb{D}^n,\mathbb{S}^{n-1})$.
    and thus the defining property of the Thom class implies $\iota^*\tau_{\mathbb{S}^n}$ is the positively oriented generator of $\mathrm{H}^{n}(B_{-p_0},B_{-p_0,\geq\epsilon_0})\cong \mathrm{H}^n(\mathbb{D}^n,\mathbb{S}^{n-1})$.
    Consider the following diagram.
    $$
    \begin{tikzcd}
        \mathrm{H}_n(\mathbb{S}^{n-1}\times I_i, \mathbb{S}^{n-1}\times J_i) \arrow[]{r}{\partial} 
        \arrow[]{d}{(h_{j,i})_*} 
        &
        \mathrm{H}_{n-1}(\mathbb{S}^{n-1}\times J_i) \arrow[]{d}{(h_{j,i})_*}
        \\
        \mathrm{H}_n(\mathbb{D}^n,\mathbb{S}^{n-1}) \arrow[]{r}{\partial} & \mathrm{H}_{n-1}(\mathbb{S}^{n-1}) .
        \end{tikzcd}
    $$
    The horizontal arrows in this diagram are the connecting morphisms in the long exact sequence of a pair.
    The diagram commutes by the naturality of the long exact sequence of pairs.
    The lower horizontal arrow in this diagram is well-known to be an isomorphism.
    Note that $\mathrm{H}_{n}(\mathbb{S}^{n-1}\times I_i,\mathbb{S}^{n-1}\times J_i) \cong \mathbb{Q}$ and $$\mathrm{H}_{n-1}(\mathbb{S}^{n-1}\times J_i)  \cong \mathrm{H}_{n-1}(\mathbb{S}^{n-1}\times \{ \tfrac{2i-1}{2m+2}-\delta_0\}) \oplus \mathrm{H}_{n-1}(\mathbb{S}^{n-1}\times \{\tfrac{2i-1}{2m+2}+\delta_0\})  \cong \mathbb{Q}^2 .$$ 
    Under these identifications the connecting morphism can be written as $\partial (x) = (x,-y)$. 
    Finally, one sees from the definition that the restrictions
    $$   (h_{j,i})_* \colon \mathrm{H}_{n-1}(\mathbb{S}^{n-1}\times \{\tfrac{2i-1}{2m+2}+\delta_0\}) \to \mathrm{H}_{n-1}(\mathbb{S}^{n-1})      $$
    is an isomorphism if $j \geq 2i-1$ and $0$ otherwise.
    Moreover, we have that
    $$   (h_{j,i})_* \colon \mathrm{H}_{n-1}(\mathbb{S}^{n-1}\times \{\tfrac{2i-1}{2m+2}-\delta_0\}) \to \mathrm{H}_{n-1}(\mathbb{S}^{n-1})      $$
    is an isomorphism if $j \geq 2i$ and $0$ otherwise.
    Piecing everything together one hence finds that $c_{j,i} = \delta_{ji}$ which shows the claim.   
\end{proof}

\begin{theorem}\label{theorem_computation_copairing_even_sphere}
    Let $n\in\mathbb{N}$ be even. 
    For all $m\geq 0$ the antipodal copairing on $\mathbb{S}^n$ satisfies
    $$    \vee_a A_{2m+1}' =  \sum_{l=1}^{m}  A_{2l-2}\times A_{2m-2l}  \quad \text{and}\quad   \vee_a B_{2m+1}' =  \sum_{l=1}^m   A_{2l-2}\times B_{2m-2l} + \sum_{l=1}^m  B_{2l-2}\times A_{2m-2l}   .  $$
\end{theorem}
\begin{proof}
Set $l = 2m+1$.
Recall that we have
$$  A_l' =  (f_l)_* \circ (p_l)_! [p_0] \quad \text{and}\quad B_l' = (f_l)_* \circ (p_l)_! [U\mathbb{S}^n]  .     $$
It is clear by construction that $(p_l)_! [U\mathbb{S}^n] = [\Gamma_l]$ where the latter class is the orientation class of $\Gamma_l$.
We want to describe the class $(p_l)_![p_0]$.
Let $F \in\mathrm{H}_{l(n-1)}(\Gamma_l)$ be the homology class which is dual to the cohmology class $\xi_1\cup \ldots \cup \xi_l$. 
Using the cohomology ring of $\Gamma_l$ we find that $(p_l)_![p_0] =  - F$.
Since both classes $A_l'$ and $B_l'$ are in the image of the map $(f_l)_*\colon \mathrm{H}_i(\Gamma_l)\to \mathrm{H}_i(\Lambda \mathbb{S}^n)$ we consider the following diagram and claim that it commutes.
$$
\begin{tikzcd}
    \mathrm{H}_i(\Gamma_l\times I, \Gamma_l\times \partial I) \arrow[]{r}{(f_l\times \mathrm{id}_I)_*}
    \arrow[]{d}{(f_l\times \mathrm{id}_I)^*\tau_a \cap}
    &
    \mathrm{H}_i(\Lambda \mathbb{S}^n\times I,\Lambda \mathbb{S}^n\times \partial I)
    \arrow[]{d}{\tau_a \cap}
    \\
    \mathrm{H}_{i-n}(U_{\Gamma_l} ) \arrow[]{r}{(f_l\times \mathrm{id}_I)_*} \arrow[]{d}{k_*}
    &
    \mathrm{H}_{i-n}(U_a) \arrow[]{d}{(R_a)_*}
    \\
    \mathrm{H}_{i-n}(\Gamma_l\times \sqcup_{i=1}^{m+1} \{\tfrac{2i-1}{2m+2}\} ) \arrow[]{r}{(f_l)_*}
    &
    \mathrm{H}_{i-n}(F_a) .
\end{tikzcd}
$$
By naturality of the cap product the first square commutes.
In the second square the map $k\colon U_{\Gamma_l}\to \Gamma_l\times \sqcup_{i=1}^{m+1} \{\tfrac{2i-1}{2m+1}\}$ is defined by setting
$$    k(\gamma,s) =  (\gamma, \tfrac{2i-1}{2m+1}) \quad \text{if }s\in (\tfrac{2i-1}{2m+1} -\delta, \tfrac{2i-1}{2m+1} + \delta) .      $$
For the second square one can see that the underlying diagram of maps commutes up to homotopy and hence the lower square commutes as well.

Using Lemma \ref{lemma_pullback_thom_class} we make the following explicit computations.
For $i\in \{1,\ldots, m+1\}$ the cap product $  ( \xi_{2i-1}\times \chi_i ) \cap ([\Gamma_l] \times [I])      $ yields the dual homology class to the class $\alpha\cup \xi_1\ldots \xi_{2i-2}\cup \xi_{2i}\cup\ldots \cup \xi_{2m+1}$.
We denote this class by $G_{2i-1}\in\mathrm{H}_{\bullet}(\Gamma_l\times \{\tfrac{2i-1}{2m+1}\})$.
Similarly denote the dual homology class to the class $\xi_1\cup \ldots \cup \xi_{2i-2}\cup \xi_{2i}\cup\ldots \cup \xi_{2m+1}$ by $F_{2i-1}\in\mathrm{H}_{\bullet}(\Gamma_l\times \{ \tfrac{2i-1}{2m+1}\})$.
Then one finds that $(\xi_{2i-1}\times \chi_i)\cap (- F\times [I])=   F_{2i-1} $.
Thus, we get
\begin{equation}\label{eq_compputation_copairing_1}
     (R_a)_* (\tau_a \cap (A_l'\times [I]))  =   \sum_{i=1}^{m+1} (f_l)_* (F_{2i-1}) \times [\tfrac{2i-1}{2m+2}]   
    \end{equation}
    and
    \begin{equation}\label{eq_computation_copairing_2}
    (R_a)_* (\tau_a \cap (B_l'\times [I])) = \sum_{i=1}^{m+1} (f_l)_* ( G_{2i-1})\times [\tfrac{2i-1}{2m+2}] .
\end{equation}  

We are thus left with understanding the effect of the cutting map on the classes $(f_l)_* G_{2i-1}$ and $(f_l)_* F_{2i-1}$.
Note that we have maps $p_{l,k}\colon \Gamma_l\to U\mathbb{S}^n\times_{\mathbb{S}^n}U\mathbb{S}^n$ given by
$$  p_{l,k} ( p,v_0,\ldots, v_l) =  (p,v_0,v_{k}) .   $$
These are fiber bundles with fiber $(\mathbb{S}^{n-1})^{l-1}$.
The important point about this map is that the classes $F_{2i-1}$ and $G_{2i-1}$ are in the image of the Gysin map $(p_{l,2i-1})_!$.
Indeed we have that the cohomology ring of $ E= U\mathbb{S}^n\times_{\mathbb{S}^n}U\mathbb{S}^n$ is $    \mathrm{H}^{\bullet}(E) \cong \Lambda_{\mathbb{Q}}(\beta,\zeta)    $ with $|\beta| =2n-1$ and $|\zeta| = n-1$.
Moreover, we can choose the generators so that $p_{l,k}^* \beta = \alpha$ and $p_{l,k}^*\zeta = \xi_k$.
Furthermore note that the dual homology class to $\beta$ is the class $\Delta'_* [U\mathbb{S}^n]$ where $\Delta'\colon U\mathbb{S}^n\to U\mathbb{S}^n\times_{\mathbb{S}^n}U\mathbb{S}^n$ is the diagonal map.
Then one checks that $(p_{l,2i-1})_! [p_0] = F_{2i-1}$ and 
$(p_{l,2i-1})_! (\Delta'_* [U\mathbb{S}^n]  ) =  G_{2i-1}$.
Consider the following diagram.
\footnotesize
$$
\begin{tikzcd}
    \mathrm{H}_{i-{\lambda_l}+(n-1)}(U\mathbb{S}^n\times_{\mathbb{S}^n}U\mathbb{S}^n) \arrow[]{r}{(p_{l,2i-1})_!} 
    \arrow[]{d}{ (\mathrm{id}\times Da)_* }
    & [3.5em]
    \mathrm{H}_{i} (\Gamma_l \times \{\tfrac{2i-1}{2m+1}\}) \arrow[]{r}{(f_l\times \mathrm{id}_I)_*} 
    \arrow[]{d}{(\Phi_1^{-1})_*}
    & [3.5em]
    \mathrm{H}_i(F_a)
    \arrow[]{d}{(\mathrm{cut})_*}
    \\
    \mathrm{H}_{i-\lambda_l + (n-1)} (U\mathbb{S}^n\times_{\Delta_a}U\mathbb{S}^n) \arrow[]{r}{(p_{{2i-2}}\times p_{l-(2i-1)})_!} 
    \arrow[]{d}{(j_1)_*}
    &
    \mathrm{H}_{i}(\Gamma_{2i-2}\times_{\Delta_a} \Gamma_{l-(2i-1)} ) \arrow[]{r}{(f_{2i-2}\times f_{l-(2i-1)})_*}
    \arrow[]{d}{}
    &
    \mathrm{H}_i(P \hphantom{i}_{\mathrm{ev}_1}\times_{\mathrm{ev}_0} P) 
    \arrow[]{d}{}
    \\
    \mathrm{H}_{i-\lambda_l + (n-1)}(U\mathbb{S}^n\times U\mathbb{S}^n) \arrow[]{r}{(p_{2i-2}\times p_{l-(2i-1)})_!}
    &
    \mathrm{H}_i(\Gamma_{2i-2}\times \Gamma_{l-(2i-1)}) \arrow[]{r}{(f_{2i-2}\times f_{l-(2i-1)})_*}
    &
    \mathrm{H}_i(P_a\mathbb{S}^n\times P_a\mathbb{S}^n) .
\end{tikzcd}
$$
\normalsize
The vertical arrows in the lower squares are induced by the respective inclusions.
We claim that this diagram commutes. 
The upper left square commute since the diagram of maps
$$
\begin{tikzcd}
    \Gamma_{2i-2}\times_{\Delta_a} \Gamma_{l-(2i-1)} \arrow[]{d}{\Phi_1} \arrow[]{r}{p_{2i-2}\times p_{l-(2i-1)}}
    & [3em]
    U\mathbb{S}^n\times_{\Delta_a} U\mathbb{S}^n \arrow[]{d}{\mathrm{id}\times Da}
    \\
    \Gamma_l \arrow[]{r}{p_{l,2i-1}} & U\mathbb{S}^n\times_{\mathbb{S}^n} U\mathbb{S}^n
\end{tikzcd}
$$
commutes and since $\mathrm{id}\times Da \colon U\mathbb{S}^n \times_{\mathbb{S}^n}U \mathbb{S}^n \to U\mathbb{S}^n\times_{\Delta_a}U\mathbb{S}^n$ and $\Phi_1$ are orientation-preserving.
The commutativity of the upper right square can be seen directly from the definitions of the involved maps, see also Lemma \ref{completing_manifolds_diffeo}.
The lower left square commutes for the following reason.
Note that the fiber bundle $p_{2i-2}\times p_{l-(2i-1)}\colon \Gamma_{2i-2}\times_{\Delta_a}\Gamma_{l-(2i-1)}\to U\mathbb{S}^n\times_{\Delta_a}U\mathbb{S}^n$ is the pullback of the bundle $\Gamma_{2i-2}\times \Gamma_{l-(2i-1)}\to U\mathbb{S}^n\times U\mathbb{S}^n$ along the inclusion $U\mathbb{S}^n\times_{\Delta_a}U\mathbb{S}^n \hookrightarrow U\mathbb{S}^n\times U\mathbb{S}^n$.
The commutativity of the lower left square can therefore be seen as in the proof of \cite[Lemma 7.3]{stegemeyer2024string}.
Finally, the commutativity of the lower right square is clear.

In order to complete the proof we need to understand what happens if we start at the upper left corner of the above diagram with the classes $[p_0] \in \mathrm{H}_0(U\mathbb{S}^n\times_{\mathbb{S}^n}U\mathbb{S}^n)$ and $\Delta'_* [U\mathbb{S}^n]\in \mathrm{H}_{2n-1}(U\mathbb{S}^n\times_{\mathbb{S}^n}U\mathbb{S}^n)$.
For the class $[p_0]$ it is clear that $(j_1)_*(\mathrm{id}\times Da)_* [p_0] = [p_0] \times [p_0] \in\mathrm{H}_0(U\mathbb{S}^n\times U\mathbb{S}^n)$.
Then we see that $$  \mathrm{cut}_* \big( (f_l)_* (F_{2i-1}) \times [\tfrac{2i-1}{2m+2}]\big) =  (f_{2i-2}\times f_{l-(2i-1)})_* (p_{2i-2}\times p_{l-(2i-1)})_! ([p_0] \times [p_0])  . $$
By \cite[Proposition 14.3]{bredon:2013} we have
$$   (p_{2i-2}\times p_{l-(2i-1)})_! [p_0] \times [p_0]  =  (-1)^{(\mathrm{dim}(U\mathbb{S}^n)  + \mathrm{dim}(\Gamma_{2i-2}))\mathrm{dim}(U\mathbb{S}^n)} (p_{2i-2})_! [p_0] \times (p_{l-(2i-1)})_! [p_0]   $$
and since $\mathrm{dim}(U\mathbb{S}^n)  + \mathrm{dim}(\Gamma_{2i-2})$ is even we see that 
$$   \mathrm{cut}_* \big( (f_l)_* (F_{2i-1}) \times [\tfrac{2i-1}{2m+2}]\big)  = A_{2i-2}\times A_{l-(2i-1)} .    $$
By combining this with equation \eqref{eq_compputation_copairing_1} we obtain the expression for $\vee_a A'_{2m+1}$.
For the class $\Delta'_* [U\mathbb{S}^n ]$ we note that the composition
$$    U\mathbb{S}^n \times_{\mathbb{S}^n} U\mathbb{S}^n \xrightarrow[]{\mathrm{id}\times Da} U\mathbb{S}^n\times_{\Delta_a}U\mathbb{S}^n \xrightarrow[]{j_1}  U\mathbb{S}^n\times U\mathbb{S}^n $$
is equal to the composition 
$$      U\mathbb{S}^n\times_{\mathbb{S}^n}U\mathbb{S}^n\hookrightarrow U\mathbb{S}^n \times U\mathbb{S}^n  \xrightarrow[]{\mathrm{id}\times Da} U\mathbb{S}^n\times U\mathbb{S}^n  .  $$
Since $Da_* $ is the identity in homology in degree $2n-1$ it holds that 
$$   (j_1)_* (\mathrm{id}\times Da)_* \big( \Delta'_* [U\mathbb{S}^n]\big) =   \Delta_* [U\mathbb{S}^n ]      $$
where $\Delta\colon U\mathbb{S}^n\to U\mathbb{S}^n\times U\mathbb{S}^n$ is the ordinary diagonal map.
From the cohomology ring of $U\mathbb{S}^n$ we see directly that 
$$    \Delta_* [U\mathbb{S}^n] =    [U\mathbb{S}^n] \times [p_0]  + [p_0] \times [U\mathbb{S}^n]   .   $$
Again going along the lower horizontal arrows of the above diagram then yields the expression for the copairing $\vee_a' B_{2m+1}$.
This completes the proof.
\end{proof}

The comodule structures can be figured out with analogous methods.
We just note the result without proof.
\begin{prop}\label{prop_computation_comodule_even_sphere}
    Consider an even-dimensional sphere with the antipodal map $a\colon \mathbb{S}^n\to \mathbb{S}^n$ and take homology with rational coefficients.
    The left and right comodule pairings $\vee_{a,l}$ and $\vee_{a,r}$ behave as follows.
    For $m\in\mathbb{N}$ we have
    $$      \vee_{a,l} A_{2m} =  \sum_{i=1}^m     A'_{2m-(2i-1)} \times A_{2i-2}      $$
    as well as
    $$   \vee_{a,l} B_{2m} = - \sum_{i=1}^m    B'_{2m-(2i-1)} \times A_{2i-2} +   \sum_{i=1}^m    A'_{2m-(2i-1)} \times B_{2i-2 }  .   $$
    Furthermore, for $m\in\mathbb{N}$ it holds that
    $$   \vee_{a,r} A_{2m} = - \sum_{i=1}^m A_{2i-2}\times A'_{2m-(2i-1)}   $$
    and
    $$   \vee_{a,r} B_{2m} =- \sum_{i=1}^m  A_{2i-2}\times B'_{2m-(2i-1)}  - \sum_{i=1}^m   B_{2i-2}\times A'_{2m-(2i-1)} .     $$
\end{prop}

We end this section by describing the dual cohomology product.
For $m\in\mathbb{N}_0$ even let $\alpha_m\in\mathrm{H}^{\bullet}(P_a\mathbb{S}^n)$ be the class dual to $A_m$ and let $\beta_m\in\mathrm{H}^{\bullet}(P_a\mathbb{S}^n) $ be the class dual to $B_m$.
Furthermore for $k\in\mathbb{N}$ odd let $\alpha_k\in\mathrm{H}^{\bullet}(\Lambda M,M)$ be the class dual to $A_k'$ and let $\beta_k\in\mathrm{H}^{\bullet}(\Lambda M,M)$ be the class dual to $B_k'$.
Then by Theorem \ref{theorem_computation_copairing_even_sphere} and Proposition \ref{prop_computation_comodule_even_sphere} as well as by the computation of the Goresky-Hingston product for even-dimensional spheres in \cite{goresky:2009} we find that we have
$$     \alpha_i \overline{\ostar} \alpha_j =  \rho^a_{i,j}\alpha_{i+j+1}  , \quad \alpha_i\overline{\ostar}\beta_j =  \rho^b_{i,j} \beta_{i+j+1} ,  $$
$$    \beta_i \overline{\ostar}\alpha_j = \rho^c_{i,j} \beta_{i+j+1} \quad \text{and}\quad \beta_i\overline{\ostar} \beta_j = 0     $$
for $i,j\in\mathbb{N}_0$.
Here, $\rho^{a}_{i,j},\rho^b_{i,j},\rho^c_{i,j} \in \{1,-1\}$ for $i,j\in\mathbb{N}_0$.
In particular we note that the ring $(\mathrm{H}^{\bullet}(\Lambda M,M)\oplus \mathrm{H}^{\bullet}(P_a\mathbb{S}^n), \overline{\ostar})$ is generated by the classes $\alpha_0$ and $\beta_0$.
Recall from Section \ref{sec_computation_spheres} that the extended Chas-Sullivan product on even-dimensional spheres behaves very similarly to the Chas-Sullivan ring of odd-dimensional spheres.
We see now that a similar phenomenon happens for the extended cohomology product.
Recall from \cite[Theorem 15.3]{goresky:2009} that the rational Goresky-Hingston ring for an odd-dimensional sphere can be written as 
$$    (\mathrm{H}^{\bullet}(\Lambda \mathbb{S}^{n},\mathbb{S}^{n};\mathbb{Q}),\ostar_{\mathrm{GH}} ) \cong  \Lambda_{\mathbb{Q}} (U) \otimes \mathbb{Q}[T]_{\geq 2}   $$
where $\mathbb{Q}[T]_{\geq 2}$ means the ideal generated by the element $T^2$ in the polynomial ring $\mathbb{Q}[T]$.
The element $T^2$ is identified with the generator of the cohomology group $\mathrm{H}^{n-1}(\Lambda \mathbb{S}^n,\mathbb{S}^n)$.
See also Table \ref{table_GH} for a visualization.
Moreover note that for $i\geq 1$ the classes $T^{2i}, T^{2i+1}, UT^{2i}, UT^{2i+1}$ originate from the same critical manifold as explained in \cite{goresky:2009}.
\begingroup
\setlength{\tabcolsep}{7pt} 
\renewcommand{\arraystretch}{1.3} 
\begin{table}[h!]
\centering
\begin{tabular}{|c | c c c c c c c c c|} 
 \hline
   degree & $0$ & $n-1$ & $n$ & $2n-2$ & $2n-1$ & $3n-3$ & $3n-2$ & $4n-4$ & $4n-3$   \\  
 \hline
 $\mathrm{H}_{\bullet}(\Lambda \mathbb{S}^n,\mathbb{S}^n)$  &   &  $T^2$   &  &  $T^3$  &   $UT^2 $ &  $T^4$ & $UT^3$ &  $T^5$  & $UT^4$  \\ 
 \hline
\end{tabular}
\caption{Goresky-Hingston product for an odd-dimensional sphere.}
\label{table_GH}
\end{table}
\endgroup

For the extended Goresky-Hingston algebra for even-dimensional spheres we have seen above that the generator $\alpha_1\in\mathrm{H}^{n-1}(\Lambda \mathbb{S}^n,\mathbb{S}^n)$ is itself a power of the class $\alpha_0\in\mathrm{H}^0(P_a\mathbb{S}^n)$, i.e. we have $\alpha_1 = \pm \alpha_0^2$.
For $i\geq 0$ define $\alpha_i = S^{i+1}$ and $\beta_i = WS^{i+1}$.
From the above we find that we can - up to signs - write the extended Goresky-Hingston ring as 
$$    \Lambda_{\mathbb{Q}}(W) \otimes \mathbb{Q}[S]_{\geq 1}    .  $$
Here, $\mathbb{Q}[S]_{\geq 1}$ is the ideal in the polynomial ring $\mathbb{Q}[S]$ generated by the element $S$.
For each $i\in\mathbb{N}$ the elements $S^i$ and $WS^i$ are induced from the same critical manifold.
The ideal generated by $S^2$ and $WS^2$ is the Goresky-Hingston ring of the even-dimensional sphere.
See Table \ref{table_extended_GH} for a visualization of the extended cohomology product.

\begingroup

\setlength{\tabcolsep}{7pt} 
\renewcommand{\arraystretch}{1.3} 
\begin{table}[h!]
\centering
\begin{tabular}{|c | c c c c c c c c c|} 
 \hline
   degree & $0$ & $n-1$ & $n$ & $2n-2$ & $2n-1$ & $3n-3$ & $3n-2$ & $4n-4$ & $4n-3$   \\  
 \hline
 $\mathrm{H}_{\bullet}(\Lambda \mathbb{S}^n,\mathbb{S}^n)$ &   &  $S^2$   &  &   &  &  $S^4$ & $WS^2$ &   &  \\ 
 \hline
 $\mathrm{H}_{\bullet}(P_a\mathbb{S}^n)\hphantom{bli} $ &  $S$  &  &  & $S^3$ & $WS$ &  & & $S^5$  & $ WS^3$ \\
 \hline
\end{tabular}
\caption{Extended Goresky-Hingston product for an even-dimensional sphere with rational coefficients.}
\label{table_extended_GH}
\end{table}

\endgroup

\section{A resonance theorem for closed geodesics on real projective space}\label{sec_resonance_theorem}

In this final section we apply the computations of the homology and cohomology products for even-dimensional spheres to obtain a resonance theorem for closed geodesics on real projective space.
We begin by introducing critical values of homology and cohomology classes and investigate how they behave with respect to the products.
We then prove the resonance theorem following the ideas of Hingston and Rademacher \cite{hingston:2013}.
For discussions of the resonance theorem and further applications of the resonance theorem in the case of spheres we refer to \cite{hingston:2013}.

We begin by considering a general situation.
Let $X$ be a Hilbert manifold and let $f\colon X\to \mathbb{R}$ a $C^1$-function which is bounded from below.
Let $A\in\mathrm{H}_{\bullet}(X)$ be a homology class where we take homology with coefficients in a commutative unital ring $R$.
We define the \textit{critical value of} $A$ to be
$$    \mathrm{cr}(A)  =  \mathrm{inf}\{   a\in \mathbb{R}\,|\,   A\in \mathrm{im}\big( \mathrm{H}_{\bullet}( X^{\leq a}) \to \mathrm{H}_{\bullet}(X)  \big) \}         $$
where the map $\mathrm{H}_{\bullet}(X^{\leq a}) \to \mathrm{H}_{\bullet}(X)$ is induced by the inclusion of the sublevel set $X^{\leq a }\hookrightarrow X$.
Similarly, in cohomology, we consider a class $\varphi\in\mathrm{H}^{\bullet}(X)$ and define the \textit{critical value of} $\varphi$ by 
$$   \mathrm{cr}(\varphi) =  \mathrm{sup}\{ a\in \mathbb{R} \,|\,   A \in\mathrm{ker}\big( \mathrm{H}^{\bullet}(X)\to \mathrm{H}^{\bullet}(X^{\leq a})\big) \} .    $$
Goresky and Hingston show in \cite[Lemmas 4.2 and 4.4]{goresky:2009} that the critical value of a homology class, resp. of a cohomology class is indeed a critical value of $f$.

Let $M$ be a closed manifold and let $f\colon M\to M$ be an involution.
Fix a Riemannian metric on $M$.
We now consider critical values of homology and cohomology classes in the free loop space, respectively in the path space $P_f M$.
Recall from Section \ref{sec_antipodal_geodesics} that we have the energy functional $E\colon PM\to [0,\infty)$ which restricts to the energy functional on $\Lambda M$ and $P_f M$, respectively.
Define the function $\mathcal{L}\colon PM\to [0,\infty)$ by $\mathcal{L}(\gamma) = \sqrt{E(\gamma)}$ for $\gamma\in PM$.
We consider this function instead of the energy functional since it behaves better under concatenation of paths.
Note that $\mathcal{L}$ is not smooth on $PM$ but it is smooth away from the constant loops and thus behaves well enough for our purposes.
From now on we always take the critical values of homology and of cohomology classes with respect to the function $\mathcal{L}$.

Recall from \cite{goresky:2009} that the Chas-Sullivan product and the Goresky-Hingston product satisfy the following inequalities.
For $X,Y\in\mathrm{H}_{\bullet}(\Lambda M)$ and $\varphi,\psi\in\mathrm{H}^{\bullet}(\Lambda M,M)$ we have
$$   \mathrm{cr}(X\wedge_{\mathrm{CS}} Y) \leq \mathrm{cr}(X) + \mathrm{cr}(Y) \quad \text{and}\quad \mathrm{cr}(\varphi\ostar_{\mathrm{GH}} \psi) \geq \mathrm{cr}(\varphi) + \mathrm{cr}(\psi) .      $$
Moreover, for the $\mathbb{S}^1$-action we have $\mathrm{cr}(B(X)) \leq \mathrm{cr}(X)$ for all $X\in\mathrm{H}_{\bullet}(\Lambda M)$.
We will now show that the analogous inequalities hold for the operations that we are considering in this paper.
\begin{prop}\label{prop_morse_theoretic_inequalities}
    Let $M$ be a closed oriented manifold and $f\colon M\to M$ a smooth involution.
    Take homology with coefficients in a commutative and unital ring $R$.
    Choose a Riemannian metric $g$ on $M$ and consider the length functional $\mathcal{L}\colon \Lambda M\to [0,\infty)$, respectively $\mathcal{L}\colon P_f M\to [0,\infty)$ with respect to $g$.
    \begin{enumerate}
        \item Let $X\in\mathrm{H}_{\bullet}(\Lambda M)$ and $A,B\in\mathrm{H}_{\bullet}(P_f M)$.
        Then 
        $$     \mathrm{cr}(A\wedge_f B) \leq \mathrm{cr}(A) + \mathrm{cr}(B) , \quad \mathrm{cr}(X*_l A) \leq \mathrm{cr}(X) + \mathrm{cr}(A) $$ and $$ \mathrm{cr}(A*_r X) \leq \mathrm{cr}(A) + \mathrm{cr}(X) .      $$
        \item Assume that $R$ is a field and let $\varphi\in\mathrm{H}^{\bullet}(\Lambda M,M)$ and $\alpha,\beta\in\mathrm{H}^{\bullet}(P_f M)$.
        Then 
        $$    \mathrm{cr}(\alpha\ostar_f\beta ) \geq \mathrm{cr}(\alpha) + \mathrm{cr}(\beta) ,\quad  \mathrm{cr}(\varphi \ostar_l \alpha) \geq \mathrm{cr}( \varphi)  + \mathrm{cr}(\alpha)      $$
        and 
        $$     \mathrm{cr}(\alpha \ostar_r \varphi) \geq \mathrm{cr}(\alpha)  
        + \mathrm{cr}(\varphi)     .   $$
    \end{enumerate}
\end{prop}
\begin{proof}
    Let $a>0$ be the minimum of $\mathcal{L}\colon P_f M\to [0,\infty)$ and let $b,c\in\mathbb{R}$ be numbers with $b,c\geq  a$.
    We consider the sublevel sets 
    $  P^{\leq b}  =  (P_f M)^{\leq b} $ and $ P^{\leq c} = (P_f M)^{\leq c}  $.
    We want to define a map $\wedge_{f;b,c}\colon \mathrm{H}_{\bullet}(P^{\leq b})\otimes \mathrm{H}_{\bullet}(P^{\leq c}) \to \mathrm{H}_{\bullet}(P^{\leq b+c+2 \epsilon})$ such that the diagram 
    \begin{equation}\label{eq_diagram_sublevel}
    \begin{tikzcd}
        \mathrm{H}_{\bullet}(P^{\leq b}) \otimes \mathrm{H}_{\bullet}(P^{\leq c}) \arrow[]{r}{} \arrow[]{d}{\wedge_{f;b,c}} 
    &
    \mathrm{H}_{\bullet}(P_f M)\otimes \mathrm{H}_{\bullet}(P_f M)\arrow[]{d}{\wedge_f}
    \\
    \mathrm{H}_{\bullet}(\Lambda M^{\leq b+c+2 \epsilon}) \arrow[]{r}{}
    &
    \mathrm{H}_{\bullet}(\Lambda M) 
    \end{tikzcd}
    \end{equation}
    commutes where the vertical arrows are induced by the inclusions.
    Here, the number $\epsilon >0$ is the one which we chose for the tubular neighborhood $U_M$, see Section \ref{sec_extensions_of_coproduct}.
    Note that $\wedge_f$ is independent of the choice of $\epsilon$.
    From Section \ref{sec_involution_path_space_string_operations} we recall the space $C_f\subseteq P_f M\times P_f M$ as well as the tubular neighborhood $U_C \subseteq P_f M\times P_f M$ of $C_f$.
    For the definition of $\wedge_{f;b,c}$ we define the spaces
    $$ U_{C;b,c} =  U_C \cap  (P^{\leq b}\times P^{\leq c}) \quad \text{and}\quad C^f_{b,c} = C_f \cap  (P^{\leq b}\times P^{\leq c}) .   $$
    There is an inclusion map $j\colon (U_{C;b,c},U_{C;b,c}\setminus C^f_{b,c})\hookrightarrow (U_C,U_C\setminus C_f)$.
    Moreover we define the \textit{optimal concatenation map} $\phi\colon C_f\to \Lambda M$ which is given by
    $$    \phi(\gamma,\sigma)(t) =   \begin{cases}
            \gamma(\!\frac{t}{\tau}) , & 0\leq t \leq \tau \\
            \sigma(\!\frac{t-\tau}{1-\tau}), & \tau\leq t\leq 1
    \end{cases}
     $$
     for $\tau = \frac{\mathcal{L}(\gamma)}{\mathcal{L}(\gamma)+\mathcal{L}(\sigma)}$.
     One checks that for $(\gamma,\sigma)\in C_f$ we have
     $   \mathcal{L}(\phi(\gamma,\sigma)) = \mathcal{L}(\gamma)  + \mathcal{L}(\sigma)    $.
     Moreover, $\phi$ is homotopic to the ordinary concatenation map $\mathrm{concat}\colon C_f \to \Lambda M$.
     Similarly one can define the optimal concatenation of several paths.
        Recall that the retraction $R^f\colon U_{C}\to C_f$ is defined as 
        $$   R^f(\gamma,\sigma) =  (\gamma, \mathrm{concat}(\overline{\gamma(1)\sigma(0)},\sigma, \overline{\sigma(1)\gamma(0)}) \quad \text{for}\,\,\,(\gamma,\sigma)\in U_{C} .     $$
    We define a map $\widetilde{R}^f\colon U_{C}\to C_f$ by defining it like $R^f$ but we use optimal concatenation instead of the usual concatenation.
    One checks that this induces a map $\widetilde{R}^f\colon U_{C;b,c}\to C^f_{b,c+2\epsilon}$.
     We define $\wedge_{f;b,c}$ as the composition
      \begin{eqnarray*}
        \wedge_{f;b,c} \colon  \mathrm{H}_{i}(P^{\leq b})\otimes \mathrm{H}_j(P^{\leq c}) &\xrightarrow{  (-1)^{n(n-i)}\times   }& \mathrm{H}_{i+j}(P^{\leq b}\times P^{\leq c})
     \\ 
     &\xrightarrow{ \hphantom{coninci}\vphantom{r}\hphantom{conci} }& \mathrm{H}_{i+j}(P^{\leq b}\times P^{\leq c},P^{\leq b}\times P^{\leq c}\setminus C^f_{b,c}) 
     \\
    &\xrightarrow{\hphantom{cii}\text{excision}\hphantom{ici}}&
    \mathrm{H}_{i+j}(U_{C;b,c}, U_{C;b,c}\setminus C^f_{b,c}) \\
    &\xrightarrow{\hphantom{coci} j^*\tau_{C}\cap \hphantom{coci} }& \mathrm{H}_{i+j-n}(U_{C;b,c})
    \\ &\xrightarrow[]{\hphantom{conci}\widetilde{R}^f_* \hphantom{conci}} & \mathrm{H}_{i+j-n}(C^f_{b,c+2\epsilon}) 
    \\ & \xrightarrow[]{\hphantom{ici}\mathrm{concat}_*\hphantom{ici}} & \mathrm{H}_{i+j-n}(\Lambda M^{ \leq b+c+2\epsilon}) \, .
    \end{eqnarray*}

    It is clear with this definition that the diagram \eqref{eq_diagram_sublevel} commutes and since the number $\epsilon$ can be chosen arbitrarily small the inequality $\mathrm{cr}(A\wedge_f B)\leq \mathrm{cr}(A) + \mathrm{cr}(B)$ follows.
    The other two inequalities in part (1) involving the module structures can be shown in an analogous manner.

    For the second part we consider the set 
    $$    (P\times P)^{\leq b+2\epsilon} :=  \{ (\gamma,\sigma)\in P_f M\times P_f M\,|\,  \mathcal{L}(\gamma) + \mathcal{L}(\sigma) \leq b+2\epsilon \}  .     $$
    We claim that the antipodal pairing factors as follows
    $$   
    \begin{tikzcd}
        \mathrm{H}_{\bullet}(\Lambda M^{\leq b}) \arrow[swap]{d}{\vee_f'}
        \arrow[]{dr}{\vee_f}
        & 
        \\
        \mathrm{H}_{\bullet}((P\times P)^{\leq b+2\epsilon}) \arrow[]{r}{} & \mathrm{H}_{\bullet}(P_f M\times P_f M) 
    \end{tikzcd}
    $$ 
    where the lower horizontal arrow is induced by the inclusion $(P\times P)^{\leq b +2\epsilon} \hookrightarrow P_f M\times P_f M$.
    Define the spaces $U_{f,b} := U_f\cap (\Lambda M^{\leq b}\times I)$ and $U_{f,b,\geq\epsilon_0} := U_{f,\geq\epsilon_0}\cap (\Lambda M^{\leq b}\times I)$.
    There is an inclusion $j\colon (U_{f,b}, U_{f,b,\geq\epsilon_0})\to  (U_f,U_{f,\geq\epsilon_0})$.
    Moreover, recall that if $(\gamma,s)\in U_f$ then we have
    $$  R_f (\gamma,s) = \big( \mathrm{concat}^s (\mathrm{concat}(\gamma|_{[0,s]},\overline{\gamma(s)\gamma(0)}), \mathrm{concat}( \overline{\gamma(0)\gamma(s)}, \gamma_{[s,1]}))  ,s\big)  $$
    where $\mathrm{concat}^s$ means that the first path is run through during the interval $[0,s]$ and the second one during the interval $[s,1]$. 
    This is uniquely defined such that the reparametrizations are affine linear maps.
    We define a map $\widetilde{R}_f\colon U_f\to F_f$ by setting
    $$     \widetilde{R}_f (\gamma,s) = \big( \mathrm{concat}^s (\phi(\gamma|_{[0,s]},\overline{\gamma(s)f(\gamma(0))}), \phi( \overline{f(\gamma(0))\gamma(s)}, \gamma_{[s,1]}))  ,s\big)     $$
    for $(\gamma,s)\in U_f$ where $\phi$ means again optimal concatenation.
    This map is homotopic to $R_f$.
    Moreover, if we have $(\gamma,s)\in F_f$ such that $\mathcal{L}(\gamma)\leq b$, then write $\mathrm{cut}(\gamma,s) = (\sigma,\eta)$.
    One checks that
    $$   \mathcal{L}(\sigma) = \sqrt{s\int_0^s |\Dot{\gamma}(u)|^2 \,\mathrm{d}s}  \quad \text{and}\quad  \mathcal{L}(\eta)  = \sqrt{(1-s) \int_{s}^1 |\Dot{\gamma}(u)|^2 \,\mathrm{d}s}   .   $$
    It is elementary to show that for numbers $a,b\in [0,\infty)$ and $\tau\in[0,1]$ it holds that $\sqrt{\tau a}  + \sqrt{(1- \tau)b }\leq \sqrt{a+b}$ and thus
    $   \mathcal{L}(\sigma) + \mathcal{L}(\eta) \leq \mathcal{L}(\gamma)     $.
    Consequently, the cutting map maps $F_{f,b} = F_f \cap (\Lambda M^{\leq b}\times I)$ to $(P\times P)^{\leq b}$.
    Putting everything together we obtain a map $\vee'_f\colon \mathrm{H}_{\bullet}(\Lambda M^{\leq b})\to \mathrm{H}_{\bullet}((P\times P)^{\leq b+2\epsilon})$ and the diagram above commutes as claimed.   
    It then follows that if we consider critical values on $P_f M\times P_f M$ with respect to the function 
    $ P_f M \times P_f M\to [0,\infty) $, $ (\gamma,\sigma)\mapsto \mathcal{L}(\gamma) +\mathcal{L}(\sigma)    $,
  then $\mathrm{cr}(\vee_f X) \leq \mathrm{cr}(X)$ for all $X\in\mathrm{H}_{\bullet}(\Lambda M)$.
   By dualizing we obtain the inequality
   $$    \mathrm{cr}(\alpha\ostar_f \beta) \geq \mathrm{cr}(\alpha) + \mathrm{cr}(\beta)  $$
   for cohomology classes $\alpha,\beta\in\mathrm{H}^{\bullet}(P_f M)$.
   The other inequalities can be shown similarly.
\end{proof}

\begin{example}
	Consider the standard metric on $\mathbb{S}^n$ for $n$ even.
	From the construction in Section \ref{sec_completing_manifolds} it follows that we have the following critical values on $\Lambda \mathbb{S}^n$ and $P_a\mathbb{S}^n$.
	Following the notation in Section \ref{sec_computation_spheres} we have for $l\in\mathbb{N}$ even and $m\in\mathbb{N}$ odd
	$$   \mathrm{cr}(A_l) = \mathrm{cr}(B_l) =  (l+1)\pi \quad \text{and}\quad \mathrm{cr}(A_m') = \mathrm{cr}(B_m' ) = (m+1) \pi .$$
	From the behavior of the extended loop product on even-dimensional spheres, see Section \ref{sec_computation_spheres}, we see that we actually have
	$     \mathrm{cr}(X \wedge Y  ) = \mathrm{cr}(X) + \mathrm{cr}(Y)       $
	for all $X,Y\in\widehat{\mathrm{H}}_{\bullet}(\mathbb{S}^n)$.
	In general, i.e. for an arbitrary metric, it will not be true that equality always holds.
	
	Already for odd-dimensional spheres with the standard metric - in the notation of Corollary \ref{cor_extendedn_product_sphere} - we have for example that
	$$   \mathrm{cr}(E\wedge_{a} E) = 0 < 2\, \mathrm{cr}(E) . 
	$$
    We have $\mathrm{cr}(E) > 0$ for any metric since the energy functional on $P_a\mathbb{S}^n$ is bounded from below by the shortest antipodal geodesic.
	Note moreover, that already for the Chas-Sullivan product itself and for the standard metric, it holds that for some products the critical values do not add up, e.g. we have
	$    \mathrm{cr}(U\wedge_{\mathrm{CS}} U) =  \mathrm{cr}(U) < 2 \mathrm{cr}(U),  $
	see \cite{goresky:2009}.
	 Note however that in the subalgebra generated by the classes $B$ and $V$ however, see Remark \ref{remark_non-trivial_subalgebra}, we have that for all $X,Y\in (B,V)$ it holds that $\mathrm{cr}(X\wedge Y) = \mathrm{cr}(X) + \mathrm{cr}(Y)$.
\end{example}

Recall from Corollary \ref{cor_antipodal_universal_covering_space} that the space of antipodal paths is the universal covering space of the component of non-contractible loops of the free loop space $\Lambda \mathbb{R}P^n$.
As before we denote this component by $\Lambda_1\mathbb{R}P^n$ and we denote the component of contractible loops by $\Lambda_0\mathbb{R}P^n$.
In particular the map $\pi_1\colon P_a\mathbb{S}^n\to \Lambda_1\mathbb{R}P^n$ given by
$\pi_1(\gamma)(t) =  \mathrm{pr}(\gamma(t))$ where $\mathrm{pr}\colon \mathbb{S}^n\to \mathbb{R}P^n$ is the standard $2$-fold covering is this universal covering map for $n\geq 3$, see Corollary \ref{cor_antipodal_universal_covering_space}.
Similarly, the map $\pi_0\colon \Lambda \mathbb{S}^n\to \Lambda_0 \mathbb{R}P^n$ given by $\pi_0(\gamma)(t) = \mathrm{pr}(\gamma(t))$ is the universal covering map.
Furthermore note that any Riemannian metric $g$ on $\mathbb{R}P^n$ induces a unique $\mathbb{Z}_2$-invariant metric on $\mathbb{S}^n$ such that the covering map $\mathrm{pr}\colon \mathbb{S}^n\to \mathbb{R}P^n$ becomes a Riemannian covering.
Denote this metric on $\mathbb{S}^n$ by $\widetilde{g}$.
We now take the length functional $\mathcal{L}_g \colon \Lambda \mathbb{R}P^n\to [0,\infty)$ on the free loop space of real projective space and we take the length functional $\mathcal{L}_{\widetilde{g}}\colon \Lambda \mathbb{S}^n\to [0,\infty)$, respectively $\mathcal{L}_{\widetilde{g}}\colon P_a\mathbb{S}^n\to [0,\infty)$.
From elementary Riemannian geometry it is clear that if $\gamma\in P_a\mathbb{S}^n$ is an antipodal geodesic then $\pi_1(\gamma)$ is a closed geodesic in $\mathbb{R}P^n$.
Analogously, $\pi_0$ maps closed geodesics on $\mathbb{S}^n$ to contractible closed geodesics on $\mathbb{R}P^n$.
Conversely, every closed geodesic on $\mathbb{R}P^n$ is the image of a critical point in $\Lambda \mathbb{S}^n$ or $P_a\mathbb{S}^n$ under the map $\pi_0$ or $\pi_1$, respectively.
Moreover, in Lemma \ref{lemma_index_and_null_preserved_under_covering} we have seen explicitly that also the Morse index is preserved, i.e. if $\gamma\in P_a\mathbb{S}^n$ is an antipodal geodesic of index $k$, then $\pi_1(\gamma)\in \Lambda \mathbb{R}P^n$ is a closed geodesic of index $k$.
Hence, the critical values of $\mathcal{L}_{\widetilde{g}}\colon \Lambda \mathbb{S}^n\sqcup P_a\mathbb{S}^n \to [0,\infty)$ are in one-to-one correspondence with the critical values of $\mathcal{L}_g\colon \Lambda \mathbb{R}P^n \to [0,\infty)$.

Following \cite{hingston:2013} we shall now apply the Morse-theoretic inequalities of Proposition \ref{prop_morse_theoretic_inequalities} to study the behavior of how the critical value of a homology or cohomology class is related to the degree of this class.
Note that for a closed Riemannian manifold $M$ the degree of a homology class $X\in\mathrm{H}_{\bullet}(\Lambda M)$ is related to the index of a closed geodesic at the level $\mathrm{cr}(X)$, see \cite{hingston:2013}. 
More precisely, if $X\in\mathrm{H}_{\bullet}(\Lambda M)$ for a closed manifold $M$ with $\mathrm{cr}(X) = a$ then there is a closed geodesic $\gamma\in \Lambda M$ at level $a$ with
$$   \mathrm{ind}(\gamma) \leq \mathrm{deg}(X) \leq \mathrm{ind}(\gamma) + \mathrm{null}(\gamma) + 1 \leq \mathrm{ind}(\gamma) + 2n-1 .    $$
Hence, for each homology classes $X\in\mathrm{H}_{\bullet}(\Lambda M)$ the pair $(\mathrm{cr}(X),\mathrm{deg}(X))$ gives us a closed geodesic at level $\mathrm{cr}(X)$ with index approximately equal to $\mathrm{deg}(X)$.
In our case we study critical values of classes in $\Lambda \mathbb{S}^n\sqcup P_a\mathbb{S}^n$ for a $\mathbb{Z}_2$-invariant metric on $\mathbb{S}^n$. 
As we discussed above the pair $(\mathrm{cr}(X),\mathrm{deg}(X))$ for $X\in\mathrm{H}_{\bullet}(\Lambda \mathbb{S}^n)\oplus \mathrm{H}_{\bullet}(P_a\mathbb{S}^n)$ gives us the length of a closed geodesic in $\mathbb{R}P^n$ and the degree of this closed geodesic.
Therefore if we find that the points $(\mathrm{cr}(X),\mathrm{deg}(X))\in \mathbb{R}_{\geq 0}\times \mathbb{R}_{\geq 0}$ lie only in a certain area, e.g. a strip, then this gives us information about the closed geodesics in $\mathbb{R}P^n$.
We will see exactly this behavior, i.e. that the pairs $(\mathrm{cr}(X),\mathrm{deg}(X))$ lie in a strip, in the following resonance theorem.
First, we need to recall some facts from \cite[Section 3]{hingston:2013} about critical values.
\begin{lemma}\label{lemma_basics_critical_values}
    Let $X$ be a Hilbert manifold and $F\colon X\to [0,\infty)$ a Morse-Bott function.
    Take homology with coefficients in a field $\mathbb{K}$.
    \begin{enumerate}
        \item Let $A\in\mathrm{H}_{\bullet}(X)$ be a homology class and $\lambda \in\mathbb{K}\setminus \{0\}$ then $\mathrm{cr}(\lambda A) = \mathrm{cr}( A)$.
        \item Let $\varphi\in\mathrm{H}^{\bullet}(X)$ and $\lambda\in \mathbb{K}\setminus \{0\}$ then $\mathrm{cr}(\lambda \varphi) = \mathrm{cr}(\varphi)$.
        \item Assume that $\mathrm{H}_i(X) \cong \mathbb{K}$ and let $A\in\mathrm{H}_i(X)$ be a generator.
        If $\varphi\in\mathrm{H}^i(X)$ is a cohomology class such that $\varphi(A) \neq 0$ then $\mathrm{cr}(\varphi)  = \mathrm{cr}(A)$.
    \end{enumerate}
\end{lemma}

\begin{theorem}\label{theorem_resonance_even}
    Let $g$ be a Riemannian metric on $\mathbb{S}^n$ for $n$ even and consider the induced energy functional on $\Lambda\mathbb{S}^n \sqcup P_a\mathbb{S}^n$.
    Consider the direct sum $\widehat{\mathrm{H}}_{\bullet}(\mathbb{S}^n) = \mathrm{H}_{\bullet}(\Lambda \mathbb{S}^n;\mathbb{Q}) \oplus \mathrm{H}_{\bullet}(P_a\mathbb{S}^n;\mathbb{Q})$.
    There are constants $\overline{\alpha},\beta > 0$ such that for all $X\in\widehat{\mathrm{H}}_{\bullet}(\mathbb{S}^n)$ we have
    $$   | \mathrm{deg}(X) - \overline{\alpha}\, \mathrm{cr}(X)  | \leq \beta .      $$
\end{theorem}
\begin{proof}

Note that every class $X\in\widehat{\mathrm{H}}_{\bullet}(\mathbb{S}^n)$ can be written as a product $X = \lambda A_0^{r} \wedge [p_0]^s \wedge B_0^t $ for $\lambda\in\mathbb{Q}$, $r,s\in\{0,1\}$ and $t\in \mathbb{N}_0$.
By Proposition \ref{prop_morse_theoretic_inequalities} we therefore get
$$   \mathrm{cr}(X) \leq r \mathrm{cr}(A_0) +  \mathrm{cr}(B_0^t)      $$
since $\mathrm{cr}([p_0]) = 0$ for every metric.
Moreover, every class $\varphi\in \widehat{\mathrm{H}}^{\bullet}(\mathbb{S}^n) = \mathrm{H}^{\bullet}(\Lambda \mathbb{S}^n,\mathbb{S}^n)\oplus \mathrm{H}^{\bullet}(P_a\mathbb{S}^n)$ can be written as $\varphi = \mu \beta_0^y \overline{\ostar}\alpha_0^z$ for $\mu\in\mathbb{Q}$, $y\in\{0,1\}$ and $z\in\mathbb{N}$.
Hence, Proposition \ref{prop_morse_theoretic_inequalities} implies that
$$  \mathrm{cr}(\varphi) \geq   y \mathrm{cr}(\beta_0) +  \mathrm{cr}(\alpha_0^z ) \geq \mathrm{cr}(\alpha_0^z) .    $$
If $\varphi(X) \neq 0$ we have $\mathrm{cr}(X) = \mathrm{cr}(\varphi)$ and we see from the explicit computations in Sections \ref{sec_computation_spheres} and \ref{sec_extensions_of_coproduct} that $z = t +1$. 
Therefore if we write $z = 2k+\ell$ with $\ell\in \{0,1\}$ we get a bound of the form
$$   \mathrm{cr}(\alpha_0^{2(k-1)}) \leq \mathrm{cr}(\alpha_0^{t+1})  \leq \mathrm{cr}(X) \leq r \mathrm{cr}(A_0) + 
\ell \mathrm{cr}(B_0) + \mathrm{cr}(B_0^{2k})  .  $$
By \cite[Proof of Theorem 1.1]{hingston:2013} we have an inequality
\begin{equation}\label{eq_inequality_using_s1_action}
    \mathrm{cr}(B_0^{2k}) \leq \mathrm{cr}(\alpha_0^{2(k+2)}) . 
\end{equation}   
The rest of the argument then follows as in the proof of \cite[Theorem 1.1]{hingston:2013}.  
\end{proof}
\begin{remark}
    Consider the situation of the above resonance theorem.
    \begin{enumerate}
        \item 
    One of the crucial points in the proof of \cite[Theorem 1.1]{hingston:2013} is the use of equation \eqref{eq_inequality_using_s1_action}.
    This inequality comes from the $\mathbb{S}^1$-action on the free loop space and its compatibility with critical values.
    As we have discussed in Remark \ref{remark_bv_algebra} there is also an $\mathbb{S}^1$-action on $P_a\mathbb{S}^n$ and an induced operator $B\colon \mathrm{H}_{\bullet}(P_a\mathbb{S}^n)\to \mathrm{H}_{\bullet+1}(P_a\mathbb{S}^n)$.
    It is not hard to see that we then also have $\mathrm{cr}(B(X))\leq \mathrm{cr}(X)$ for $X\in\mathrm{H}_{\bullet}(P_a\mathbb{S}^n)$ if the metric on $\mathbb{S}^n$ is $\mathbb{Z}^2$-invariant.
    In order to use this inequality for the proof of the above theorem we would however need an explicit computation of the operator $B$ for $\mathrm{H}_{\bullet}(P_a\mathbb{S}^n)$ which we did not include in this paper.
    In the proof above however we only use the inequality $\mathrm{cr}(B(X))\leq \mathrm{cr}(X)$ for classes $\mathrm{H}_{\bullet}(\Lambda \mathbb{S}^n)$ in the free loop space where the explicit computations of the $\mathbb{S}^1$-action are available.
    \item     Let $\omega\in\mathrm{H}^{n-1}(\Lambda \mathbb{S}^n,\mathbb{S}^n)$ be the non-nilpotent class for the extended Goresky-Hingston product.
    Set $\mu = \lim_{m\to\infty}\tfrac{\mathrm{cr}(\omega^{\ostar  m})}{m}$.
    By the results of Section \ref{sec_extensions_of_coproduct} we see that $\mu$ agrees with $ \lim_{m\to\infty}\tfrac{\mathrm{cr}(\alpha_0^{\overline{\ostar}m})}{m}$.
    As in \cite{hingston:2013} we find that the constant $\overline{\alpha}$ in Theorem \ref{theorem_resonance_even} is $\overline{\alpha} = \tfrac{2(n-1)}{\mu}$.
    \item In \cite[Section 2.6]{cieliebak2020poincar} Cieliebak, Hingston and Oancea introduce the notion of the \textit{string point invertibility} of a closed manifold and claim that certain simply connected compact rank one symmetric spaces satisfy this condition.
    They further claim that string point invertibility implies a resonance theorem for homology classes in the free loop space $\Lambda M$ and they claim that the string point invertible compact rank one symmetric spaces satisfy a resonance theorem.
    According to \cite{cieliebak2020poincar}, the real projective spaces are not string point invertible at least with $\mathbb{Z}_2$-coefficients.
    A resonance theorem can hence not be concluded in this case.
    \end{enumerate}
\end{remark}
We now use Theorem \ref{theorem_resonance_even} to show a statement about the behavior of closed geodesics on real projective space.
As in \cite{hingston:2013} we call the number $\overline{\alpha}$ from Theorem \ref{theorem_resonance_even} the \textit{global mean frequency}.
If $\gamma$ is a closed geodesic in $\mathbb{R}P^n$ we let $\alpha_{\gamma}$ be the \textit{average index} which is defined by
$$  \alpha_{\gamma} = \lim_{n\to\infty} \frac{\mathrm{ind}(\gamma^n)}{n} .   $$
Furthermore we define the \textit{mean frequency} of $\gamma$ to be $\overline{\alpha}_{\gamma} = \tfrac{\alpha_{\gamma}}{\mathcal{L}(\gamma)}$.
The following Theorem is a version of Hingston and Rademacher's \textit{density theorem} \cite[Theorem 1.2]{hingston:2013} for even-dimensional real projective space.
\begin{theorem}\label{theorem_density}
    Let $g$ be a Riemannian metric on $\mathbb{R}P^n$ for $n$ even and let $\overline{\alpha}$ be the global mean frequency from Theorem \ref{theorem_resonance_even}.
    Let $\epsilon > 0$ and let $\mathcal{S}$ be the set of closed geodesics in $\mathbb{R}P^n$ with mean frequency in $(\overline{\alpha}-\epsilon,\overline{\alpha}+\epsilon)$.
    Then we have
    $$ \sum_{\gamma\in\mathcal{S}} \frac{1}{\alpha_{\gamma}}  \geq \frac{1}{n-1}  .    $$
\end{theorem}
\begin{proof}
    We follow the proof of \cite[Theorem 1.2]{hingston:2013}, see Section 7 of that paper.
    Let $N(d)$ be the number of distinct critical levels of homology classes in $\widehat{\mathrm{H}}_{\bullet}(\mathbb{S}^n) = \mathrm{H}_{\bullet}(\Lambda \mathbb{S}^n)\oplus \mathrm{H}_{\bullet}(P_a\mathbb{S}^n)$ for classes of degree less than or equal to $d$.
    As in \cite{hingston:2013} it is sufficient to show that $\lim_{d\to\infty}\tfrac{N(d)}{d}\geq \tfrac{1}{n-1}$.
    We claim that the classes in the sequence $A_0,A_1',A_2,A_3',\ldots$ all have distinct critical values.
    Indeed, we have that for each $i\in\mathbb{N}$ the homology class $A_i$ or $A_i'$, respectively, is dual to $\alpha_i$ or $\alpha_i'$, respectively.
    But $\alpha_{2m} = \pm {\alpha_0}^{2m} $ and $\alpha'_{2k+1} = \pm \alpha_0^{2k+1}$ with respect to the extended Goresky-Hingston product.
    Moreover, since $\alpha_0\in\mathrm{H}^{\bullet}(P_a\mathbb{S}^n)$ we have $\mathrm{cr}(\alpha_0 ) =  \mathrm{cr}(A_0) > 0$ and by Proposition \ref{prop_morse_theoretic_inequalities} we thus have
    $$   \mathrm{cr}(\alpha'_{2k+1}) \geq \mathrm{cr}(\alpha_{2k}) + \mathrm{cr}(A_0) > \mathrm{cr}(\alpha_{2k})   $$ as well as $$ \mathrm{cr}(\alpha_{2m}) \geq \mathrm{cr}(\alpha'_{2m-1}) + \mathrm{cr}(A_0) > \mathrm{cr}(\alpha'_{2m-1}) .     $$
    Hence we get an increasing sequence $$\mathrm{cr}(A_0) < \mathrm{cr}(A_1') < \mathrm{cr}(A_2) < \ldots  $$ and since taking the product with $\alpha_0$ increases the degree by $n-1$ the claim follows.    
\end{proof}

\begin{remark}
    Note that in the above proof we make use of taking the extended Goresky-Hingston product with the class $\alpha_0\in\mathrm{H}^0(P_a\mathbb{S}^n)$ to obtain an increasing sequence of critical values.
    In \cite{hingston:2013} in the case of odd-dimensional spheres Hingston and Rademacher argue with a sequence of homology classes where they use the $\mathbb{S}^1$-action to compare the critical values.
\end{remark}

We conclude this section by showing that a resonance theorem for odd-dimensional real projective space follows quite easily from the resonance theorem of Hingston and Rademacher in \cite{hingston:2013}.
\begin{lemma}\label{lemma_homotopy_equivalences_respect_energy}
    Let $n$ be odd and let $g$ be a Riemannian metric on $\mathbb{S}^n$.
    Consider the homotopy equivalence $\Phi\colon \Lambda \mathbb{S}^n\to P_a\mathbb{S}^n$ as well as its homotopy inverse $\Psi\colon P_a\mathbb{S}^n \to \Lambda \mathbb{S}^n$ as in Section \ref{sec_involution_path_space_string_operations}.
    Then there are constants $c,d>0$, depending on the metric, such that
    $$      \mathrm{cr}(  \Phi_* X) \leq \mathrm{cr}(X) + c \quad \text{and}\quad \mathrm{cr}(\Psi_* Y) \leq \mathrm{cr}(Y) + d     $$
    for all $X\in\mathrm{H}_{\bullet}(\Lambda \mathbb{S}^n), Y\in\mathrm{H}_{\bullet}(P_a\mathbb{S}^n)$.
    In particular there are constants $\eta,\mu>0 $ such that
    $$    | \mathrm{cr}(\Phi_* X) - \mathrm{cr} (X) |   \leq \eta \quad \text{and} \quad |\mathrm{cr}(\Psi_* Y) - \mathrm{cr}(Y) | \leq \mu     $$
    for all $X\in\mathrm{H}_{\bullet}(\Lambda \mathbb{S}^n)$ and $Y\in\mathrm{H}_{\bullet}(P_a\mathbb{S}^n)$.
\end{lemma}
\begin{proof}
    Recall that $\Phi\colon \Lambda \mathbb{S}^n\to P_a\mathbb{S}^n$ is given as the map
    $$  \Phi(\gamma) = \mathrm{concat}(\gamma,\eta(\gamma(0)))     $$
    where $\eta(p)\in P_a\mathbb{S}^n$ is a path from $-p$ to $p$ for $p\in\mathbb{S}^n$.
    Moreover, as in previous sections we can replace the concatenation at time $t = \tfrac{1}{2}$ by the \textit{optimal concatenation} to obtain a map $\widetilde{\Phi}\colon \Lambda \mathbb{S}^n\to P_a\mathbb{S}^n$ which is homotopic to $\Phi$.
    We then obtain that 
    $$  \mathcal{L}(\widetilde{\Phi}(\gamma)) \leq \mathcal{L}(\gamma) + \mathcal{L}(\eta(\gamma(0)))  .   $$
    Now, let $c = \max\{ \mathcal{L}(\eta(p))\,|\, p\in \mathbb{S}^n \}$, then we have that $\mathcal{L}(\widetilde{\Phi}(\gamma))\leq \mathcal{L}(\gamma) + c$.
    The inequality $\mathrm{cr}(\Phi_* X) \leq \mathrm{cr}(X) + c$ for $X\in\mathrm{H}_{\bullet}(\Lambda \mathbb{S}^n)$ follows directly. 
    The other inequalities can be shown similarly.
\end{proof}
\begin{theorem}\label{theorem_resonance_odd}
    Let $g$ be a Riemannian metric on $\mathbb{S}^n$ for $n$ odd and consider the induced energy functional on $\Lambda\mathbb{S}^n \sqcup P_a\mathbb{S}^n$.
    Consider the direct sum $\widehat{\mathrm{H}}_{\bullet}(\mathbb{S}^n) = \mathrm{H}_{\bullet}(\Lambda \mathbb{S}^n;\mathbb{Q}) \oplus \mathrm{H}_{\bullet}(P_a\mathbb{S}^n;\mathbb{Q})$.
    There are constants $\alpha,\beta > 0$ such that for all $X\in\widehat{\mathrm{H}}_{\bullet}(\mathbb{S}^n)$ we have
    $$   | \mathrm{deg}(X) - \alpha\, \mathrm{cr}(X)  | \leq \beta .      $$
\end{theorem}
\begin{proof}
    By \cite[Theorem 1.1]{hingston:2013} there are constants $\alpha_0,\beta_0$ such that 
    $$     |\mathrm{deg}(X) - \alpha_0 \mathrm{cr}(X) | \leq \beta_0     $$
    holds for all $X\in\mathrm{H}_{\bullet}(\Lambda \mathbb{S}^n)$.
    Let $Y\in\mathrm{H}_{\bullet}(P_a\mathbb{S}^n)$, then the result of Lemma \ref{lemma_homotopy_equivalences_respect_energy} and the triangle inequality show that
    $$   |\mathrm{deg}(Y) 
     - \alpha  \, \mathrm{cr}(Y) | \leq \beta_0 + \alpha_0\mu $$
    where $\mu>0$ is the constant as in Lemma \ref{lemma_homotopy_equivalences_respect_energy}.
    Hence, setting $\alpha = \alpha_0$ and $\beta = \beta_0 + \alpha\mu$ and noting that $Y$ was arbitrary yields the claim.
\end{proof}

\appendix
\section{Orientations of the completing manifolds}\label{sec_appendix_orientations}

In this appendix we explain how we choose the orientations of the completing manifolds and we prove that the diffeomorphisms $\Phi_0$ and $\Phi_1$ that are introduced in Section \ref{sec_completing_manifolds} are orientation-preserving.
We only consider the case of even-dimensional spheres.

First, we need to specify the orientations.
We fix an orientation of the standard sphere, indeed we can just take the usual orientation which says that $e_1,\ldots, e_n$ is a positively oriented basis for $T_{e_{n+1}}\mathbb{S}^n$, i.e. this means that the isomorphism $(p) \oplus  T_p\mathbb{\mathbb{S}^n}\cong \mathbb{R}^{n+1}$ is orientation-preserving.
Moreover, for each $p\in\mathbb{S}^n$ we denote by $\mathbb{S}^{n-1}_p$ the hypersphere
$$  \mathbb{S}^{n-1}_p =  \{ q\in \mathbb{S}^n\,|\, \langle q,p\rangle = 0\}.        $$
We orient it as follows.
Inside $\mathbb{S}^n$ is has a nowhere vanishing normal vector field given by the vector $-p$ and we orient $T_q\mathbb{S}^{n-1}_p$ by demanding that $(-p) \oplus T_q \mathbb{S}^{n-1}_p  \cong T_q \mathbb{S}^n$ is orientation-preserving.  
Note that the tangent space $T_q\mathbb{S}^{n-1}_p$ is
$$   T_q \mathbb{S}^{n-1}_p =  \{  Y\in\mathbb{R}^{n+1}\,|\, \langle Y,p \rangle = 0 = \langle Y, q\rangle \} .     $$

We turn to the unit tangent bundle $U\mathbb{S}^n$ which can be understood as a submanifold
$$   U\mathbb{S}^n = \{(p,v)\in \mathbb{R}^{n+1}\times \mathbb{R}^{n+1} \,|\,   \|p\|^2 = 1 = \| v\|^2 , \langle p,v\rangle = 0\} .     $$
The tangent space of $U\mathbb{S}^n$ can be described as follows.
We have
$$   T_{(p,v)} U\mathbb{S}^n =  \{ (X,Y)\in\mathbb{R}^{n+1}\times \mathbb{R}^{n+1} \,|\, \langle X,p\rangle = 0 = \langle Y,v\rangle , \,\, \langle X,v \rangle + \langle Y,p \rangle = 0\}  .   $$
We find an isomorphism
$    \varphi\colon  T_{(p,v)} U\mathbb{S}^n \to T_p \mathbb{S}^n \oplus T_v \mathbb{S}^{n-1}_p         $
by setting
$  \varphi (X,Y) =  (X,  Y + \langle X,v\rangle p )     $.
In fact, the inverse is given by $\psi\colon T_p\mathbb{S}^n\oplus  T_v \mathbb{S}^{n-1}_p\to T_{(p,v)} U\mathbb{S}^n$, $\psi(U,W) = (U,W-\langle U,v\rangle p)$.
We define the orientation of $U\mathbb{S}^n$ by declaring this isomorphism to be orientation preserving.
Similarly, for the higher fiber products $U\mathbb{S}^n\times_{\mathbb{S}^n} \ldots \times_{\mathbb{S}^n}U\mathbb{S}^n$ we have
\begin{eqnarray*}
        T_{(p,v_0,\ldots, v_l)}\Gamma_l  =  \{ (X,Y_0,\ldots, Y_l) \in (\mathbb{R}^{n+1})^{l+2}  &   | &   \langle X,p \rangle = 0 = \langle Y_i,v_i\rangle \,\,\, \text{for all}\,\,\, i\in \{0,\ldots, l\} 
            \\ & 
             &    \langle X,v_i \rangle + \langle Y_i, p \rangle = 0 \,\,\,\text{for all}\,\,\, i\in\{0,\ldots,l\} \} .
\end{eqnarray*}
Hence, we have an isomorphism 
$$  T_{(p,v_0,\ldots, v_l)} \Gamma_l \to T_p \mathbb{S}^n \oplus T_{v_0} \mathbb{S}^{n-1}_p \oplus \ldots \oplus T_{v_l}\mathbb{S}^{n-1}_p      $$
given by
$$  (X,Y_0,\ldots, Y_l) \mapsto (X,Y_0 + \langle X,v_0\rangle p , \ldots, Y_l + \langle X,v_l\rangle p )  .    $$
We define an orientation on $\Gamma_l$ by declaring this isomorphism to be orientation-preserving.

There is another way of determining an orientation of $\Gamma_l $ by looking at the embedding $\Gamma_l\hookrightarrow (U\mathbb{S}^n\times U\mathbb{S}^n)^{l+1}$.
Note that the normal bundle $N(\Gamma_l \hookrightarrow (U\mathbb{S}^n)^{l+1})$ is isomorphic to the pullback bundle $(\pi_{l+1})^* N(\Delta \mathbb{S}^n\hookrightarrow (\mathbb{S}^n)^{l+1})$.
The normal bundle 
$N(\Delta \mathbb{S}^n\hookrightarrow (\mathbb{S}^n)^{l+1})$ is given 
by
$$ N(\Delta \mathbb{S}^n\hookrightarrow (\mathbb{S}^n)^{l+1}) =  \{ (u_0,\ldots, u_l)\in (T_p\mathbb{S}^n)^{l+1} \,|\, \sum_i u_i = 0\} .      $$
There is an isomorphism
$$   N_{(p,\ldots,p)}  (\Delta^l\mathbb{S}^n\hookrightarrow (\mathbb{S}^n)^{l+1}) \cong (T_p\mathbb{S}^n)^{l}     $$
which is given by
$$  (u_0,\ldots, u_{l-1},u_l) \mapsto (u_0,\ldots, u_{l-1})           $$
and since $n$ is even this is an orientation-preserving isomorphism.
Hence, one finds that the normal space $\Gamma_l\hookrightarrow (U\mathbb{S}^n)^{l+1}$
is
$$    N_{(p,v_0,\ldots, v_l)} (\Gamma_l\hookrightarrow (U\mathbb{S}^n)^{l+1}) =  \{ (u_0,0,u_1,\ldots, u_l,0)\in (T_p U\mathbb{S}^n)^{l+1}\} .        $$

\begin{lemma}
    The orientation of $\Gamma_l$ induced as a submanifold of $(U\mathbb{S}^n)^{l+1}$ agrees with the orientation that we defined above.
\end{lemma}
\begin{proof}
    We have an orientation-preserving isomorphism
    $$
            T_{(p,v_0,\ldots, v_l)}\Gamma_l  \cong   T_p \mathbb{S}^n \oplus   T_{v_0} \mathbb{S}^{n-1}_p \oplus \ldots \oplus T_{v_l} \mathbb{S}^{n-1}_p .
    $$
    For the submanifold orientation we have an orientation-preserving isomorphism
    $$   N_{((p,v_0,\ldots, v_l)}(\Gamma_l\hookrightarrow (U\mathbb{S}^n)^{l+1}) \oplus  T_{(p,v_0,\ldots, v_l)} \Gamma_l \cong   T_{(p,v_0),\ldots, (p,v_l)} ( U\mathbb{S}^n)^{l+1} .   $$
    By definition the isomorphism
    \begin{eqnarray*}
        T_{(p,v_0),(p,v_l)} ( U\mathbb{S}^n)^{l+1} \cong T_p \mathbb{S}^n \oplus T_{v_0} \mathbb{S}^{n-1}_p \oplus \ldots \oplus T_p \mathbb{S}^n\oplus T_{v_l} \mathbb{S}^{n-1}_p 
     \end{eqnarray*}
     is orientation-preserving.
     Now, since $n$ is even we have that the swapping
    \begin{eqnarray*}  &   & 
          T_p \mathbb{S}^n \oplus T_{v_0} \mathbb{S}^{n-1}_p \oplus \ldots \oplus T_p \mathbb{S}^n\oplus T_{v_l} \mathbb{S}^{n-1}_p  \to   \\   &   &  T_p\mathbb{S}^n \oplus T_p\mathbb{S}^n\oplus \ldots \oplus T_p\mathbb{S}^n \oplus T_{v_0}\mathbb{S}^{n-1}_p \oplus T_{v_1}\mathbb{S}^{n-1}_p   \oplus \ldots \oplus T_{v_l}\mathbb{S}^{n-1}_p 
    \end{eqnarray*}
     is orientation-preserving.
    Furthermore we have that 
    $$N_{(p,v_0,\ldots,v_l)} (\Gamma_l\hookrightarrow (U\mathbb{S}^n)^{l+1}) \cong T_p (\mathbb{S}^n)^{l}$$ is orientation-preserving and again since $n$ is even we see that there is an orientation-preserving isomorphism
    \begin{eqnarray*} 
           &   &   N_{(p,v_0,\ldots, v_k)}(\Gamma_l\hookrightarrow (U\mathbb{S}^n)^{l+1}) \oplus T_{p,v_0,\ldots, v_l} \Gamma_l \cong  \\  &   &   T_p\mathbb{S}^n \oplus T_p\mathbb{S}^n\oplus \ldots \oplus T_p\mathbb{S}^n \oplus T_{v_0}\mathbb{S}^{n-1}_p \oplus T_{v_1}\mathbb{S}^{n-1}_p   \oplus \ldots \oplus T_{v_l}\mathbb{S}^{n-1}_p .
    \end{eqnarray*}
    Hence, the orientations agree.
\end{proof}
The fiber product $\Gamma_l\times_{\Delta_a}\Gamma_m$ is the pullback of $p_l\times p_m\colon \Gamma_l\times \Gamma_m \to \mathbb{S}^n\times \mathbb{S}^n$ along the map $\Delta_a\colon\mathbb{S}^n\to \mathbb{S}^n\times \mathbb{S}^n$, $p\mapsto (p,-p)$.
The normal bundle of $\Delta_a\mathbb{S}^n$ inside $\mathbb{S}^n\times \mathbb{S}^n$ is given by
$$   N_{(p,-p)}(\Delta_a\mathbb{S}^n\hookrightarrow \mathbb{S}^n\times \mathbb{S}^n) = \{(u,u) \in T_p\mathbb{S}^n \,|\, u\in T_p\mathbb{S}^n\} .    $$
The isomorphism $$ N_{(p,-p)}(\Delta_a\mathbb{S}^n\hookrightarrow \mathbb{S}^n\times \mathbb{S}^n) \to T_p\mathbb{S}^n$$ given by $(u,u)\mapsto u$ is orientation-preserving.
The normal bundle of $\Gamma_l\times_{\Delta_a}\Gamma_m$ inside $\Gamma_l\times \Gamma_m$ thus is
\begin{eqnarray*}
     &   &   N_{((p,v_0,\ldots,v_l),(-p,w_0,\ldots,w_m))}(\Gamma_l\times_{\Delta_a}\Gamma_m \hookrightarrow \Gamma_l\times \Gamma_m) = \\  &   &   \{ ((u,0,\ldots ,0),(u,0,\ldots 0))\in T_{(p,v_0,\ldots,v_l),(-p,w_0,\ldots,w_m))} \Gamma_l\times \Gamma_m  \,|\,  u\in T_p\mathbb{S}^n  \} .      
\end{eqnarray*}
Now, $\Gamma_l\times \Gamma_m$ sits inside $(U\mathbb{S}^n)^{l+1}\times (U\mathbb{S}^n)^{m+1}$, again with suitable orientation as a submanifold and $\Phi_1$ is the restriction of the orientation-preserving map 
$$   \Psi = (\mathrm{id}_{U\mathbb{S}^n})^{l+1} \times (Da)^{m+1 }\colon (U\mathbb{S}^n)^{l+1}\times (U\mathbb{S}^n)^{m+1} \to (U\mathbb{S}^n)^{l+m+2} .    $$
Note that $\Psi$ is orientation-preserving since $Da$ is orientation-preserving, see Lemma \ref{lemma_twisted_intersection_product}.
\begin{prop}
    The diffeomorphism $\Phi_1\colon \Gamma_l\times_{\Delta_a}\Gamma_m \to \Gamma_{l+m+1}$ is orientation-preserving.
\end{prop}
\begin{proof}
    It just remains to check that the isomorphism induced by $\Psi$ on the normal spaces is orientation-preserving.
    This can be checked directly by considering the explicit description of the normal spaces above.
\end{proof}
With the same strategy one can see that the diffeomorphism $\Phi_0$ is orientation-preserving.

\bibliography{lit}
 \bibliographystyle{amsalpha}
 
\end{document}